\documentclass[preprint,12pt,authoryear]{elsarticle}

\usepackage{amsfonts}
\usepackage{amsmath}
\usepackage{amssymb}
\usepackage{amsthm}
\usepackage{graphicx}
\usepackage{xcolor}
\usepackage{epsfig}
\usepackage{url}
\usepackage{bm}
\usepackage{subcaption}
\usepackage{natbib}
\usepackage{multirow}
\usepackage{setspace}
\usepackage{etoolbox}
\usepackage{graphicx}
\usepackage{blindtext}
\usepackage{enumitem}
\usepackage{geometry}
\definecolor{blue2}{rgb}{0.0,0.4,0.8}
\usepackage[ruled,vlined,linesnumbered]{algorithm2e}

\SetCommentSty{mycommfont}
\SetNoFillComment
\allowdisplaybreaks
\usepackage{booktabs}
\usepackage[flushleft]{threeparttable}
\usepackage{longtable}
\usepackage{makecell}
\usepackage{float}
\usepackage{tikz}

\usepackage[bookmarks=false,colorlinks]{hyperref}
\AtBeginDocument{%
	\hypersetup{
		linkcolor=blue,   
		citecolor=red,
		urlcolor=magenta}}

\restylefloat{table}
\usepackage[export]{adjustbox}
\usetikzlibrary{shapes.geometric, arrows}
\tikzstyle{startstop} = [rectangle, rounded corners, minimum width=3cm, minimum height=1cm,text centered, draw=black, fill=white!30]
\tikzstyle{io} = [trapezium, trapezium left angle=70, trapezium right angle=110, minimum width=3cm, minimum height=1cm, text centered, draw=black, fill=white!30]
\tikzstyle{process} = [rectangle, minimum width=3cm, minimum height=1cm, text centered, text width=5cm, draw=black, fill=white!30]
\tikzstyle{decision} = [diamond, aspect=3, minimum width=3cm, minimum height=1cm, text centered,text width=5cm, draw=black, fill=white!30]
\tikzstyle{arrow} = [thick,->,>=stealth]

\usepackage{caption,setspace}
\newtheorem{modl}{Model}

\newenvironment{model}{\begin{samepage}\begin{modl}}{\end{modl}\end{samepage}}

\newcommand{\N}{\mathcal{N}}
\renewcommand{\P}{\mathcal{P}}
\newcommand{\D}{\mathcal{D}}
\newcommand{\K}{\mathcal{K}}
\newcommand{\M}{\mathcal{M}}
\newcommand{\IB}{\mathcal{IB}}
\newcommand{\OB}{\mathcal{OB}}
\newcommand{\Z}{\mathcal{Z}}
\newcommand{\PM}{\mathcal{P}_m}
\newcommand{\NM}{\mathcal{N}_m}
\newcommand{\PZP}{\mathcal{P}_{z_p}}
\newcommand{\PZD}{\mathcal{P}_{z_d}}
\newcommand{\V}{\mathcal{V}}
\newcommand{\U}{\mathcal{U}}
\newcommand{\A}{\mathcal{A}}
\newcommand{\G}{\mathcal{G}}
\newcommand{\add}{$\mathsf{ADD}$\xspace}
\newcommand{\rem}{$\mathsf{REM}$\xspace}
\newcommand{\swap}{$\mathsf{SWAP}$\xspace}
\newcommand{\ra}{$\mathsf{RA}$\xspace}
\newcommand{\ri}{$\mathsf{RI}$\xspace}

\begin{document}

\begin{frontmatter}

\title{A chance-constrained dial-a-ride problem with utility-maximizing demand and multiple pricing structures}

\author[1]{Xiaotong Dong}
\author[2]{Joseph Y. J. Chow}
\author[1]{S. Travis Waller}
\author[1]{David Rey\corref{mycorrespondingauthor}}
\cortext[mycorrespondingauthor]{Corresponding author}

\address[1]{Research Centre for Integrated Transport Innovation (rCITI)\\
	School of Civil and Environmental Engineering\\
	University of New South Wales, Sydney, NSW 2052, Australia}
\address[2]{C2SMART University Transportation Center,\\ New York University Tandon School of Engineering, Brooklyn, NY, USA}

\begin{abstract}
The classic Dial-A-Ride Problem (DARP) aims at designing the minimum-cost routing that accommodates a set of user requests under constraints at an operations planning level, where users' preferences and revenue management are often overlooked. In this paper, we present a mechanism for accepting/rejecting user requests in a Demand Responsive Transportation (DRT) context based on the representative utilities of alternative transportation modes. We consider utility-maximizing users and propose a mixed-integer programming formulation for a Chance Constrained DARP (CC-DARP), that captures users' preferences in the long run via a Logit model. We further introduce class-based user groups and consider various pricing structures for DRT services. A customised local search based heuristic is developed to solve the proposed CC-DARP. We report numerical results for both DARP benchmarking instances and a realistic case study based on New York City yellow taxi trip data. Computational experiments performed on 105 benchmarking instances with up to 96 nodes yield an average optimality gap of $2.69\%$ using the proposed local search heuristic. The results obtained on the realistic case study reveal that a zonal fare structure is the best strategy in terms of optimising revenue and ridership. The proposed CC-DARP formulation provides a new decision-support tool to inform on revenue and fleet management for DRT systems at a strategic planning level. 
\end{abstract}


\end{frontmatter}


\section{Introduction}
\label{S:1}

The Dial-a-Ride Problem (DARP) is a highly-constrained combinatorial optimization problem initially designed for providing door-to-door transportation for people with limited mobility (the elderly or disabled). It consists of routing and scheduling a fleet of capacitated vehicles to service a set of requests with specified pickup and drop-off locations as well as a tight time window imposed on ride time. With the details of requests obtained either beforehand (static DARP) or en-route (dynamic DARP), dial-a-ride service providers strive to deliver an efficient and yet high quality transport solution under a set of constraints. 

What makes DARP distinct from general vehicle routing problems is that service quality must be considered, either as constraints or as a part of the objective function. To balance the trade-off between users' inconvenience and system routing cost is the biggest challenge dial-a-ride service operators are facing when designing an efficient solution. Generally, a DARP request can be categorised as either an inbound trip (with a tight time window specified on arrival time) or an outbound trip (with a tight time window specified on departure time) \citep{cordeau2003tabu}. Additionally, a maximum ride time is imposed as a constraint for each user.

The standard formulation of DARP \citep{cordeau2006branch} has been extensively studied over decades. Research on DARP tends to incorporate the problem with various real-life characteristics to extend its applicability on practical matters (e.g. rich vehicle routing problems). One of the principals of classic (static) DARP is to serve all requests, and to optimise the routing and scheduling of vehicles in order to service all users. Sometimes, in order to serve one or two requests, significant adjustments in routing must be made, without necessarily benefiting either system routing cost or the overall service quality of users. Demand Responsive Transportation (DRT) providers/operators are often ridership-dependent. Persistently serving all requests at the cost of compromising users' experiences of on-board passengers can be rather myopic and unsustainable, considering service quality is the major incentive for stable future ridership \citep{yooneffect}. Therefore, strategic planning, including designing fare policies for a fleet that operates as a DARP, is also critical (albeit overlooked) in DRT services. 

Furthermore, existing studies on DARP fall short of taking users' preferences into consideration. Since most DARP models impose identical time window and ride time constraints to all customers, users' preferences including their tolerances towards travel time, schedule delay, fare, etc. are left unaccounted. Most DRT operators using classic DARP models assume that user demand is inelastic even when it is not the case \citep{sayarshad2015scalable}. Consequences including risks of no-show and/or cancellation after booking the service (as customers might not be happy with the trip offer) may lead to poor service quality and reduce ridership in the long-run.

We propose a DARP model integrated with users' preferences from DRT operators' perspective, facilitating them to optimise their overall profit while maintaining service quality.  In this model, we consider users' utility and take into account users' preferences within a dial-a-ride problem. We propose a Chance-Constrained Programming (CCP) \citep{charnes1959chance} approach, which aims to identify users that are better off travelling by DRT services than using a reserve travel mode. The proposed CCP formulation incorporates operator's accept/reject decisions of users' requests. Only accepted users are served by the DRT operator, while the rest of the requests are rejected and left out of consideration of routing and scheduling. This CCP accept/reject mechanism is expected to guide service providers on a strategic planning level to avoid blindly accepting all requests at the cost of driving down system profit or the user experiences of on-board passengers.

This study also explores the design of revenue/fleet management and pricing differentiation \citep{talluri2004revenue} under the proposed chance-constrained DARP model. Ridership response studies often define price differentiation policies targeting different segments (classes) based on user groups and trip characteristics \citep{cervero1990transit}. Users from different segments (classes) of age, income, auto access and trip purpose tend to respond with various sensitivities to fare changes. Accordingly, we introduce multiple user classes in the chance-constrained DARP model. The proposed multi-class, chance-constrained DARP is formulated as a Mixed-Integer Linear Program (MILP) and a customized heuristic algorithm based on local search which iteratively explores neighbourhoods of both the selection of requests and routing solutions is developed. By embedding class-based formulation into various fare structures, the proposed model also helps design and evaluate pricing strategies and revenue management. Numerical results based on Cordeau's benchmarking DARP instances \citep{cordeau2006branch} are reported to show the strong performance of the customised heuristic, followed by realistic DRT scenarios generated from the Yellow taxi trip data of New York City ({\small \url{ https://data.cityofnewyork.us/Transportation/2016-Yellow-Taxi-Trip-Data/k67s-dv2t}}).



The remainder of this paper is organised as follows: in Section \ref{lit}, we briefly review the existing research on DARP models and solution methodologies and list out our contributions. In Section \ref{math}, we introduce the mathematical formulation of the proposed multi-class chance-constrained DARP model. A customized local search heuristic is presented in Section \ref{algo}. Numerical experiments and results are reported in Section \ref{num}. Our work is concluded with a discussion of findings and future research directions in Section \ref{con}.
 
\section{Literature Review}
\label{lit}
\subsection{DARP and VRP variants}
 
A Vehicle Routing Problem (VRP) concerns designing a minimized cost routing for a fleet of vehicles that depart from the same depot so that they can visit a set of customers with known demand \citep{dumas1991pickup}. One of the generalizations of VRP is the pickup and delivery problem with time windows (PDPTW), where each request is constrained with an earliest and latest service time. As a special case of PDPTW, the Dial-a-Ride Problem (DARP) looks into transportation of people instead of goods. Hence, requirements related to service quality are satisfied either through additional constraints of the model or with extra terms in the objective function. \\

The factor that makes DARP challenging is the trade-off between operational cost and user inconvenience \citep{cordeau2003tabu} \citep{cordeau2006branch}. The most frequently used objective function in the DARP is to minimize the operational cost, subject to full demand satisfaction and other service design constraints \citep{molenbruch2017typology}. However, it is sometimes unrealistic to satisfy all demand, as the operational cost would be too large, and service level would drop significantly because of the detours required to serve all user requests. \\

Although studies over the years have developed a series of exact and approximate solution method applied to the DARP, solving the DARP and its variants remains challenging, given that it is NP-hard. For exact solution methods, \cite{cordeau2006branch} proposed a branch-and-cut algorithm, which was further strengthened by other studies via adding valid inequalities \citep{ropke2007models} \citep{parragh2011introducing} \citep{hall2009integrated}. The column generation technique, embedded in a branch-and-bound scheme as a branch-and-price method, is also widely applied to solve larger instances \citep{parragh2012models} \citep{garaix2011optimization}. Approximate solution methods, on the other hand, are favoured by many researchers since they can generate close to optimal solutions with a cheaper computational cost. Most heuristics for the DARP are based on \cite{jaw1986heuristic}'s time-sequential insertion heuristic \citep{diana2004new} \citep{toth1997exact}. To better escape local optima, more advanced heuristics including metaheuristics (e.g. Tabu Search) and matheuristics for DARP have also been studied by many \citep{cordeau2003tabu} \citep{ho2011local} \citep{archetti2006tabu}. \\

 Alternatively, a variant of TSP assigns profit values to links and/or nodes in the network but relaxes the constraint of visiting all nodes. \cite{ feillet2005traveling}'s work is among the first ones that links TSPs with Profits. This generalisation of TSP differs from classic TSP as each customer has an associated profit value that can be gained if this user is visited. In \cite{feillet2005traveling}’s study, the two conflicting objectives of TSPs with profits, namely to collect more profit and to minimise travel cost, are formulated either as a weighted sum, or one of the two objectives is constrained with a bound value. In the first case, the problem is also identified as the Profitable Tour Problem (PTP) \citep{ dell1995prize}, which occurs more often as a subproblem in a column generation based solution algorithm in various routing problems. The problems of the second case are also known as the Orienteering Problem (OP) or the Prize-Collecting TSP (PCTSP, also known as quota TSP), depending on which one of the objectives is stated as a constraint. \cite{archetti2009capacitated} explored two versions of capacitated VRPs with profits: an extension of the capacitated OP with multiple tours (Capacitated Team Orienteering Problem, shortened to CTOP), as well as the capacitated PTP (CPTP) with a single vehicle. A column generation based exact approach (branch-and-price) and various metaheuristics are implemented to solve the problems. A branch-and-cut algorithm for CPTP is later developed by \citep{jepsen2014branch} with a number of valid inequalities introduced. \citep{chow2012generalized} introduced a new class of problems, named the generalised PTPs (GPTPs), in order to cater to the need of location-based activity routing systems. In GPTPs, nodes are categorised into different clusters for different activity types. The objective is to maximize the summed utility of all selected activities and minimize the disutility of travelling, with a set of constraints including time windows and fuel constraints, etc. In a separate paper, \cite{feillet2005profitable} also introduced a similar routing problem called the profitable arc tour problem (PATP), where the profit is associated with arcs, instead of vertices.  \\

Although profit-driven TSPs and VRPs have been studied by many, DARP with profit has hardly ever been explored. \cite{schonberger2003combined} first link profit and the Pickup and Delivery Problem with Time Windows (PDPTW) by proposing the Pickup and Delivery Selection Problem (PDSP), which is an extension of the PDPTW with an accept decision of the requests. With a revenue associated to each request, service providers aim to maximise the total profit by identifying a selection of requests to service. 
\cite{parragh2015dial} introduced the concept of DARP with Split Requests and Profits (DARPSRP). A path-based formulation of the DARPSRP is proposed, where a revenue is associated to each request and the objective is to maximise the total profit. Unlike classic DARPs, requests can be rejected by DRT operators if they are not profitable. A branch-and-price algorithm for the DARPSRP is developed alongside a Variable Neighbourhood Search (VNS) metaheuristic method. To the best of our knowledge, the existing studies on Dynamic DARP request accept/reject decisions are simply made simultaneously as routing is being optimised. No exploration has been done aiming to optimise the selection of requests as an independent procedure on account of users' utility and preferences, mainly because it can be tricky to observe and determine utility in a dynamic setting. However, in our case, utility theory can be well applied in strategic planning and revenue management.

\subsection{Choice Models and VRPs}

Service attributes including in-vehicle travel time, waiting time, price, safety and reliability influence travellers when deciding whether to take a trip by public transport. When modelling and analysing mode choice decisions, multinomial logit models have been widely applied for estimating mode split \citep{train2009discrete}. Utility functions are estimated in logit models to quantify users' satisfaction level in each travel mode. Basic utility functions of public transport travel mode take into consideration of travel time, waiting time, price and headway between two consecutive services. Follow-up studies introduce additional characteristics that might also affect users’ choices including density of standees per square meter, availability of information panels, proportion of seats been used, etc \citep{tirachini2014multimodal} \citep{dell2011quality}. Parameters of utility functions, including value of travel time (VOT), value of waiting time and schedule delay and their correlation can vary with different time of the day, demographic characteristics of an area, etc \citep{athira2016estimation} \citep{lisco1968value} \citep{mohring1987values} \citep{small1982scheduling}. \\

Transit pricing research has received a wide attention for decades, with objectives such as generating more revenue and encouraging users to switch from private transport modes to public transportation. \cite{cervero1990transit} examined how the changes in pricing policy, including fare levels, fare structures and forms of payment, affect ridership and revenue. Ridership response studies are often broken down into segments by user groups and trip characteristics. Users from different segments of age, income, auto access, trip purpose and trip lengths tend to respond with various sensitivities to fare changes. Furthermore, the operating environment such as land use (density, composition, etc) and location (CBD, suburb, etc) also play a role on fare elasticities. For example, \cite{paulley2006demand} demonstrate in their work that people in areas with low population densities rely more on cars and less on public transport, and therefore are more likely to switch to private travel when fares increase. Fare elasticity values of off-peak time periods are about twice the peak values. Also, travellers with higher incomes tend to have higher elasticity. Another aspect of public transport pricing scheme is fare structure. As a component within public transport systems, fare structure refers to the spatial structure reflecting the relationship between fare levels and distance travelled \citep{streeting2006developments}. There are four commonly used fare structures in public transport system: flat fare, distance-based fare, time-based fare and zone-based fare \citep{chin2016nonadditive}. A well-suited fare structure provides a balance of efficiency, ridership and equity benefits. Reforming fare structure can be a less costly way of boosting ridership in the short term \citep{liu2019evaluating}.\\

Fares that maximize the overall social welfare can be rather low and might require a subsidy that covers the operating cost \citep{ tirachini2014multimodal}, due to the need of easing congestion, promoting social equity and environmental considerations. A low fare can attract users of other travel modes to switch to public transport, but at the cost of increasing travel time of all users, more stopping delay and more discomfort due to overcrowding services. \cite{cervero1990transit} point out that transit users are approximately twice more sensitive in service quality changes than in fare modifications, which implies that public transport services would be able to compete with private travel mode better when service quality is improved significantly, rather than subsidizing fares. On the other hand, from an economic perspective, optimal pricing policies that maximize social welfare as a whole vary based on the economic logic behind (first-best pricing, second best pricing, etc) \citep{tirachini2012multimodal} \citep{gentile2005advanced}. \\
 
However, quality of service can be difficult to define since the measurable and tangible elements only represent a small part of the service, whereas most elements are intangible and heterogeneous. Knowing exactly what improvements need to be made depend largely on the characteristics and particularities of the targeted area and individuals \citep{redman2013quality}. Customer based quality is firstly mentioned by \cite{falcocchio1979methodology} with three main dimensions: convenience, comfort and safety. This list was further extended to eight dimensions with 41 attributes. According to the survey conducted by \cite{gunay2016estimation}, the reasons why DRT services are not favoured among travellers are chiefly proximity (25\%) and fare (23\%). The modal shift from automobiles to DRT services is estimated to be 52\% once an intelligent transport service is introduced which provides flexible and reliable services.
 
 \subsection{Our Contributions}
As outlined in our review of the literature, existing formulations for DARP and its variants focus on operational planning, and largely overlook revenue management considerations. In this paper, we approach revenue and fleet management in DRT services via a dial-a-ride problem formulation that accounts for the preferences of utility-maximizing users in the long run. We propose a multi-class Chance Constrained DARP (CC-DARP) model which equips DRT operators with an accept/reject mechanism for user selection. The proposed class-base feature of the model takes into account users socio-demographic attributes, which is unique for analysing various fare structures and revenue management at a strategic level. We present a MILP formulation for the proposed CC-DARP model and develop a customized local search based heuristic to solve the CC-DARP. The proposed solution method combines local search moves on both user selection and route optimization. The proposed local search based heuristic yields an average optimality gap of $2.69\%$ in competitive time across 105 benchmarking instances of the DARP literature, when compared to the best integer solution obtained via the exact MILP approach. This highlights the ability of the proposed approach for the CC-DARP to compromise solution quality and computing time. We also implement the proposed solution method on a realistic case study synthesized from the yellow taxi trip data of New York City (NYC) which provides insights into revenue management strategies for DRT operators. 

\section{Mathematical Formulation}
\label{math}

We follow the formulation of \cite{cordeau2006branch} for the base classic DARP model and briefly recall the main elements of this formulation hereafter. Let $n$ be the number of requests. Our problem is defined on a complete directed graph $\G=(\N,\A)$, where $\N= \P \cup \D \cup \{0, 2n+1\}$. $\A$ is the set of arcs $(i,j)$, where $i \in \N$ and $j \in \N$. Subset $\P =\{1,\dotsc,n\}$ contains pick-up nodes and subset $\D=\{n+1,\dotsc,2n\}$ contains drop-off nodes, while nodes 0 and $2n+1$ represent the origin and destination depots respectively. Each node pair $(i,n+i)$ represents a travel request from origin node $i$ to destination node $n+i$.  

We assume that $\K$ is the set of vehicles that are providing DRT services in this network. Each vehicle $k \in \K$ has a capacity $Q_k$ and its maximal vehicle travel time is $T_k$. Each request $i$ has an associated load $q_i$ (number of passengers) and a non-negative service duration $d_i$. We assume that $q_0 = q_{2n+1} = d_0 = d_{2n+1} =0$, and $q_i = - q_{n+i}$. The earliest and latest time that service may begin at node $i$ is the time window $[e_i, l_i]$. For each arc $(i,j) \in \A$ in the network, $c_{ij}$ represents its travel cost and $t_{ij}$ represents its direct travel time. We denote $L$ as the maximum acceptable travel time. 

For each arc $(i,j) \in \A$ and each vehicle $k \in \K$, let $x_{ij}^{k}$ be a binary decision variable equal to 1 if the shared mobility service vehicle $k$ is used from node $i$ to node $j$. For each node $i \in \N$, let $B_i^k$ be the time at which node $i$ is served by the DRT system and let $Q_i^k$ be the load of vehicle $k$ after visiting node $i$. We denote $L_i^k$ as the travel time of user $i$ using the DRT service if assigned to vehicle $k$. Note that for homogeneous vehicles, one can reduce the number of constraints and variables by using aggregated time and load variables $B_i$ and $Q_i$ at every node except for the depot nodes $\{0,2n+1\}$ \citep{cordeau2006branch}. We focus on homogeneous vehicles with identical capacities. A summary of the mathematical notations used throughout the paper is provided in Table \ref{notation}.

\begin{table}
\begin{center}
\begin{tabular}{ ll } 
\toprule
\multicolumn{2}{l}{\textbf{Classic DARP}} \\
$x_{ij}^k$ & Binary variable equal to 1 if vehicle $k$ travels from node $i$ to $j$\\  
$B_i^k$ & The time at which vehicle $k$ begins/finishes services at depot\\ 
$B_i$ & The time at which vehicle begins services user $i$\\
$[e_i,l_i]$ &The time window within which the service may begin at node $i$\\
$L_i$&The ride time of user $i$\\
$t_{ij}$ & The travel time between node $i$ and $j$\\
$c_{ij}$ & The routing cost between node $i$ and $j$ \\
$q_i/q_{im}$ & The load of request $i$ (in class $m$) \\
\midrule
\multicolumn{2}{l}{\textbf{Choice model}} \\
$V_i$ & Representative utility value of DRT service for user $i$;\\
$\hat{V_i}$ & Representative utility value of 'private travel' mode for user $i$;\\
$\beta_T, \beta_S$ & Parameters that convert travel time and schedule delay into monetary units\\
$f_i/f_{im}$ & Fare charged to user $i$ (in class $m$) according to a specified fare structure\\
$\hat{c_i}$ &  Travel cost of private travel mode for user $i$ \\
$\beta_F$ & Coefficient for travel cost\\
$y_i$ & Binary variable equal to 1 if user $i$ is served and 0 otherwise\\
\midrule
\multicolumn{2}{l}{\textbf{Sets}} \\
$\N$ & Set of nodes, $\N = \P \cup \D \cup \{0,2n+1\}$ \\
$\A$ & Set of arcs $(i,j)$\\
$\P$ & Set of pick-up nodes or set of requests, $\P=\IB \cup \OB$ \\
$\D$ & Set of drop-off nodes \\
$\IB$ & Set of origins of inbound trips\\
$\OB$ & Set of origins of outbound trips \\
$\K$ & Set of vehicles \\
$\M$ & Set of classes \\
$\Z$ & Set of geographic zones\\
\bottomrule
\end{tabular}
\end{center}
\caption{Notations}
\label{notation}
\end{table}

We next present the utility functions and chance constraint formulation in Section \ref{util}. We introduce a multi-class formulation of users' requests in Section \ref{class} before summarizing the resulting mixed-integer linear formulation in Section \ref{ccdarp}.

\subsection{Utility Functions and Chance Constraint}
\label{util}

For simplicity purposes, we consider two travel modes available for users to choose between a DRT service and the best option available to users other than DRT. In this discussion, we assume the other travel mode to be a private travel option with its utility function defined accordingly. The representative utility for each travel mode are defined as functions of travel time, schedule delay and monetary cost \citep{dong2020dial}.

For a request $i \in \P$ (note that $\P$ is used as both the set of pick-up nodes and the set of requests), the representative utility for private travel is:
\begin{equation}\label{uhat}
\hat{V_i}= - \beta_T t_{i,n+i} - \beta_F \hat{c_i}, \qquad\forall i \in \P.
\end{equation} 

For the DRT service, as in traditional DARP models, we distinguish between inbound and outbound user requests, denoted by sets $\IB$ and $\OB$, respectively. Each user mark their trip as either an inbound (with a tight time window on departure) or an outbound (with a tight time window on arrival) trip when they send out their request. Note that $\P = \IB \cup \OB$ and $\IB \cap \OB = \emptyset$. The  representative utility function for the DRT service are: 
\begin{align}
&V_i= - \beta_F F - \beta_T L_i  - \beta_S (B_i  - e_i), \qquad&&\forall i \in \IB, \label{uib}\\
&V_i=-\beta_F F-\beta_T L_i -\beta_S (l_{n+i} -B_{n+i}), \qquad&&\forall i \in \OB. \label{uob}
\end{align}

In utility functions (\ref{uhat})-(\ref{uob}), $\beta_T$ and $\beta_S$ are parameters that convert travel time and schedule delay into monetary units. For private travel, we assume that travel time $t_{i,n+i}$ is fixed and that schedule delay is null. For the DRT service, travel time (ride time) is represented by variable $L_i$, which depends on the route serving user $i$. Schedule delay for DRT is defined based on trip types. For inbound trips (trips with a tight time window imposed on departure time), schedule delay is defined as waiting time at pickup locations $B_i-e_i$. Whereas for outbound trips, schedule delay is defined as the amount of time that user $i$ arrives earlier than expected $l_{n+i}-B_{n+i}$, assuming the upper bound of time window $l_{n+i}$ is the preferred arrival time. Let $\beta_F$ be the coefficient for monetary travel cost, which is the fare $f_{i}$ for DRT service and cost $\hat{c_i}$ (including gas, maintenance, registration, etc.) for private travel mode.\\

In this work, we assume that the difference in utility $\Delta{U_i}=(\hat{V_i}-V_i)+\epsilon$ is a random variable following a logistic distribution $\Delta{U_i}\sim(\overline{\Delta{U_i}},s)$; wherein $\overline{\Delta{U_i}}$ is the deterministic utility gap defined as follows:
\begin{align}
\overline{\Delta{U_i}}= \hat{V_i}-V_i = \begin{cases}
\beta_T(L_i-t_{i,n+i})  + \beta_S ( B_i-  e_i)+ \beta_F (F-\hat{c_i}),&\text{if } i \in \IB,\\
\beta_T(L_i-t_{i,n+i})  + \beta_S ( l_{n+i}-  B_{n+i})+ \beta_F (F-\hat{c_i}),&\text{if } i \in \OB, \label{deltau}
\end{cases}
\end{align}

and $s$ is the scale parameter proportional to its standard deviation. Observe that $\overline{\Delta{U_i}}$ also serves as a decision variable in the proposed CC-DARP, since it is a function of variables $L_i$ and $B_i$. \\

The sign of $\Delta{U_i}$ indicates the preference of user $i$: if $\Delta{U_i} \leq 0$, then in the long-run, user $i$ is expected to experience a higher utility using DRT than the private travel alternative. To capture users' preferences in the long-run, we introduce the following chance constraint to ensure that $\Delta{U_i}\leq0$ , $\forall i \in \P$ holds with a confidence level of $p$:
\begin{equation}
\mathrm{Pr}(\Delta{U_i}\leq 0) \geq p,\qquad \forall i \in \P.\label{ncc1}
\end{equation}

Let $F_{\Delta{U_i}}$ be the cumulative distribution function of random variable$\Delta{U_i}$:
\begin{equation}
\mathrm{Pr}(\Delta{U_i}\leq0) = F_{\Delta{U_i}}(0),\qquad \forall i \in \P.\label{pdf}
\end{equation}

Combined with (\ref{pdf}), for each $i \in \P$, chance constraint (\ref{ncc1}) can be written as:
\begin{align}
 F_{\Delta{U_i}}(0) &\geq p,\nonumber\\
\Leftrightarrow \qquad 0 &\geq Q_{\Delta{U_i}}(p),\nonumber\\
\Leftrightarrow \qquad 0 &\geq\overline{\Delta{U_i}}+s\ln\left(\frac{p}{1-p}\right), \label{nncc}
\end{align}

where $ Q_{\Delta{U_i}}(p)$ is the inverse cumulative distribution function (quantile function) of the logistic distribution $\Delta{U_i}\sim(\overline{\Delta{U_i}},s)$. Combining (\ref{nncc}) and utility functions (\ref{deltau}):
\begin{align}
&\beta_T(L_i-t_{i,n+i})  + \beta_S (B_i-  e_i) + \beta_F (F-\hat{c_i}) \leq -s\ln\left(\frac{p}{1-p}\right), &&\forall i \in \IB, \label{cc1}\\
&\beta_T(L_i-t_{i,n+i})  + \beta_S (l_{n+i}-  B_{n+i}) + \beta_F (F-\hat{c_i})\leq -s\ln\left(\frac{p}{1-p}\right), &&\forall i \in \OB.\label{cc2}
\end{align} 		

The goal of the proposed chance-constrained formulation is to ensure that the DRT service operator only serves users' which are better-off using DRT services than their private travel alternative with a confidence level of $p$. To incorporate the proposed chance constraints (\ref{ncc11}) and (\ref{ncc22}) in the formulation, we introduce a binary variable $y_i \in \{0,1\}$ representing if the DRT service operator's chooses to serve user request $i \in \P$ (1) or not (0). Let $W_i$ be a large constant (the value of this constant is discussed later). Using variable $y_i$, we rewrite Eqs. (\ref{ncc11}) and (\ref{ncc22}) as on/off constraints as follows:
\begin{align}
&\beta_T(L_i-t_{i,n+i})  + \beta_S (B_i - e_i) + \beta_F (F-\hat{c_i})-(1-y_i)W_i \leq -s\ln\left(\frac{p}{1-p}\right), &&\forall i \in \IB, \label{ncc11}\\
&\beta_T(L_i-t_{i,n+i})  + \beta_S (l_{n+i} - B_{n+i}) + \beta_F (F-\hat{c_i})-(1-y_i)W_i \leq -s\ln\left(\frac{p}{1-p}\right), &&\forall i \in \OB.\label{ncc22}
\end{align}

To link user requests accept/reject decisions with routing decision variables $x_{ij}^k$, we rewrite the classical DARP constraint requesting that all users are served as:
\begin{equation}\label{eq:link}
\sum_{j \in \N} \sum_{k \in \K} x_{ij}^k = y_i, \qquad \forall i \in \P.
\end{equation}

The proposed chance-constrained formulation links user requests accept/decisions to the random variable $\Delta{U_i}$ representing the relative utility of user $i$. This mechanism is sensitive to the required confidence level $p$ and the scale parameter $s$ of the assumed logistic distribution. The impact of these parameters on model performance are examined via sensitivity analyses in our numerical experiments (see Section \ref{saps}).

\subsection{Class-based Fare Structures}
\label{class}

In the discussion of Section \ref{util}, we assumed that all users are homogeneous, and that a flat fare $F$ is charged for the DRT service. In this section, we group users into various classes by their socio-demographic characteristics (age, income, car ownership, etc) and introduce a class-based pricing system with various potential fare structures for DRT services, including distance-based fare and zone-based fare. \\

Let $\M=\{1,2,...,M\}$ be a set of classes. Each request $i$ can be categorised into one of the mutually exclusive class of requests $\P_m$ (i.e. $\cup_{m \in \M} \P_m = \P$). Furthermore, we redefine utility functions on a class basis with class-specific parameters $\beta_{i}^m,\beta_{F}^m,\beta_{T}^m,\beta_{S}^m$ and user-based fare $f_{im}$. Specifically, we rewrite the deterministic part of the utility functions \eqref{uhat}-\eqref{uob} as follows. \\

For private travel:
\begin{equation}\label{uhat2}
\hat{V_i}=- \beta_T^m t_{i,n+i} - \beta_F^m \hat{c_i}, \qquad \forall m \in \M, \forall i \in \PM.
\end{equation} 

For DRT service: 
\begin{align}
&V_i=- \beta_F^m f_{im} - \beta_T^m L_i  - \beta_S^m (B_i  - e_i), \quad&&\forall m \in \M, \forall i \in \IB \cap \PM, \label{cbuib2}\\
&V_i=-\beta_F^m f_{im}-\beta_T^m L_i -\beta_S^m (l_{n+i} -B_{n+i}), \quad&&\forall m \in \M,\forall i \in \OB\cap \PM.\label{cbuob2}
\end{align}

Accordingly, we assume that the scale parameter of the random variable is also class-dependent and the random variable is rewritten as $\Delta{U_i}\sim(\overline{\Delta{U_i}},s_m)$, for each request $i \in \P_m$. Recall that $\overline{\Delta{U_i}} = \hat{V}_i - V_i$. Observe that $\hat{V}_i$ is a constant, and that $V_i$ is a real bounded variable which upper and lower bounds, denoted $\overline{V}_i$ and $\underline{V}_i$ respectively, can be determined by taking the maximum and minimum values of the right-hand sides of Eqs. \eqref{cbuib2}-\eqref{cbuob2}. Details on the calculation of these bounds can be found in \cite{dong2020dial}. Let $W_{im} = \hat{V}_i - \underline{V}_i +s_m\ln\left(\frac{p_m}{1-p_m}\right)$. The chance constraints \eqref{ncc11} and \eqref{ncc22} are redefined with a specified confidence level $p_m$ assigned to each class $m$  and written compactly as: 
\begin{equation}
0\geq \hat{V}_i - V_i +s_m\ln\left(\frac{p_m}{1-p_m}\right)-(1-y_i)W_{im}, \qquad \forall m \in\M , \forall i \in \PM.\label{cbnccc}
\end{equation}

We next present two fare structures that can be used to determine user-based trip fares $f_{im}$ in DRT services. 

\subsubsection{Distance-based Fare Structure}
One of the most common fare structures for public transport is a distance-based fare structure, as widely applied in taxi services. Under distance-based fare structure, customised fares are charged to passengers based on the distances of their trips. Integrating user classes into this fare structure, we define the fare $f_i^{db}$ charged to user $i$ as follows:
\begin{equation}
f_i^{db} = \alpha_m Dist_{i,n+i}, \qquad \forall m \in \M, \forall i \in \PM.\label{dbfare}
\end{equation}

In (\ref{dbfare}), $Dist_{i,n+i}$ represents the direct travel distance between the pickup and drop off locations of request $i$, while $\alpha_m$ is a class-based cost parameter that converts distances into monetary units.

\subsubsection{Zone-based Fare Structure}
Zone-based or zonal-based fare structures divide the service area into several districts (zones), and the fare for request $i$ is calculated according to which zones the pickup location $i \in \P$ and drop-off location $(i+n) \in \D$ lie in. Let $\Z=\{1,2,...,Z\}$ be a set of geographical zones. The pickup and drop off locations of request $i$ can be each categorised into one of the mutually exclusive zone $\P_z$ (i.e. $\cup_{z \in \Z} \P_z = \P$). Assume that for request $i$, its pickup location $i$ lies in zone $\P_{z_p}$ (i.e. $i \in \P_{z_p}$), and its drop-off location $i+n$ belongs to zone $\P_{z_d}$ (i.e. $(i+n) \in \P_{z_d}$). The zone-based fare value of user request $i$, $f_i^{zb}$, is then calculated as follows:
\begin{equation}
f_i^{zb} = \theta_{{z_p},{z_d}} f_m, \qquad \forall m \in \M , \forall z_p, z_d \in \Z, \forall i \in \PM \cap \PZP, \forall (i+n) \in \PZD .\label{zbfare}\\
\end{equation}

In (\ref{zbfare}), $f_m$ is a base fare for user class $m$, whereas $\theta_{{z_p},{z_d}}$ is a weight parameter considering the corresponding zones $z_p$ and $z_d$ that include pickup location $i$ and drop up location ($n+i$) of request $i$.

\subsection{Chance-constrained DARP (CC-DARP) Formulation}
\label{ccdarp}

We incorporate the proposed class-based, chance-constrained DARP (CC-DARP) formulation with the classic DARP model. To account for user request accept/reject decisions, we propose a profit maximization formulation for the objective function:
\begin{equation}
\max z=\sum_{m \in \M} \sum_{i \in \P_m}f_{im} q_{im} y_i -\sum_{k \in \K}\sum_{i \in \N}\sum_{j \in \N}  c_{ij} x_{ij}^k  
\label{obj}
\end{equation}

The first term of the objective function \eqref{obj} represents the total revenue of the DRT service operator and the second term represents the total routing cost of the operations. Note that $q_{im}$ represents the load of request $i$ in class $m$, considering different pricing strategies on the number of passengers could be applied to different classes. With users' preferences and chance constraint incorporated into the classic three-index DARP formulation, we summarize the mixed-integer linear programming formulation for the CC-DARP model in Model \ref{md5}.

\begin{model}
	\allowdisplaybreaks
	\begin{subequations}
	\begin{align}
	&\max z=\sum_{m \in \M} \sum_{i \in \P_m}f_{im} q_{im} y_i  -\sum_{k \in \K}\sum_{i \in \N}\sum_{j \in \N} c_{ij} x_{ij}^k, \label{1a}\\
	&\text{\emph{\underline{subject to:}}} \notag\\		
	&\text{\emph{Classic DARP constraints:}}\notag\\
	&\sum_{j \in \N} x_{ij}^k - \sum_{j \in \N}x_{n+i,j}^k = 0, &&\forall i \in \P,  k \in \K,\label{3}\\
	&\sum_{j \in \N} x_{0j}^k=1,  \qquad &&\forall k \in \K,\label{4}\\
	&\sum_{i \in \N} x_{i,2n+1}^k=1, \qquad &&\forall k \in \K,\label{5}\\
	&\sum_{j \in \N} x_{ji}^k - \sum_{j \in \N} x_{ij}^k = 0, \quad   &&\forall i \in \P\cup \D, k \in \K,\label{6}\\
	&B_j\geq(B_0^k+d_0+t_{0j})x_{0j}^k, \qquad   &&\forall  j \in \N, k \in \K, \label{36}\\
	&B_j\geq(B_i+d_i+t_{ij})\sum_{k \in \K}x_{ij}^k, \qquad   &&\forall i , j \in \N, \label{37}\\
	&B_{2n+1}^k\geq(B_i+d_i+t_{i,2n+1})x_{i,2n+1}^k, \qquad   &&\forall i \in \N, k \in \K, \label{38}\\
	&Q_j\geq(Q_0^k+q_{jm})x_{0j}^k, \qquad   &&\forall m \in \M,  j \in \NM, k \in \K,  \label{8a}\\
	&Q_j\geq(Q_i+q_{jm})\sum_{k \in \K}x_{ij}^k, \qquad   &&\forall i \in \N,  m \in \M,  j \in \NM, \label{40}\\
	&Q_{2n+1}^k\geq(Q_i+q_{2n+1,m})x_{i,2n+1}^k, \qquad   &&\forall i  \in \N, k \in \K,  \label{41}\\
	&L_i = B_{n+i}-(B_i+d_i), \qquad &&\forall i \in \P, \label{39}\\
	&B_{2n+1}^k - B_0^k \leq T,  \qquad &&\forall k \in \K, \label{10}\\
	&e_i\leq B_i\leq l_i, \qquad &&\forall i \in \P\cup \D,\label{11}\\
	&t_{i,n+i} \leq L_i \leq L,  \qquad &&\forall i \in \P,\label{12}\\
	&\max \left\{0,q_{im}\right\} \leq Q_i \leq \min \left\{Q, Q+q_{im}\right\}, \quad &&\forall m \in \M,  i \in \NM, k \in \K,\label{13}\\		
	&x_{ij}^k \in \{0,1\}, \qquad &&\forall i \in \N, j\in \N, k \in \K,\label{21}\\
	&y_i \in \{0,1\}, \qquad &&\forall i \in \P, \label{22}\\
	&\underline{V}_i \leq V_i \leq \overline{V}_i, \qquad &&\forall i \in \P, \label{vbounds}\\
	&\text{\emph{Utility Functions:}} \quad (\ref{cbuib2}) - (\ref{cbuob2}),\nonumber\\
	&\text{\emph{Chance Constraints:}} \quad (\ref{cbnccc}), \notag\\		
	&\text{\emph{Accept/Reject:}} \quad (\ref{eq:link}).\notag
	\end{align}
	\end{subequations}
	\label{md5}
\end{model}	

Aligned with Cordeau's classic DARP model \citep{cordeau2006branch}, constraints (\ref{3}) guarantees that pick-up and drop-off nodes for each request are visited by the same vehicle, whereas constraints (\ref{4}) and (\ref{5}) ensure that each route starts at the origin depot and ends at the destination depot. Constraints (\ref{6}), (\ref{36})--(\ref{38}) and (\ref{8a})--(\ref{41}) guarantee flow conservation, time and load consistency respectively. Constraint (\ref{39}) defines the ride time of each user, which is bounded by (\ref{12}). The duration of routes are bounded by constraint (\ref{10}), while (\ref{11}) and (\ref{13}) are the time window and capacity constraints respectively.\\

We acknowledge that from a mathematical programming point of view, binary variable $y_i$ is redundant and can be replaced by $\sum_{j \in \N} \sum_{k \in \K} x_{ij}^k$ everywhere. Similarly, variable $L_i$ could also be omitted and represented by $B_i$ in Cordeau's model \citep{cordeau2006branch}. Extensive computational experiments show that the CPU time does not increase significantly with the presence of this additional variable. Hence, we deliberately kept $y_i$, since it makes the connection between routing decision and user selection more intuitive and organic. 

\section{Solution Method}
\label{algo}

We propose a local search based descent heuristic algorithm for the CC-DARP. For simplicity, we assume the passengers from one request cannot be split and serviced by different vehicles.  Starting with an initial solution generated by a customised insertion heuristic based on \cite{jaw1986heuristic}'s sequential insertion heuristics (Section \ref{ini}), the algorithm iterates between executing local search on the search space of user selection (Section \ref{lsy}) and the search space of routing optimisation (Section \ref{lsx}). The algorithm terminates when no better solution can be found and an overview of the entire heuristic is provided in Section \ref{heu}.

\subsection{Construction of an Initial Solution}
\label{ini}

\cite{jaw1986heuristic} proposed a sequential insertion heuristic, which we deem to be a decent starting point for generating initial solutions with a few modifications applied. Requests are indexed and sorted in an increasing order according to their earliest departure time. Note that for an outbound trip $i$ with a tight time window on arrival time (i.e.  $ e_{n+i}\leq B_{n+i} \leq l_{n+i}$), a pickup time window can be inferred using the direct travel time $t_{i,n+i}$ and the maximum ride time $L$:
\begin{align*}
&B_i \leq l_{n+i} - L, && \forall i \in \OB, \nonumber\\
&B_i \geq e_{n+i} - t_{i,n+i},  && \forall i \in \OB. \nonumber
\end{align*}

With specified information on vehicle capacity and fleet, the algorithm processes requests sequentially, inserting one request into vehicles' work-schedule at a time, checking all feasible insertions, and choosing the insertion where its additional cost is minimised. We check for violations of chance constraint \eqref{cbnccc} after each insertion, and it must be satisfied for an insertion to be labelled 'feasible'. Evidently, some requests might not be inserted in this step, given the limited number of vehicles. These temporarily rejected requests are saved in a pool and for possible later insertion within the following local search step.

Although the literature shows that the quality of initial results does not have a significant influence on final results, structuring a good initial solution can still reduce computational costs \citep{xiang2006fast}. In this regard, adjustments are made in Jaw's insertion heuristic. To avoid situations where more profitable requests could not be inserted in a later stage, we adjust the sorting rule. Intuitively, requests with relatively higher expected profits should be inserted first. Breaking down the profitability of each request by looking into constraint (\ref{cbnccc}) and the utility functions in the proposed CC-DARP model, we identify two factors that determine profitability: direct trip distance and decentralisation. To take into account decentralisation, \citep{diana2004new} proposed a decentralisation index $D_i$, which can be used to our algorithm as follows:
\begin{equation*}
D_i = \frac{\sum\limits_{j \in \N \setminus \{i,n+i\} }t_{ij} + \sum\limits_{j \in \N \setminus \{i,n+i\}} t_{n+i,j}}{\sum\limits_{i' \in \P}\sum\limits_{j \in \N \setminus \{i',n+i'\}} t_{i'j} +\sum\limits_{i' \in \P}\sum\limits_{j \in \N \setminus \{i',n+i'\}} t_{n+i',j}}, \qquad \forall i \in \P.
\end{equation*}

In this way, the decentralisation of request $i$ from both the depot and other requests are taken into account. The higher $D_i$ is, the more decentralised request $i$ is. Therefore, requests with a higher decentralisation index $D_i$ should be inserted in priority, given that requests that are more decentralised can be more difficult to be inserted later on. 

Similarly, to take into consideration the impact of direct trip length, we propose a direct travel time index $TT_i$ defined as:
\begin{equation*}
TT_i = \frac{t_{i,n+i}}{\sum\limits_{i' \in \P}t_{i',n+i'}},\qquad \forall i \in \P.
\end{equation*}

Requests with smaller direct travel time index $TT_i$ should be inserted with priority since these requests are more likely to be profitable according to \ref{deltau} and harder to be inserted successfully at a later stage. To combine the impact of these two factors, we design a general index $G_i$ as a weighted sum of $D_i$ and $TT_i$, which utterly determines the adjustments of sequential sorting indices by increasing $G_i$ values (\ref{gi}):
\begin{equation}
G_i = \omega (1-D_i) + (1-\omega)TT_i.
\label{gi}
\end{equation}

The greater the parameter $\omega\in (0,1)$ is, the more weight is put on the influence of decentralisation. After the initial sorting by $e_i$, we consider again each pair of consecutive request $(i-1)$ and $i$, according to the default sequential ranking order. The two requests get swapped if the following inequality is verified:
\begin{equation}
G_{i-1}-G_i\geq \delta. \label{gic}
\end{equation}

Then we consider the $i$-th and $(i-2)$-th requests. The parameter $\delta$ also takes on calibrated values between 0 and 1. With a smaller $\delta$, we emphasise more the impact of original time sequential sorting.\\

In our numerical experiments, parameters $\omega$ and $\delta$ are calibrated using benchmarking instances for the DARP. Additionally, to refine the initialisation process of the routes, we initialise the routes by assigning the first $|\K|$ seed requests with smallest $G_i$, one to each empty vehicle. The seed requests are, considered by our standard, more difficult to insert later on and supposedly more profitable.\\

The remaining requests are inserted one by one following the adjusted insertion order $InOr$ in Algorithm \ref{algo1}. Initial solution $s_0$ is instructed following \cite{jaw1986heuristic}'s insertion method: Schedule blocks (when vehicles are active) and slack periods (when vehicles are idling) of each vehicle are created at initialisation, and updated after each successful insertion. When inserting a request $i$, the algorithm systematically scans all feasible insertions of the request into the existing work schedules of all vehicles, and choose the best insertion with minimum additional cost. At each attempt of inserting request $i$, an insertion is only labelled as feasible if it satisfies time window constraints (\ref{11}), maximum ride time constraints (\ref{12}), load constraints (\ref{13}) and chance constraints (\ref{cbnccc}) for all existing requests in current work schedule. New schedule blocks can also be created when examining for feasible insertions. For a more detailed description of the insertion heuristic we refer to \cite{jaw1986heuristic}'s paper. The goal of this customised insertion heuristic to construct an initial solution $s_0 = (\bm{x},\bm{y})$, which contains both initial user selection and routing decisions required for starting off local search afterwards. The algorithm for our insertion heuristic is summarised in Algorithm \ref{algo1}.\\ 

\begin{algorithm}[ht]
\KwIn{User request data, DRT service data}
\KwOut{Initial solution $s_0$,visited requests $\V$, request pool $\U$}
$\U \gets \emptyset$,$\V \gets \emptyset$\\
\For{$i \in \OB$}{
	update $e_i = e_{n+i} - t_{i,n+i}$
}
Create sequential inserting order list $InOr =[\arg\min_{i \in \P} e_i, ..., \arg\max_{i \in \P} e_i]$ by increasing $e_i$\\
\For{$i \in \P$}{
    Calculate general index $G_i$ using \eqref{gi}\\
}
Replace the first $|\K|$ requests in $InOr$ by the $|\K|$ requests with the smallest $G_i$ (seeds) \\ 
Adjust the remaining $InOr$ accordingly if (\ref{gi}) holds\\
\For{$i \in \P$}{
	\For {$k \in \K$}{
		Try to insert $InOr[i]$ into the work-schedule of vehicle $k$, record the additional cost of all feasible insertions
	}
	\If{No feasible insertion found}{
	Add request $InOr[i]$ into the unvisited request pool $\U$
	}
	\Else{
	Find the feasible insertion of lowest additional cost out of all vehicles and insert request $InOr[i]$\\
	Add request $InOr[i]$ into set of visited requests $\V$ and update $s_0$
	}
}
\caption{Construction of an Initial Solution}
\label{algo1}
\end{algorithm}
 
\subsection{Local Search}

With an initial solution $s_0$ obtained, we conduct local search iteratively on both user selection ($\bm{y}$) and routing variables ($\bm{x}$) until no improvement in the solution can be found.

\subsubsection{Local Search on User Selection ($\bm{y}$)}
\label{lsy}

For an incumbent solution $s = (\bm{x},\bm{y}) $, the neighbourhood $N(s)$ of solution $s$ is defined as the set of all solutions that can be obtained by applying a single operator, namely \add, \rem and \swap. Operators \add and \rem are defined as: 
\begin{itemize}
	\item \add: Add (insert) a request into vehicle work-schedules. This operation can increase the number of operating vehicles, if it is cheaper to assign an empty vehicle to service this additional request. The target vehicles for insertion can be specified with a subscript (e.g. \add$_k$ indicates the operation of only inserting the request into the work schedule of vehicle $k$). Note that there might be more than one feasible insertions when applying this operator. In this work, the \textit{best insertion} is adopted, i.e. the operator scans all feasible insertions and returns the insertion $s_i$ with minimum additional cost into the neighbourhood set. Otherwise, the operator returns null if no feasible insertion can be found. 
	\item \rem: Remove a request from its corresponding vehicle work-schedule and return the remaining work schedule. The rest of the routing sequence does not change, as the trip simply skips the removed request. This operation can reduce the number of operating vehicles as well as active schedule blocks in a vehicle. 
\end{itemize}

A two-step local search is carried out on user selection. In step one, we first explore the neighbourhood $N_1(s)$ of solution $s$, where the elements in this set are formed by applying either \rem or \add on solution $s$. In other words, all elements in $N_1(s)$ have a different $\sum_{i \in \P} y_i$ value compared to solution $s$ (either one less or one more). A subsequent local search step on search space $N_1(s)$ is then applied to $s$, yielding $s'$. If $s'$ is a better solution than $s$ in terms of objective function $z$ (profit), it replaces $s$ and the search continues on $N_1(s')$ until no better solution can be found.\\

For step two, we define an additional operator \swap:
\begin{itemize}
\item \swap: This operation explores all feasible combinations of simultaneously removing one existing request and adding one request from the pool of unserved requests $\U$, ensuring the total number of serviced requests remains unchanged. Only the best work schedule $s_{ij}$is returned if more than one feasible insertions are found.
\end{itemize} 

Note that \swap is essentially \add and \rem executed simultaneously. We add this operator simply to broaden the search space when searching for the best selection of requests, as solutions with the same size of serviced requests but with swapped elements cannot be obtained by applying \add or \rem on $s$. With a local optimal solution $s$ obtained from the local search on $N_1(s)$, The algorithm then moves to step two, where it explores the extended neighbourhood $N_2(s)$ of solutions obtained by applying \swap on $s$. In $N_2(s)$, all solutions have the same $\sum_{i \in \P} y_i$ as $s$. Similarly, a local search step is applied to $s$ on $N_2(s)$. The algorithm loops over step one and step two until no improvement can be found in both steps. This two-step local search algorithm is summarized in Algorithm \ref{algo2}.

\begin{algorithm}[H]
\KwIn{initial solution $s_0$, visited requests $\V$, request pool $\U$}
\KwOut{$s_1$}
\tcc{Initialisation}
$s$ $\gets$ $s_0$\\
$z(s'') = z(s) + 0.0001$\\
\tcc{Local Search}
\While{$z(s'') > z(s)$}{
  	\tcc{Step one}
  	\Repeat{$z(s') = z(s)$}{
      	$N_1(s) \gets \{s\}$\\
      	\For{$i \in \U$}{
      		\If{Feasible insertion(s) found}{
      		    $s_i \gets$ Apply \add to $i$\\
      		    $N_1(s) \gets N_1(s) \cup \{s_i\}$\\
      		}
      	}
        \For{$j \in \V$}{
            $s_j \gets$ Apply \rem to $j$ \\
      	    $N_1(s) \gets N_1(s) \cup \{s_j\}$\\
      	}
      	$s' \gets \arg\max\{z(s') : s' \in N_1(s)\}$\\
      	\If{$z(s') > z(s)$}{
      		$s$ $\gets$ $s'$\\
      		Update $\U$ and $\V$
      	}
    }
    \tcc{Step two}
    \Repeat{$z(s'') = z(s')$}{
    	$N_2(s) \gets \{s'\}$\\
    	\For{$i \in \V$, $j \in \U$}{
			\If{Feasible insertion(s) found}{
				$s_{ij} \gets$ Apply \swap to $(i,j)$\\
				$N_2(s) \gets N_2(s) \cup \{s_{ij}\}$\\
			}
    	}
    	$s'' \gets \arg\max\{z(s'') : s'' \in N_2(s)$\}\\
    	\If{$z(s'') > z(s')$}{
    		$s'$ $\gets$ $s''$\\
    		Update $\U$ and $\V$
    	}
    }
}
$s_1 \gets s''$
\caption{Local Search on User Selection}
\label{algo2}
\end{algorithm}
 
\subsubsection{Local Search on Routing Optimisation ($\bm{x}$)}
\label{lsx}

Once the local search on user selection terminates, a solution with the current best selection of users is passed on to the second part of local search: local search on routing optimisation. In this local search process, we fix the set of selected users and optimise the routing variables ($\bm{x}$) in a classic DARP manner, where all selected requests must be serviced. In this local search, we define the following two operators for a request $i$ that is assigned to vehicle $k$ in the current solution $s$:

\begin{itemize}
\item \ra: This operator Re-Assigns request $i$ to another vehicle. This can be broken down into the combination of operators \rem$_k$ and \add$_{K \setminus k}$ applied on $i$ (i.e. to remove request $i$ from vehicle $k$ and to insert it into any other vehicle). This operator conducts inter-tour exchanges.
\item \ri: This operator tries to Re-Insert request $i$ into its originally assigned vehicle $k$. The operator examines all feasible insertions in this vehicle, searching for a possibly better insertion rather than just inserting it back into its original place. Similarly, this can be broken down into the combination of operators \rem$_k$ and \add$_{k}$ applied on $i$ (i.e. to remove request $i$ from vehicle $k$ and then to insert it back into vehicle $k$). This operator takes into account intra-tour exchanges.
\end{itemize}

For each visited request $i$ in solution $s$, we first apply \ra  on the request, searching for feasible ways of inserting this request into another trip. If \ra fails to find any feasible insertion, \ri is then carried out to insert the request back to its original trip. In both operations, \textit{best insertion} is adopted. The neighbourhood $N_3(s)$ of solution $s$ is composed of all solutions that can be obtained from $s$ by applying one of these two operators. This way, local search is conducted with both inter-tour and intra-tour exchanges. With \ra and \ri carried out in parallel, the efficiency of local search is greatly improved \citep{xiang2006fast}. For each solution $s$, a local search is conducted on $N_3(s)$, yielding solution $s'$. If $s'$ is a better solution than $s$ in terms of objective function $z$ (profit), $s'$ replaces $s$ and the local search restarts on the new incumbent solution. The algorithm terminates when no better solution can be found. The detailed algorithm of local search on routing optimisation is summarised in Algorithm \ref{algo3}.\\

\begin{algorithm}[H]
\KwIn{$s_1$ from Algorithm \ref{algo2}, $V$}
\KwOut{$s_2$}
\tcc{Initialisation}
$s \gets s_1$\\
$z(s') = z(s) + 0.0001$\\
\tcc{Local Search}
\Repeat{$z(s') = z(s)$}{
    $N_3(s) \gets \{s\}$\\
    \For{$i \in \V$}{
        \tcc{Inter-tour Search}
        \If{Feasible insertion(s) found}{
			$s_{i} \gets$ Apply \ra to $i$\\
			$N_3(s) \gets N_3(s) \cup \{s_{i}\}$\\            
        }
        \Else{
            \tcc{Intra-tour Search}
			$s_{i} \gets$ Apply \ri to $i$\\
			$N_3(s) \gets N_3(s) \cup \{s_{i}\}$\\
        }
    }
    $s' \gets \arg\max\{z(s') : s' \in N_3(s)\}$\\
    \If{$z(s') > z(s)$}{
        $s \gets$ $s'$\\
    }
}
$s_2 \gets s'$\\
\caption{Local Search on Routing Optimisation}
\label{algo3}
\end{algorithm}

\subsection{Algorithm Summary}
\label{heu}

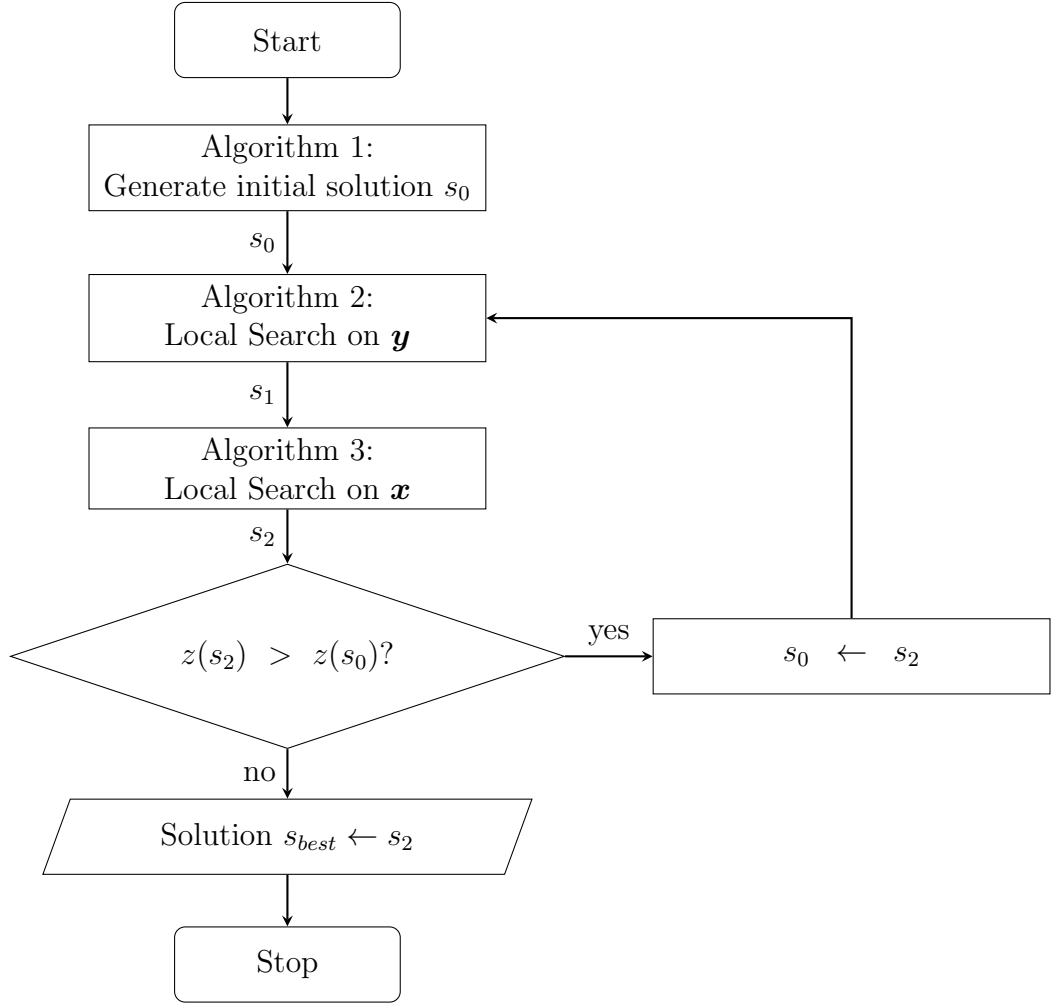
\begin{figure}[H]
\begin{tikzpicture}[node distance=1.5cm]
\node (start) [startstop] {Start};
\node (pro1) [process, below of=start,yshift=-.2cm] {Algorithm 1:\\ Generate initial solution $s_{0}$};
\node (pro2) [process, below of=pro1, yshift=-.5cm] {Algorithm 2:\\ Local Search on $\bm{y}$};
\node (pro3) [process, below of=pro2, yshift=-.5cm] {Algorithm 3:\\ Local Search on $\bm{x}$};
\node (dec1) [decision, below of=pro3, yshift=-1cm] {$z(s_2)>z(s_0)$?};
\node (pro2b) [process, right of=dec1, xshift=6cm] {$s_0\gets s_2$};
\node (out1) [io, below of=dec1, yshift=-0.9cm] {Solution $s_{best} \gets s_2$};
\node (stop) [startstop, below of=out1,yshift=-.2cm] {Stop};
\draw [arrow] (start) -- (pro1);
\draw [arrow] (pro1) -- node[anchor=east] {$s_0$}(pro2);
\draw [arrow] (pro2) -- node[anchor=east] {$s_1$}(pro3);
\draw [arrow] (pro3) -- node[anchor=east] {$s_2$}(dec1);
\draw [arrow] (dec1) -- node[anchor=east] {no} (out1);
\draw [arrow] (dec1) -- node[anchor=south] {yes} (pro2b);
\draw [arrow] (out1) -- (stop);
\draw [arrow] (pro2b) |- (pro2);
\end{tikzpicture}
\caption{LS-H Solution Method Overlook}
\label{fig1}
\end{figure}

After an initial solution is constructed, we iteratively conduct local search on user selection and routing optimisation in a loop until there is no improvement. We name this customised heuristic solution method 'Local Search-based Heuristic (\textbf{LS-H})'. A summary of LS-H is depicted in Figure \ref{fig1}.

We acknowledge that even though the proposed LS-H is proved to be capable of jumping out of local optima with the help of its layered loops, this solution method is substantially still a classic heuristic approach. Therefore, it can still get stuck in local optima occasionally. Metaheuristics, including Tabu Search and VNS, may yield results that are closer to optimal solutions. Since this work aims at exploring various fare structure formulations and their influence on class-based ridership, we leave the development of more advanced solution algorithms for future research.

\section{Numerical Experiments}
\label{num}
We conduct numerical experiments on both benchmarking instances from the DARP literature \citep{cordeau2006branch} and using a realistic case study generated based on extracted data from NYC yellow taxi trip data provided by NYC OpenData. We use DARP benchmarking instances to assess the performance of the proposed local search based heuristic by comparing with the best integer solution returned by an exact implementation of Model 1 using MILP. All algorithms and formulations are implemented in Python on a personal computer with 16GB of RAM and an Intel i7 processor at 2.9GHz. and CPLEX 12.10 is used to solve the CC-DARP MILPs. 

\subsection{Experiments Design}

\begin{figure}
\begin{center}
\includegraphics[scale=0.5]{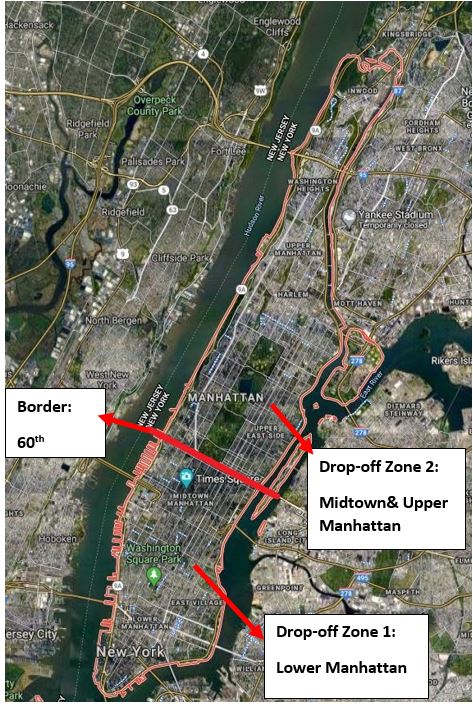}
\caption{Drop-off zones in Manhattan}
\label{mht}
\end{center}
\end{figure}

\begin{figure}
\centering
\begin{subfigure}{0.55\linewidth}{\includegraphics[width=0.95\textwidth,valign=c]{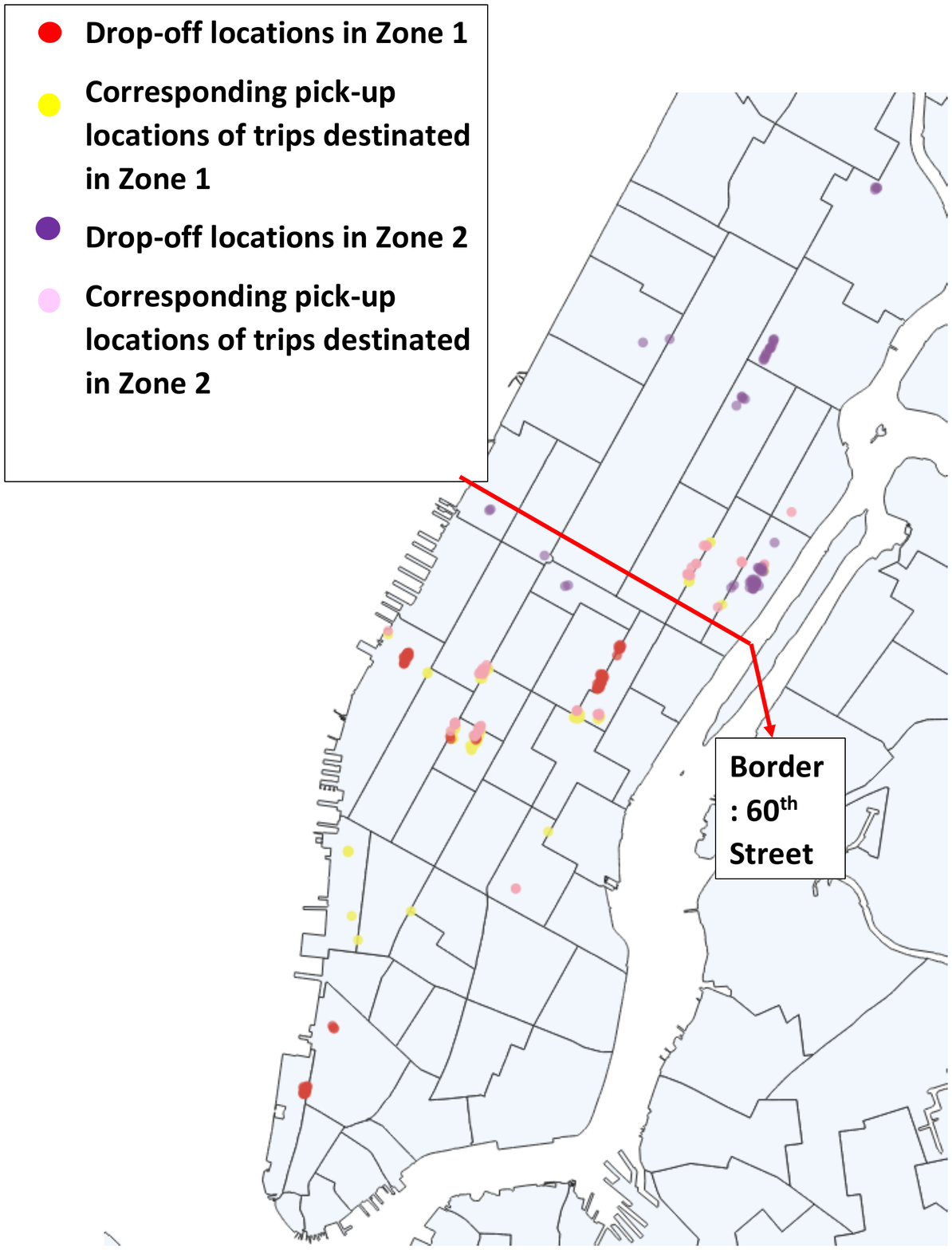}\caption{Map of pickup and dropoff node locations. \label{map2}}}
\end{subfigure}
\begin{subfigure}{0.4\linewidth}{\includegraphics[width=0.95\textwidth,valign=c]{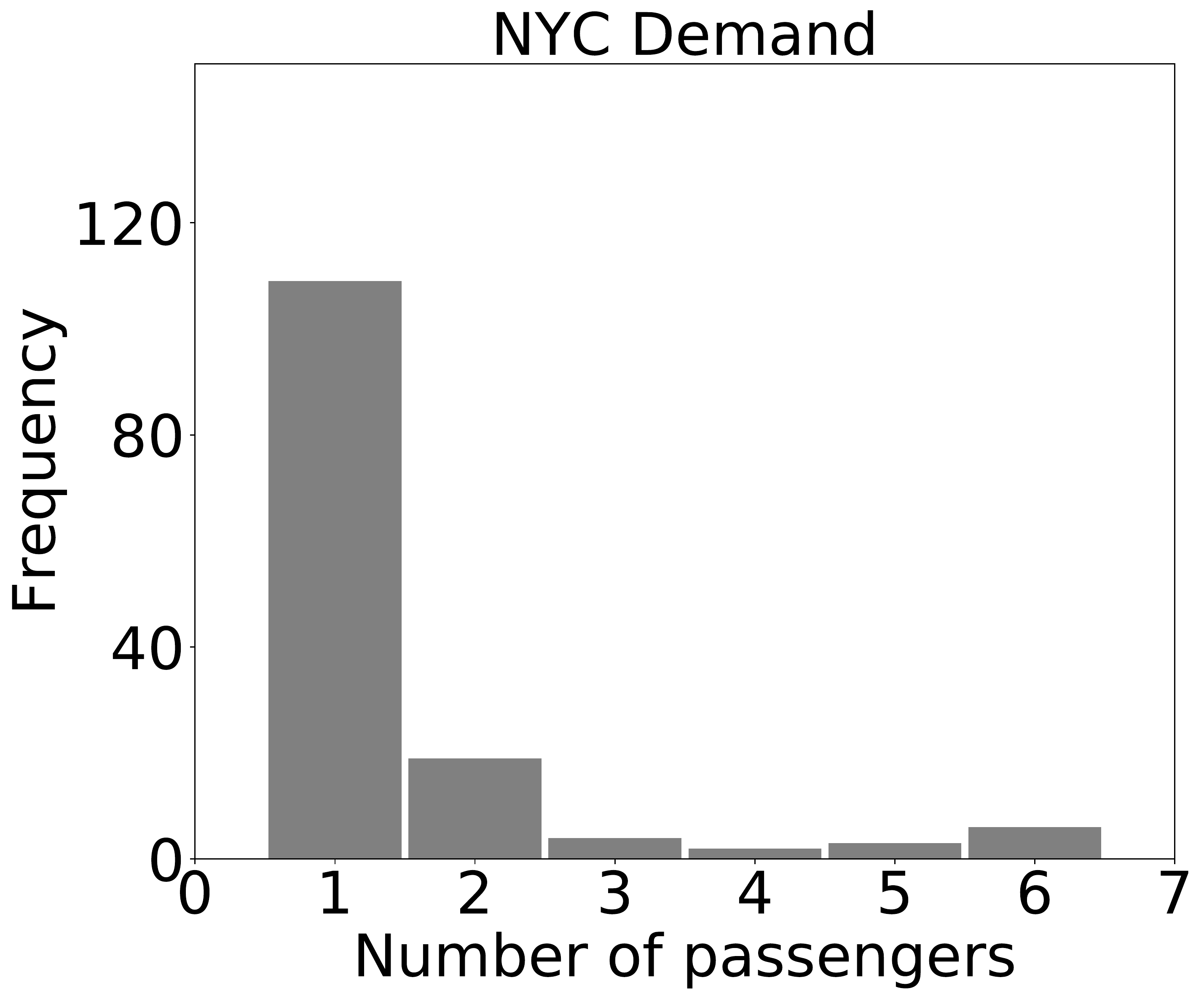}\caption{Histogram of number of passengers per request.\label{hist}}}
\end{subfigure}
\caption{NYC Dataset (143 requests, 286 nodes)}
\label{freq}
\end{figure}

The first set of instances we use in our experiments corresponds to the 'A' instances of \cite{cordeau2006branch} ({\small \url{http://neumann.hec.ca/chairedistributique/data/darp/}}) with homogeneous vehicles with identical capacities $Q=3$. Instances in this dataset contains 16-96 requests covered by 2-8 vehicles. The requests in each instance are evenly divided into 50\% outbound trips and 50\% inbound trips. The coordinates of the pickup and drop-off locations are randomly generated and uniformly distributed in the [-10, 10] $\times$ [-10, 10] square. For each request, a 3-min service time and a 15-min time window is imposed on its departure or arrival time, depending on if it is an inbound or outbound trip. The maximum ride time $L$ is set to be 30 minutes, while the length of planning horizon $T$ is 60 minutes.  \\

The second set of instances is extracted from Yellow Taxi Trip Data in NYC (2006) ({\small \url{ https://data.cityofnewyork.us/Transportation/2016-Yellow-Taxi-Trip-Data/k67s-dv2t}}). We extracted 23,772 trip records of the morning peak hours (6am-9am) on a random weekday (1st of February 2016), which contain the information of the pickup and drop off coordinates, number of passengers, pickup and drop off time of each trip. For preparation of the further discussion on fare structures, we divide Manhattan into two areas as shown in Figure \ref{mht}, and filter out the trips with a drop-off location outside zone 1 and zone 2. Note that we only filter out the trips with destinations outside Manhattan, the pickup locations of the trips are not limited. On the map, we split the area into 120,000 35m $\times$ 30m mutually exclusive sub-regions, and label each trip with two sub-regions: one for pickup coordinates and one for drop off coordinates. We then obtain the most visited sub-regions by analysing the frequency of each sub-region is labelled. By mapping the 50 most popular pickup sub-regions overall and the 20 most popular drop off sub-regions in each drop-off zone, we finally select 143 eligible trips with most frequently visited origins and destinations and various numbers of passengers in each request (see Figure \ref{freq}). These trip records were further reformed into instances of requests for testing our CC-DARP model. Note that as these trips are observed data from taxi meter records, taxi is proven to be these passengers' best travel option, which fits in our model where these trips are compared against DRT as the alternative mode.

\subsection{LS-H Performance Benchmark}

We first benchmark our Local Search based heuristic (LS-H) using Cordeau's instances under a flat fare structure. We conduct sensitivity analysis on $p$, $s$ of chance constraint (\ref{nncc}) and flat fare $f$. In the following experiments, we assume that users are all rational and seeking to maximise their own utility. Parameters $\beta_T$ is set to $\$10.6$/h based on a cost and benefit analysis conducted in Australia and New Zealand \citep{litman2009transportation}. $\beta_S$ is set to $\$21.2$/h, given that the value of waiting time is approximately two to three times of the value of travel time \citep{mohring1987values}. $\beta_F$ is set to $10/\$$ after calibration. $\omega$ and $\delta$ of LS-H are calibrated to be 0.1 and 0.03, respectively. Note that in this subsection the monetary unit '\$' in this paper refers to Australian Dollar (AUD).

Each instance of the CC-DARP is solved with an exact MILP solver and LS-H separately, CPU times and optimality gaps are then reported. For the exact solution counterpart, we solve CC-DARP using CPLEX with default options and apply some of the pre-processing user cuts (namely 'Bounds on time and load variables', 'Capacity constraints', 'Precedence Constraints' and 'Generalised order constraints') proposed by \cite{cordeau2006branch}. We use a time limit of 1 hour for all instances and report time outs with the symbol $^\dagger$.

\subsubsection{Sensitivity Analysis on $p$ and $s$}
\label{saps}

Tables \ref{sap} and \ref{sas} summarize the sensitivity analysis results of $p$ and $s$ on Cordeau's instances. Outcomes of both solution methods including objective value (Profit), CPU time and the number of requests accepted denoted $\sum y_i^\star$ and $\sum \hat{y}_i$ for CPLEX and LS-H, respectively, are reported and compared. The best lower bound (integer feasible solution) and the relative optimality gap (Gap) are reported if CPLEX could not solve the problem to optimality within the available time limit. Profits are highlighted in bold if the optimal solution is found by LS-H. In the Profit gap (\%) column, a negative value indicates that LS-H outperforms CPLEX's best lower bound within the time limit. These results show that the performance of our heuristic is rather stable when chance constraints (\ref{cbnccc}) are binding (i.e. higher $p$ and $s$ values). The average profit gap for Table \ref{sap} and \ref{sas} are $3.34\%$ and $2.88\%$, respectively. Furthermore, the last column of both tables shows that for user selection, LS-H almost always manages to select the optimal subset of requests to service within a few minutes. This suggests that although our heuristic may fail to identify optimal solution by a small margin, LS-H is frequently capable of identifying the optimal set of user requests to accept, which is valuable data for DRT operators.

\begin{table}[H]
\caption{Sensitivity Analysis on $p$ ($s=10, f=20$)}\label{sapopt}
\begin{center}
	\footnotesize
	\resizebox{\textwidth}{!}{
		\begin{tabular}{l l l l l l l l l l l }
			\toprule
			\multirow{2}{*}{Instance} &\multirow{2}{*}{$p$} & \multicolumn{4}{l}{CPLEX}& \multicolumn{3}{l}{LS-H} & \multicolumn{2}{l}{Comparison} \\
			\cmidrule(l){3-6}
			\cmidrule(l){7-9}
			\cmidrule(l){10-11}
			&  & Profit (\$) & Gap (\%) & CPU time (s) & $\sum y_i^\star$ &  Profit (\$) & CPU time (s) & $\sum \hat{y}_i$   &Profit gap (\%)  & $\sum y_i^\star - \sum \hat{y}_i$\\ \midrule
			\multirow{3}{*}{a2-16} 
			 &0.8  &\textbf{114.84}  &$< 1e^{-4}$  &0.16 &8  &\textbf{114.84}  &4.60  &8 &0.00   &0 \\ 
			 &0.95 &\textbf{100.95}	 &$<1e^{-4}$  &0.08  &7  &\textbf{100.95}	&3.82  &7 &0.00   &0 \\ 
			 &0.99 &\textbf{100.95}  &$<1e^{-4}$  &0.09  &7  &\textbf{100.95}	&3.70  &7 &0.00  &0\\ \cmidrule(l){1-11}
			 \multirow{3}{*}{a2-20}
			 &0.8  &156.95  &$<1e^{-4}$  &0.08  &10   &148.78 &4.79 &10&5.21  &0 \\ 
			 &0.95 &123.13  &$<1e^{-4}$   &0.06 &8    &120.73 &4.52 &8 &1.95  &0  \\ 
			 &0.99  &84.48 &$<1e^{-4}$  &0.03  &6    &83.83  &3.81 &6 &0.78  &0\\ \cmidrule(l){1-11}
			 \multirow{3}{*}{a2-24} 
			 &0.8  &225.29  &$<1e^{-4}$   &1.27  &15  &223.51  &7.73  &15 &0.79   &0\\ 
			 &0.95 &206.04  &$<1e^{-4}$   &1.61  &14  &205.61  &7.17  &14 &0.21   &0  \\ 
			 &0.99  &174.61 &$<1e^{-4}$   &0.36 &12  &174.47   &6.93  &12 &0.08   &0\\ \cmidrule(l){1-11}
			 \multirow{3}{*}{a3-18} 
			 &0.8  &142.26  &$<1e^{-4}$   &0.42  &9  &137.92  &4.74  &9 &3.05 &0 \\ 
			 &0.95 &142.26  &$<1e^{-4}$   &0.48  &9  &137.92  &4.59  &9 &3.05 &0\\ 
			 &0.99  &94.32  &$<1e^{-4}$   &0.23  &6  &89.97   &3.85  &6 &4.61 &0\\ \cmidrule(l){1-11}
			 \multirow{3}{*}{a3-24} 
			 &0.8  &150.66  &$<1e^{-4}$   &1.25  &10 &149.17  &6.62  &10 &0.98   &0 \\ 
			 &0.95 &135.90  &$<1e^{-4}$   &0.36  &9  &133.40  &5.82  &9  &1.84  &0  \\ 
			 &0.99  &135.90 &$<1e^{-4}$   &0.38  &9  &133.40  &5.99  &9  &1.84   &0\\ \cmidrule(l){1-11}
			 \multirow{3}{*}{a3-36} 
			 &0.8  &277.80  &$<1e^{-4}$  &10.08 &18  &273.09  &14.04  &18 &1.70    &0\\ 
			 &0.95 & 245.67 &$<1e^{-4}$  &11.38 &16  &240.35  &12.86  &16 &2.17    &0 \\ 
			 &0.99 &209.83 &$<1e^{-4}$   &4.53  &14  &206.03  &12.09  &14 &1.81   &0\\ \cmidrule(l){1-11}
			 \multirow{3}{*}{a4-16} 
			 &0.8  &157.79  &$<1e^{-4}$   &1.98  &10 &155.54  &4.64  &10 &1.42   &0 \\ 
			 &0.95 &142.07  &$<1e^{-4}$   &0.77  &9  &139.74  &4.46  &9  &1.64   &0 \\ 
			 &0.99  &92.33  &$<1e^{-4}$   &0.25  &6  &90.25   &3.49  &6  &2.26   &0\\ \cmidrule(l){1-11}
			\multirow{3}{*}{a4-32} 
			&0.8  &237.22  &$<1e^{-4}$   &46.31  &15  &224.23  &15.03  &15  &5.48   &0 \\ 
			&0.95 &203.32  &$<1e^{-4}$   &5.09  &13  &192.34  &13.20  &13  &5.40  &0 \\ 
			&0.99  &203.32 &$<1e^{-4}$   &5.06  &13  &192.34  &13.63  &13  &5.40  &0\\ \cmidrule(l){1-11}
			\multirow{3}{*}{a4-40}
		    &0.8  &325.68  &$<1e^{-4}$   &357.78 &20  &306.30  &31.13  &20 &5.95   &0 \\ 
			&0.95 &261.43  &$<1e^{-4}$   &4.73  &16  &247.83  &19.38  &16 &5.20   &0  \\ 
			&0.99  &243.77  &$<1e^{-4}$  &2.34  &15 &230.41  &17.20  &15  &5.48 &0 \\ \cmidrule(l){1-11}
			\multirow{3}{*}{a4-48}
		    &0.8  &281.02  &$<1e^{-4}$  & 178.66  &18  &270.80  &41.80  &18 &3.64 &0\\ 
			&0.95 &263.49  &$<1e^{-4}$  &63.08   &17  &254.61  &25.60  &17 &3.37 &0 \\
			&0.99  &196.61 &$<1e^{-4}$  &3.23   &13   &190.87  &17.23  &13  &2.92  &0\\ \cmidrule(l){1-11}
			\multirow{3}{*}{a5-50}
			&0.8  &400.63  &5.49   &3600$^\dagger$ &25  &352.20  &57.03  &22 &12.13  &3\\ 
			&0.95 &383.00  &5.35  &3600$^\dagger$  &24  &298.82  &28.21  &19 &21.98   &5 \\ 
			&0.99  &322.72 &3.16  &3600$^\dagger$  &21  &295.94 &50.73  &19  &11.05  &2\\ \cmidrule(l){1-11}
			\multirow{3}{*}{a6-60}
			&0.8  &571.13  &8.87  &3600$^\dagger$  &35  &546.59  &99.56  &35 &4.30  &0\\ 
			&0.95 &465.03  &7.28  &3600$^\dagger$  &29  &445.80  &78.26  &29 &4.14  & 0\\ 
			&0.99  &381.03 &6.19  &3600$^\dagger$  &24  &367.17  &59.58  &24  &3.64  &0\\ \cmidrule(l){1-11}				
			\multirow{3}{*}{a6-72}
			&0.8  &563.97  &7.99  &3600$^\dagger$  &35  &514.89  &81.33  &33 &8.70   &2\\ 
			&0.95 &516.36  &6.26  &3600$^\dagger$  &32  &458.67  &159.25  &30 &11.17 & 2\\ 
			&0.99  &443.71 &3.63  &3600$^\dagger$  &28  &409.30  &69.51  &27  &7.76   &1\\ \cmidrule(l){1-11}				
			\multirow{3}{*}{a7-84}
			&0.8  &786.59  &15.95  &3600$^\dagger$  &50  &792.41  &292.87  &50 &-0.74    &0\\ 
			&0.95 &710.42  &13.78  &3600$^\dagger$  &45  &710.45  &557.46  &46 &0.00    &-1 \\ 
			&0.99  &588.80 &11.55  &3600$^\dagger$  &37 &573.17  &284.86  &37  &2.65  &0\\ \cmidrule(l){1-11}			
			\multirow{3}{*}{a8-96}
			&0.8  &760.81  &16.63  &3600$^\dagger $  &49  &779.66  &276.95  &49 &-2.48  &0\\ 
			&0.95 &624.68  &35.98  &3600$^\dagger$  &40  &738.16  &315.50  &47 &-18.17    &-7 \\ 
			&0.99  &628.39 &13.82  &3600$^\dagger$  &40 &623.24  &335.50  &40  &0.82   &0\\ \bottomrule
			
			 $\dagger$: CPLEX timed out
			
		
		\end{tabular}
        \label{sap}
	}
\end{center}
\end{table}

\begin{table}[H]
\caption{Sensitivity Analysis on $s$ ($p=0.95, f=20$)}\label{sasopt}
\begin{center}
\footnotesize
\resizebox{\textwidth}{!}{
	\begin{tabular}{l l l l l l l l l l l }
			\toprule
			\multirow{2}{*}{Instance} &\multirow{2}{*}{$s$} & \multicolumn{4}{l}{CPLEX}& \multicolumn{3}{l}{LS-H} & \multicolumn{2}{l}{Comparison} \\
			\cmidrule(l){3-6}
			\cmidrule(l){7-9}
			\cmidrule(l){10-11} 
		&  & Profit (\$) & Gap (\%) & CPU time (s) & $\sum y_i^\star$ &  Profit (\$) & CPU time (s) & $\sum \hat{y}_i$   &Profit gap (\%)  & $\sum y_i^\star - \sum \hat{y}_i$\\ \cmidrule(l){1-11}
			\multirow{3}{*}{a2-16} 
			&1  &134.86  &$<1e^{-4}$  &0.14 &9   &129.98 		&4.86  &9 &3.62   &0 \\ 
			&10 &\textbf{100.95}  &$<1e^{-4}$  &0.08  &7  &\textbf{100.95}	&3.58  &7 &0.00   &0 \\ 
			&50 &\textbf{54.38} &$<1e^{-4}$   &0.03  &4  &\textbf{54.38}	&2.64  &4 &0.00   &0\\ \cmidrule(l){1-11}
			\multirow{3}{*}{a2-20}
			&1  &156.95  &$<1e^{-4}$  &0.08  &10   &148.78		 &4.64 &10&5.21  &0 \\ 
			&10 &123.13  &$<1e^{-4}$   &0.06  &8    &120.73 		 &4.46 &8 &1.95  &0  \\ 
			&50 &\textbf{37.71}  &$<1e^{-4}$   &0.03  &3    &\textbf{37.71}  &2.60 &3 &0.00  &0\\ \cmidrule(l){1-11}
			\multirow{3}{*}{a2-24} 
			&1  &244.98  &$<1e^{-4}$   &2.55  &16  &241.71  		&7.35  &16 &1.33   &0\\ 
			&10 &206.04  &$<1e^{-4}$   &1.59  &14  &205.61  		&7.35  &14 &0.21   &0  \\ 
			&50 &\textbf{11.43}  &$<1e^{-4}$     &0.02  &1  &\textbf{11.43}   &2.35  &1 &0.00   &0\\ \cmidrule(l){1-11}
			\multirow{3}{*}{a3-18} 
			&1  &157.55  &$<1e^{-4}$   &0.7  &10  &153.68 		 &5.06 &10 &2.46 &0 \\ 
			&10 &142.26  &$<1e^{-4}$   &0.48  &9  &137.92  		&4.85   &9 &3.05 &0\\ 
			&50  &\textbf{12.11}  &$<1e^{-4}$   &0.02  &1  &\textbf{12.11}   &3.07 &1 &0.00 &0\\ \cmidrule(l){1-11}
			\multirow{3}{*}{a3-24} 
			&1  &151.94  &$<1e^{-4}$   &1.03  &10 &149.17 		 &6.84  &10 &1.82   &0 \\ 
			&10 &135.90  &$<1e^{-4}$   &0.38  &9  &133.40  		&6.08  &9  &1.84   &0  \\ 
			&50  &\textbf{23.32}  &$<1e^{-4}$   &0.05  &2  &\textbf{23.32}  &2.68  &2  &0.00   &0\\ \cmidrule(l){1-11}
			\multirow{3}{*}{a3-36} 
			&1  &311.95  &$<1e^{-4}$  &22.33 &20  &308.37 		 &16.60  &20 &1.15    &0\\ 
			&10 &245.67 &$<1e^{-4}$  &11.19 &16  &240.35		 &13.11  &16 &2.17    &0 \\ 
			&50 &\textbf{80.55} &$<1e^{-4}$   &0.11  &6  &\textbf{80.55}  &5.74    &6  &0.00    &0\\ \cmidrule(l){1-11}
			\multirow{3}{*}{a4-16} 
			&1  &173.95  &$<1e^{-4}$   &4.14  &11 &171.25  &4.76  &11 &1.55   &0 \\ 
			&10 &142.07  &$<1e^{-4}$   &0.75  &9  &139.74  &4.38  &9  &1.64   &0 \\ 
			&50  &\textbf{13.01}  &$<1e^{-4}$   &0.02 &1  & \textbf{13.01}	 &2.18  &1  &0.00   &0\\ \cmidrule(l){1-11}
			\multirow{3}{*}{a4-32} 
			&1 &237.22  &$<1e^{-4}$   &31.11  &15  &224.23 		 &14.89  &15  &5.48   &0 \\ 
			&10 &203.32 &$<1e^{-4}$  &5.08 &13  &192.34 		 &13.21  &13  &5.40   &0 \\ 
			&50 &\textbf{25.70}  &$<1e^{-4}$  &0.06 &2  &\textbf{25.70}  &2.89  &2  &0.00   &0\\ \cmidrule(l){1-11}
			\multirow{3}{*}{a4-40}
			&1  & 343.38 &$<1e^{-4}$   &778.42 &21  &326.23  &20.71  &21 &4.99    &0\\ 
			&10 &261.43  &$<1e^{-4}$  &4.78  &16   &247.83  &17.37  &16 &5.20   & 0 \\ 
			&50 &\textbf{56.05}   &$<1e^{-4}$  &0.09  &4    &\textbf{56.05}   &3.99  &4  &0.00  & 0\\ \cmidrule(l){1-11}
			\multirow{3}{*}{a4-48}
			&1  &349.28  &$<1e^{-4}$ &2162.97  &22 &319.20  &27.44  &21 &8.61   & 1\\ 
			&10 &263.40  &$<1e^{-4}$ &34.28  &17  &254.61  &21.82  &17 &3.34   & 0 \\ 
			&50 &85.29  &$<1e^{-4}$  &0.17  &6   &71.38  &4.69     &5  &16.31  &1 \\ \cmidrule(l){1-11}
			\multirow{3}{*}{a5-50}
			&1  &400.23  &5.57  &3600$^\dagger$   &25  &352.02  &48.51  &22   &12.04 &3\\ 
			&10 &383.00  &5.36  &3600$^\dagger$  &24  &298.82  &23.93  &19 & 21.98  &5  \\ 
			&50 &\textbf{37.98}    &$<1e^{-4}$ &0.16     &3   &\textbf{37.98}  &4.38  &3  &0.00  &0 \\ \cmidrule(l){1-11}
			\multirow{3}{*}{a6-60}
			&1  &607.86  &9.01  &3600$^\dagger$   &37  &580.82  &86.55  &37 &4.45   &0 \\ 
			&10 &465.03  &7.31  &3600$^\dagger$  &29  &455.80  &68.04  &29 &4.14   & 0 \\ 
			&50  &81.29  &$<1e^{-4}$  &0.27   &6    &80.32  &7.60  &6  &1.19  &0 \\ \cmidrule(l){1-11}				
			\multirow{3}{*}{a6-72}
			&1  &597.80  &8.32   &3600$^\dagger$  &37  &544.63  &75.45  &35 &8.89   &2 \\ 
			&10 &516.53  &6.29  &3600$^\dagger$  &32  &458.67  &137.42  &30 &11.20   &2  \\ 
			&50 &106.17  &$<1e^{-4}$  &0.80    &8   &105.31  &11.74  &8  &0.81  &0 \\ \cmidrule(l){1-11}			
			\multirow{3}{*}{a7-84}
			&1 &790.40  &15.37   &3600$^\dagger$  &50  &789.75  &199.19  &50 &0.08   &0 \\ 
			&10 &717.06 &12.69  &3600$^\dagger$   &45  &694.26  &407.79  &45 &3.18   &0  \\ 
			&50 &98.03  &$<1e^{-4}$  &0.53     & 7   &94.97  &13.54  &7  &3.12  &0 \\ \cmidrule(l){1-11}			
			\multirow{3}{*}{a8-96}
			&1  &800.49  &20.45  &3600$^\dagger$  &51  &848.07  &232.45  &53& -5.94  &2 \\ 
			&10 &624.68  &35.94    &3600$^\dagger$  &40  &738.16  &210.54  &47 &-18.17   &7 \\ 
			&50 &155.02  &$<1e^{-4}$  &3.6  &11      &147.75  &29.51  &11  &4.69  &0 \\ 
		\bottomrule
		 $\dagger$: CPLEX timed out
	\end{tabular}
	\label{sas}
}
\end{center}
\end{table}

\subsubsection{Sensitivity Analysis on flat fare $f$}

\begin{table}[H]
\caption{Sensitivity Analysis on $f$ ($p=0.95, s=1$)}\label{safopt}
\begin{center}
\footnotesize
\resizebox{\textwidth}{!}{
	\begin{tabular}{l l l l l l l l l l l }
			\toprule
			\multirow{2}{*}{Instance} &\multirow{2}{*}{$f$} & \multicolumn{4}{l}{CPLEX}& \multicolumn{3}{l}{LS-H} & \multicolumn{2}{l}{Comparison} \\
			\cmidrule(l){3-6}
			\cmidrule(l){7-9}
			\cmidrule(l){10-11} 
		&  & Profit (\$) & Gap (\%) & CPU time (s) & $\sum y_i^\star$ &  Profit (\$) & CPU time (s) & $\sum \hat{y}_i$   &Profit gap (\%)  & $\sum y_i^\star - \sum \hat{y}_i$\\ \cmidrule(l){1-11}
			\multirow{3}{*}{a2-16} 
			&10.0  &83.98  &$<1e^{-4}$  &0.83    &14   &69.73  &5.51  &12 &16.97   &2 \\ 
			&20.0 &134.86  &$<1e^{-4}$  &0.14   &9  &129.98	&3.78  &9 &3.62   &0 \\ 
			&30.0 &\textbf{118.47} &$<1e^{-4}$   &0.05  &5  &\textbf{118.47}  &3.17  &5 &0.00   &0\\ \cmidrule(l){1-11}
			\multirow{3}{*}{a2-20}
			&10.0  &101.94  &$<1e^{-4}$   &0.98  &16   &95.69	  &11.82 &16 &6.13  &0 \\ 
			&20.0 &156.95  &$<1e^{-4}$   &0.08  &10    &148.78 	&3.55 &10 &5.21  &0  \\ 
			&30.0 &\textbf{67.71}  &$<1e^{-4}$   &0.03  &3  &\textbf{67.71}  &2.61 &3 &0.00  &0 \\ \cmidrule(l){1-11}
			\multirow{3}{*}{a2-24} 
			&10.0 &126.62  &$<1e^{-4}$   &18.33  &21  &116.71  &7.04  &21  &7.82   &0\\ 
			&20.0 &244.98  &$<1e^{-4}$   &2.53  &16  &241.87  &4.77  &16 &1.33   &0  \\ 
			&30.0 &119.30  &$<1e^{-4}$     &0.05  &5  &118.38   &3.06  &5 &0.77  &0\\ \cmidrule(l){1-11}
			\multirow{3}{*}{a3-18} 
			&10.0  &110.13  &$<1e^{-4}$   &203.64  &17  &96.83 &13.97 &16 &12.08 &1 \\ 
			&20.0 &157.55  &$<1e^{-4}$   &0.69  &10  &153.68 &3.55   &10 &2.46 &0\\ 
			&30.0  &\textbf{45.17}  &$<1e^{-4}$   &0.05  &2  &\textbf{45.17}  &2.50 &2 &0.00 &0 \\ \cmidrule(l){1-11}
			\multirow{3}{*}{a3-24} 
			&10.0  &109.64  &$<1e^{-4}$   &76.7  &17 &103.84  &8.36  &17 &5.29   &0 \\ 
			&20.0 &151.94  &$<1e^{-4}$   &1.03  &10  &149.17  &4.43  &10  &1.82   &0  \\ 
			&30.0  &120.02  &$<1e^{-4}$   &0.06  &5  &119.42  &3.63  &5  &0.50   &0 \\ \cmidrule(l){1-11}
			\multirow{3}{*}{a3-36} 
			&10.0  &205.54    &$<1e^{-4}$      &2300.64   &32  &194.33 &13.07  &32 &5.46    &0\\ 
			&20.0 &311.95 &$<1e^{-4}$  &18.81 &20  &308.35	&9.38  &20 &1.15    &0 \\ 
			&30.0 &211.58 &$<1e^{-4}$   &0.11  &9  &220.31  &7.42    &9  &0.57    &0\\ \cmidrule(l){1-11}
			\multirow{3}{*}{a4-16} 
			&10.0  &87.05  &$<1e^{-4}$   &55.78  &14 &82.80  &4.41  &14 &4.87   &0 \\ 
			&20.0 &173.95  &$<1e^{-4}$   &3.75  &11  &171.25  &3.34  &11  &1.55   &0 \\ 
			&30.0  &73.93  &$<1e^{-4}$   &0.03 &3  & 71.57	 &2.75  &3  &3.19   &0 \\ \cmidrule(l){1-11}
			\multirow{3}{*}{a4-32} 
			&10.0 &192.64  &11.45   &3600$^\dagger$  &29  &158.47 	&43.33  &25  &17.74   &4 \\ 
			&20.0 &237.22 &$<1e^{-4}$  &32.53 &15  &224.23 	&8.50  &15  &5.48   &0 \\ 
			&30.0 &176.06  &$<1e^{-4}$  &0.13 &7  &143.46  &5.20 &6  &18.52   &1\\ \cmidrule(l){1-11}
			\multirow{3}{*}{a4-40}
			&10.0  &226.60 &9.01   &3600$^\dagger$ &33  &205.84      &65.38  &32 &9.16    &1\\ 
			&20.0 &343.38  &$<1e^{-4}$  &781.50 &21 &326.23  &11.90  &21 &4.99   & 0 \\ 
			&30.0 &145.88   &$<1e^{-4}$  &0.11  &6    &140.64   &6.20  &6  &3.60  & 0\\ \cmidrule(l){1-11}
			\multirow{3}{*}{a4-48}
			&10.0  &272.82  &15.78 &3600$^\dagger$  &40 &214.84  &141.60  &33      &21.25  & 7\\ 
			&20.0 &349.28  &$<1e^{-4}$ &2451.31  &22  &319.20  &28.36  &21   &8.61   & 1 \\ 
			&30.0 &195.49  &$<1e^{-4}$  &0.30  &8   &195.40  &21.80     &8  &0.05  &0 \\ \cmidrule(l){1-11}
			\multirow{3}{*}{a5-50}
			&10.0  &279.11  &18.76  &3600$^\dagger$   &41  &263.72 &95.45  &41   &5.52 &0\\ 
			&20.0 &400.23  &5.63  &3600$^\dagger$  &25  &352.02  &49.16  &55 & 12.04  &3  \\ 
			&30.0 &293.73    &$<1e^{-4}$ &5.06   &12   &291.93  &13.97  &12  &0.62  &0 \\ \cmidrule(l){1-11}
			\multirow{3}{*}{a6-60}
			&10.0  &299.81  &40.70  &3600$^\dagger$   &48  &311.10  &94.42  &51 &-3.76   &-3 \\ 
			&20.0 &607.89  &9.01  &3600$^\dagger$  &37  &580.82  &85.64  &37 &4.45  & 0 \\ 
			&30.0 &323.99  &$<1e^{-4}$  &5.55   &13  &310.52  &20.22  &13  &4.16  &0 \\ \cmidrule(l){1-11}				
			\multirow{3}{*}{a6-72}
			&10.0  &316.71  &30.32  &3600$^\dagger$  &51  &289.98  &243.34  &49 &8.44   &2 \\ 
			&20.0 &598.61 &8.19  &3600$^\dagger$  &37  &544.93  &75.99  &35 &9.02   &2  \\ 
			&30.0 &429.90  &$<1e^{-4}$  &15.77    &17   &371.73  &23.93  &15  &13.53  &2 \\ \cmidrule(l){1-11}			
			\multirow{3}{*}{a7-84}
			&10.0 &290.45  &101.39  &3600$^\dagger$  &46  &437.22  &154.30  &68 &-50.53   &-22 \\ 
			&20.0 &780.98 &16.87  &3600$^\dagger$   &50  &789.75  &199.28  &50 &-1.22   &0  \\ 
			&30.0 &538.99  &$<1e^{-4}$  &2243.91     & 21   &497.94  &64.74  &20  &7.62  &1 \\ \cmidrule(l){1-11}			
			\multirow{3}{*}{a8-96}
			&10.0  &250.36  &194.28 &3600$^\dagger$  &41  &556.59 &235.40  &85& -122.32  &-44 \\ 
			&20.0 &815.36  &18.39    &3600$^\dagger$  &52  &848.07  &228.91  &53 &-4.01   &-1 \\ 
			&30.0 &482.47  &1.74  &3600$^\dagger$  &19    &408.37  &45.52  &17  &15.36 &2 \\ 
		\bottomrule
		
			    
		 $\dagger$: CPLEX timed out
	\end{tabular}
	\label{saf}
}
\end{center}
\end{table}

In Table \ref{saf}, we report the results of the sensitivity analysis on fare $f$, assuming a flat fare structure is imposed where every passenger pays the same fare for the service. The overall average profit gap for Table \ref{saf} is $1.54\%$. From the results we can see that the increase in fare does not yield a monotonic change in profit, indicating that some fare values result in higher profits with a better balance in the trade-off between higher routing cost and more revenue. The relationship between system profit and fare will be further assessed in Section \ref{nyc} under various fare structures. 

The results also indicate that LS-H struggles to converge to optima when a cheap fare is charged. When the fare is too low, revenue generated from servicing more requests can barely cover routing costs, which leads to higher computational costs and larger problem size that CPLEX also struggles to cope with. Furthermore, with a cheap fare, chance constraints (\ref{cbnccc}) are not binding since every passenger is better off not driving (based on utility functions), which makes the CC-DARP more challenging.  This is reflected in both the performance of CPLEX which times out on instance a4-32 for $f=10$ and for several larger instances too. In comparison, while the proposed heuristic yields occasionally large profit gaps on low fare instances, it also substantially outperform the best integer solution found by CPLEX on large instances, e.g. a7-84 and a8-96. Note that the last column of Table \ref{saf} only represents the difference between $\sum \hat{y}_i$ and the best lower bound returned by CPLEX $\sum y_i^\star$ within time limit, which is not necessarily the optimal user selection when the instances are too large. It is also worth mentioning that evidently it is more difficult for heuristics to jump out of local optima when tight time windows are imposed \citep{feillet2005traveling}. \\

Overall, the sensitivity analysis results (Tables \ref{sap} - \ref{saf}) show that LS-H can solve instances up to 96 users within reasonable CPU times with an overall average optimality gap of $2.69\%$, and is particularly efficient when chance constraints are binding. 

\subsection{Experiments on NYC Data} \label{nyc}

In this section, we explore the behaviour of the class-based CC-DARP model on real data extracted from NYC yellow taxi dataset. We categorise all 143 requests into two classes based on their drop-off locations: Class $M_1$ if the trip destination is in Zone 1 and Class $M_2$ if trip destination is in Zone 2 (as defined in Figure \ref{mht}). According to \cite{he2020evaluation}, the average Value of Time (VoT) in Manhattan is $\$29$/h.  \cite{small1982scheduling} pointed out that the parameter of schedule delay is 1.5-2 times of travel time. Accordingly, we define Class $M_1$ as the base class of Manhattan section with $\beta_T^{M_1}=\$29$/h and  $\beta_S^{M_1}=\$58$/h. For Class $M_2$, we set $\beta_T^{M_2} = 0.9\beta_T^{M_1}$ and $\beta_S^{M_2}=0.9\beta_S^{M_1}$, as we assume a slightly lower VoT for outside downtown area. Note that this is purely for the purpose of distinguishing two user classes. \\

We query travel time for all 143 requests (286 nodes) to obtain $t_{ij}  (\forall i \in \N,\forall j \in \N)$ from Google Maps API during the morning peak (6am-9am) on a weekday. The travel cost $c_{ij} (\forall i \in N,\forall j \in N)$ is then set proportional to $t_{ij}$:  $ c_{ij} = 0.1 t_{ij} (\forall i \in \N,\forall j \in \N) $. We then obtain our test instance named ``NYC'' with $n=143$ and $|\N|=288$ by adding two artificial depot nodes $\{0,2n+1\}$. Travel cost from depot to all nodes is set to $\$2.5$ as a remedy for not having a fixed cost penalty for dispatching an additional vehicle in the objective function. For each request $i \in \P=\{1,...,143\}$, we define the trip to be an inbound trip ($i \in \IB$) if $i$ is an odd number, outbound trip ($i \in \OB$) otherwise. For each inbound/outbound trip, we set a 15-minute time window on its departure/arrival time, with its recorded pickup/drop-off time as the midpoint of the time window. For the record, the cost of transit in NYC is \$2.75 per trip. The average taxi fare of morning peak on the sample day is \$12.93.\\

\subsubsection{Flat Fare Structure}

\begin{figure}[H]
\centering
\begin{subfigure}{0.5\linewidth}{\includegraphics[width=0.95\textwidth,valign=c]{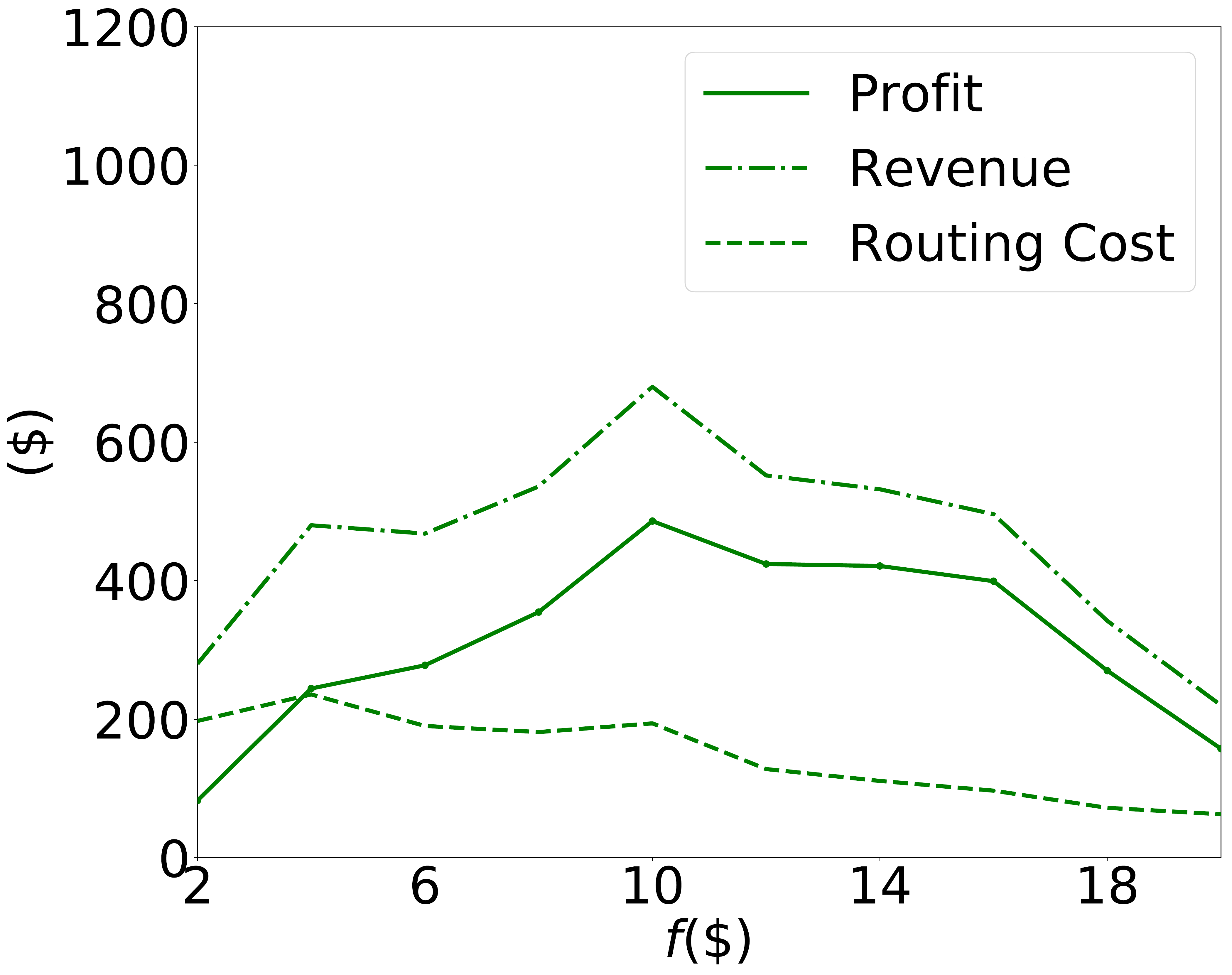}\caption{Profit, revenue and routing costs - $f$.\label{ffp}}}
\end{subfigure}\\
\begin{subfigure}{0.9\linewidth}{\includegraphics[width=0.95\textwidth,valign=c]{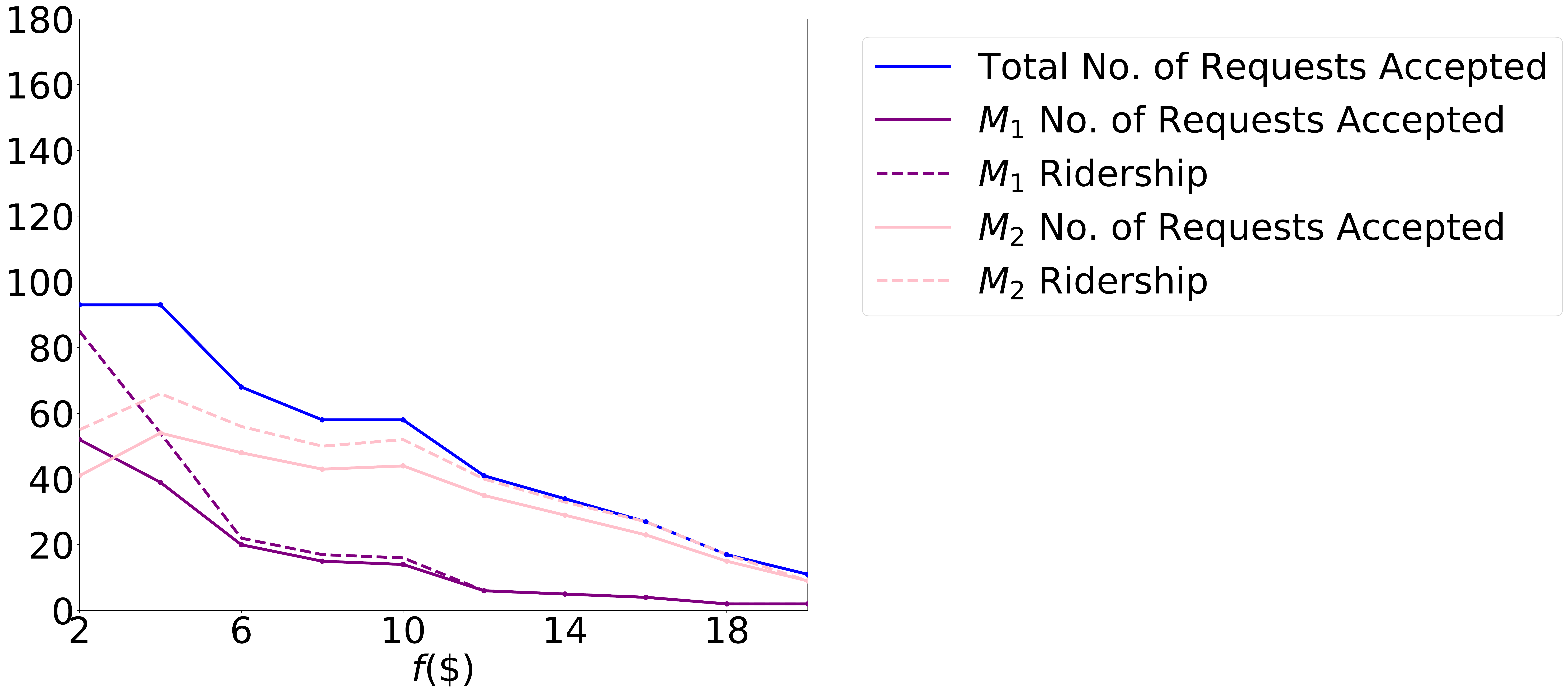}\caption{Number of accepted requests - $f$.\label{ffr}}}
\end{subfigure}
\caption{Results from the NYC instance under a flat fare structure $ f_{im}=f$ for all requests $i \in \P$ and user class $m \in \M$.}
\label{ff}
\end{figure}

We first impose a flat fare structure for DRT services. A flat fare structure charges a single fare for all passengers regardless of their travel times or travel distances. When a flat fare structure is applied, the base fare level plays a crucial rule in ridership and system profit. Figure \ref{ffp} shows how profit, revenue and routing cost of the system change with an increasing fare $f$. As $f$ increases, system profit does not increase or decrease monotonically. The optimal fare level for maximum profit lies around the peak of the profit curve. Figure \ref{ffr} demonstrates how ridership falls in both classes. Note that in our model, the fare is charged per person instead of per request. As observed from Figure \ref{ffr}, passengers with higher VoT (Class $M_1$) tend to have a higher fare elasticity, as ridership in $M_1$ drops more sharply than ridership in $M_2$. This indicates that the bindingness of chance constraint has a more significant effect on the user group with higher VoT under flat fare structures. 

\begin{figure}[H]
\centering
\begin{subfigure}{0.48\linewidth}{\includegraphics[width=0.95\textwidth,valign=c]{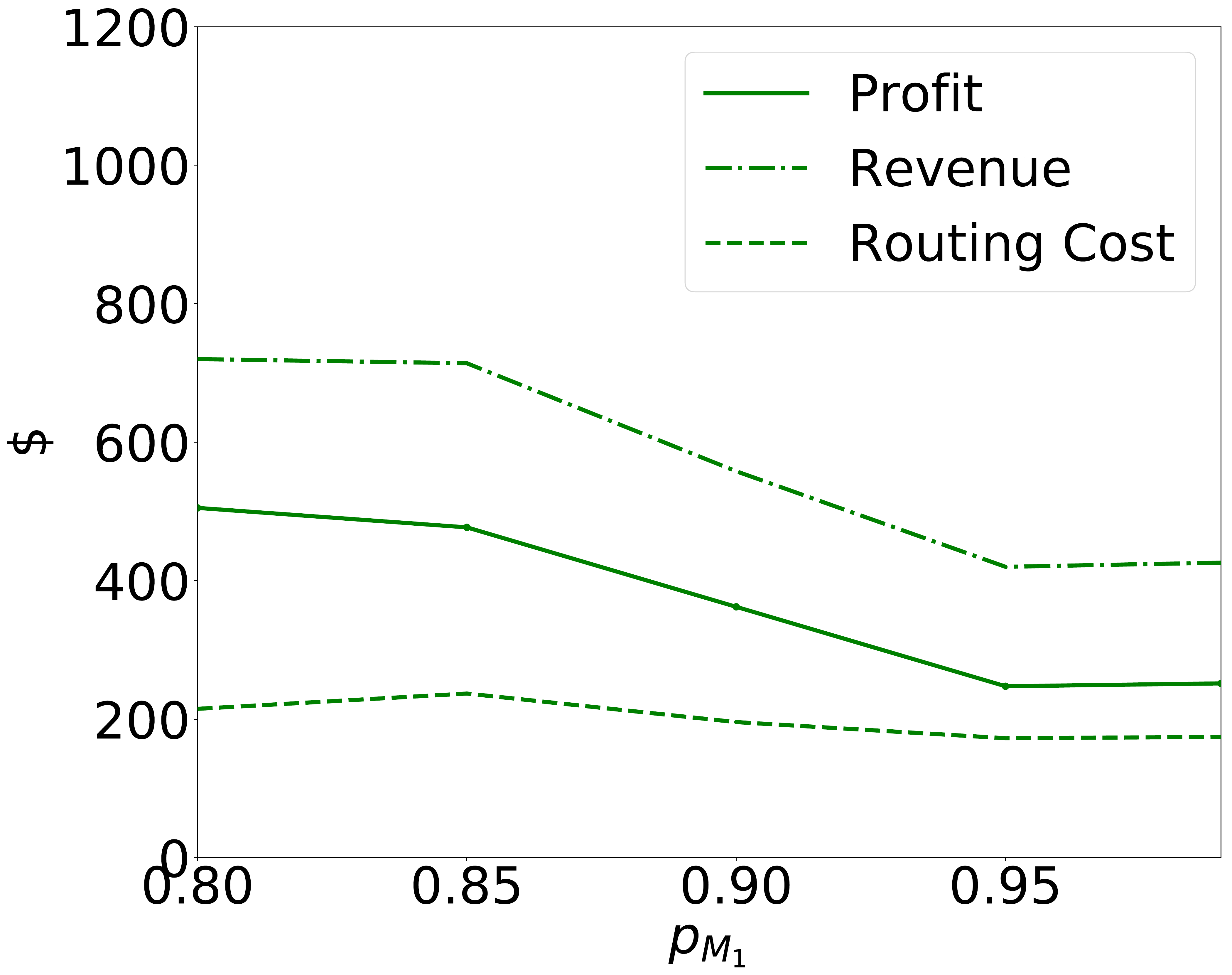}\caption{Profit, revenue and routing costs - $p_{M1}$.\label{fp1p}}}
\end{subfigure}
\begin{subfigure}{0.48\linewidth}{\includegraphics[width=0.95\textwidth,valign=c]{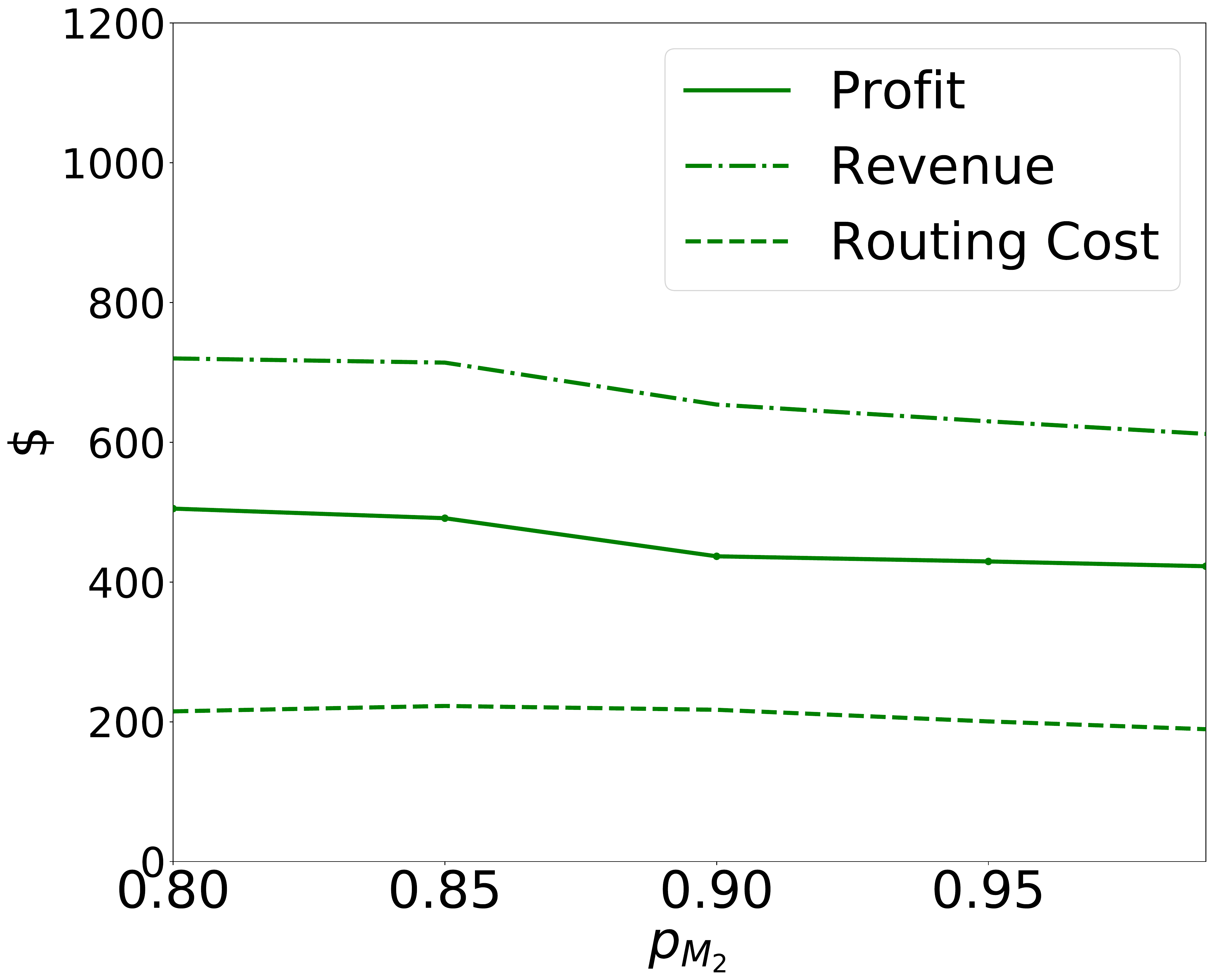}\caption{Profit, revenue and routing costs - $p_{M2}$.\label{fp2p}}}
\end{subfigure}\\
\begin{subfigure}{0.9\linewidth}{\includegraphics[width=0.95\textwidth]{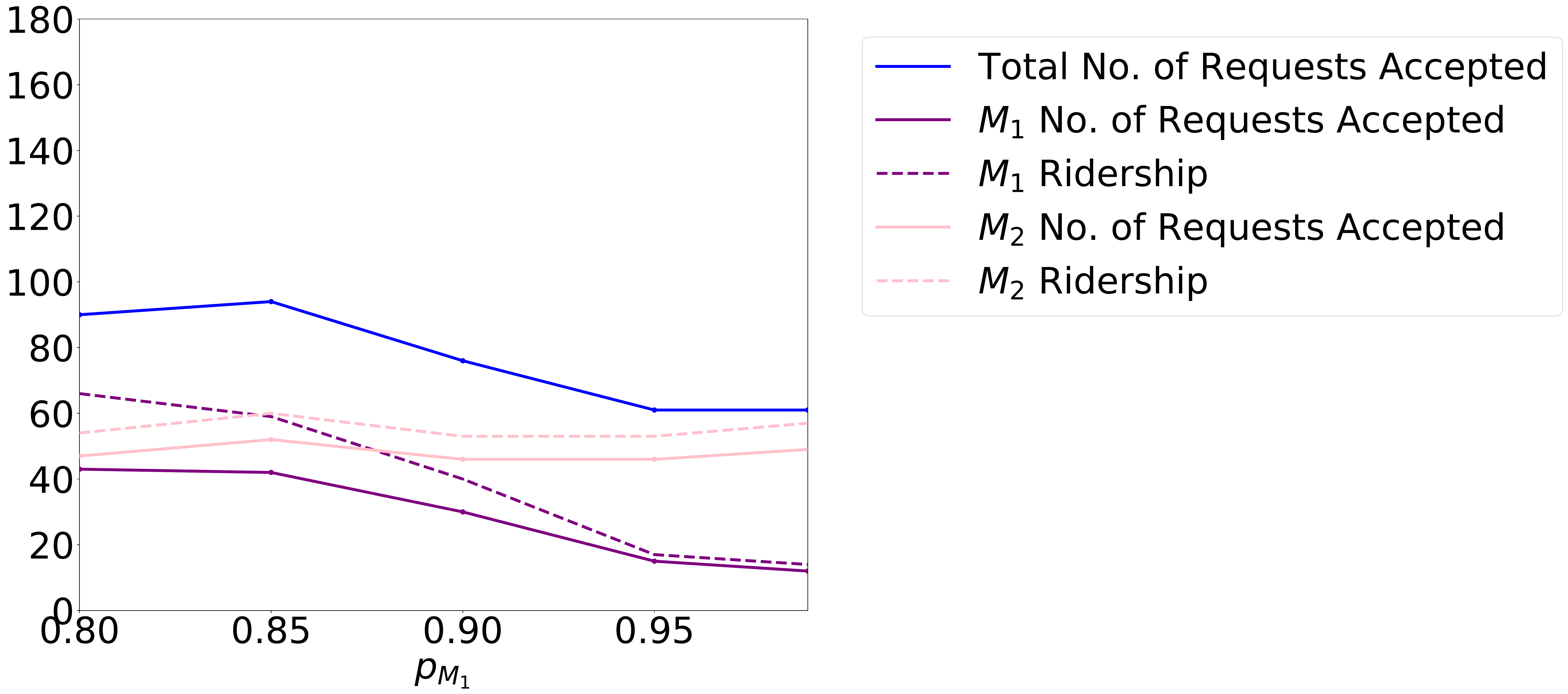}\caption{Number of accepted requests - $p_{M1}$.\label{fp1r}}}
\end{subfigure}
\begin{subfigure}{0.9\linewidth}{\includegraphics[width=0.95\textwidth]{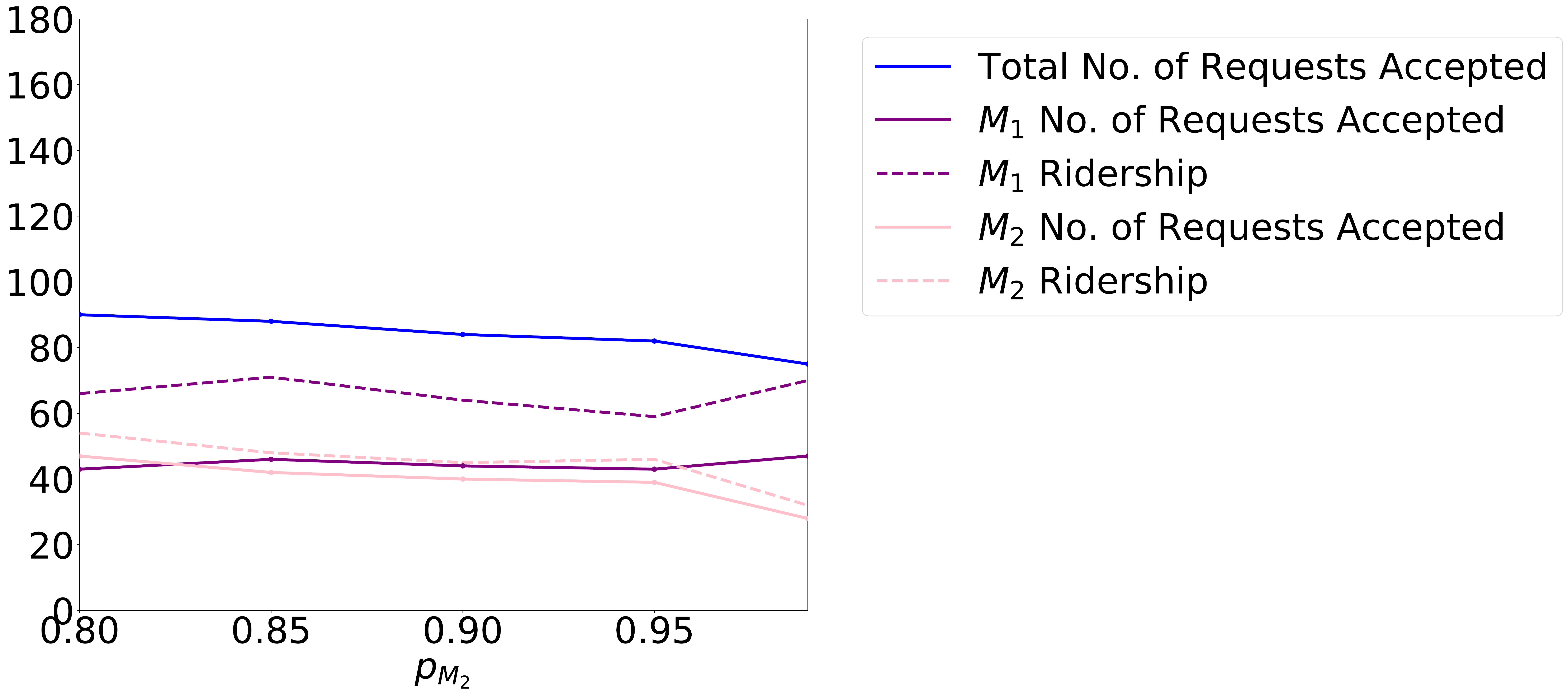}\caption{Number of accepted requests - $p_{M2}$.\label{fp2r}}}
\end{subfigure}\\
\caption{Results from the NYC instance when changing $p_{M1}$ and $p_{M2}$ under a flat fare structure. }
\label{fp}
\end{figure}
Figure \ref{fp} displays the effect of varying $p_{M_1}$ and $p_{M_2}$ in the range of $(0.80, 0.99)$ as $f$ is set to \$6.0. While profit falls as we tighten up $p_{M_1}$ and $p_{M_2}$ (Figure \ref{fp1p} and \ref{fp2p}), Figure \ref{fp1r} and \ref{fp2r} indicate that an increase in $p$ tend to cause a sharper drop in ridership in the group with higher VoT, as more requests in Class $M_1$ are filtered out by chance constraints when $p_{M_1}$ increases. Note that when increasing $p$ in one class, the ridership in the other class increases as declining more requests in one class enables capacitated vehicles to accept more requests from the other class.

\begin{figure}[H]
\centering
\begin{subfigure}{0.9\linewidth}{\includegraphics[width=0.95\textwidth,valign=c]{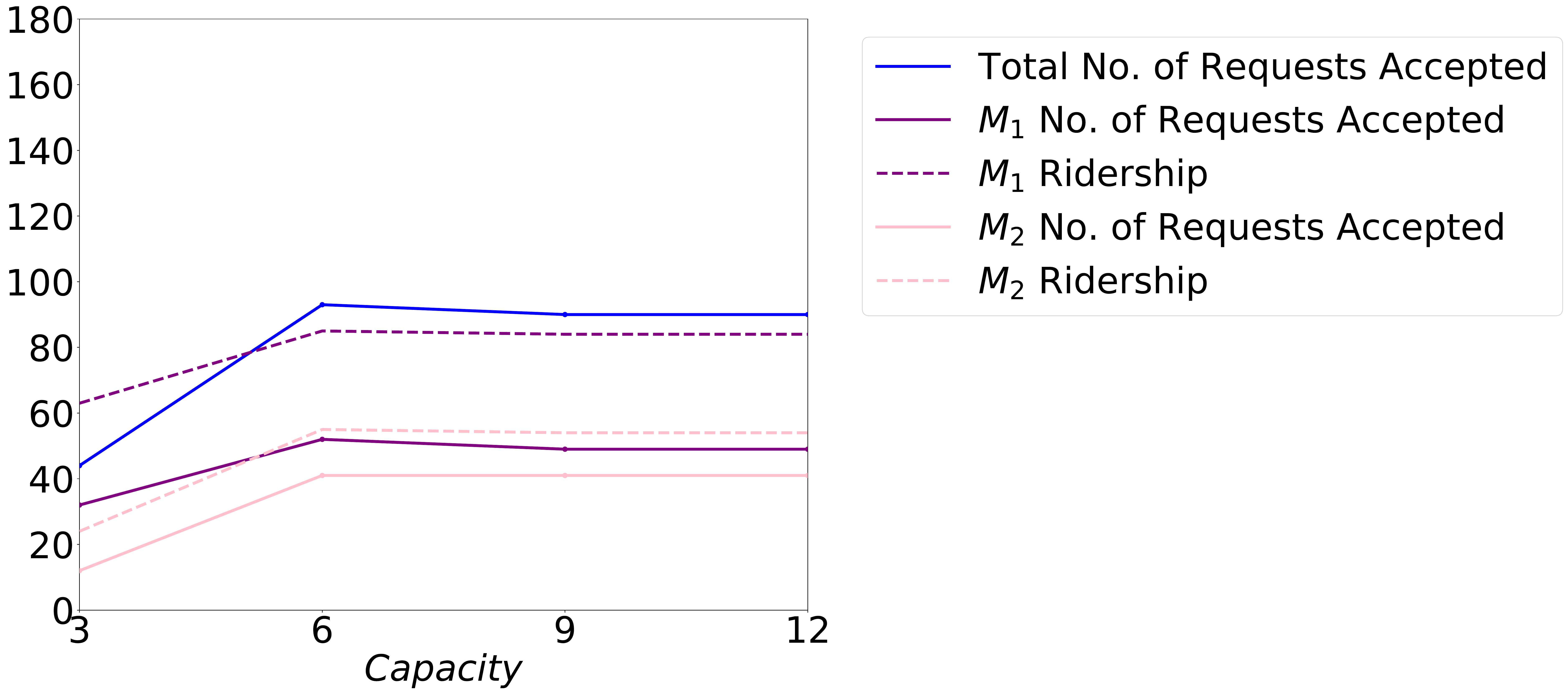}\caption{Number of accepted requests - Capacity.\label{fcr}}}
\end{subfigure}
\begin{subfigure}{0.5\linewidth}{\includegraphics[width=0.95\textwidth,valign=c]{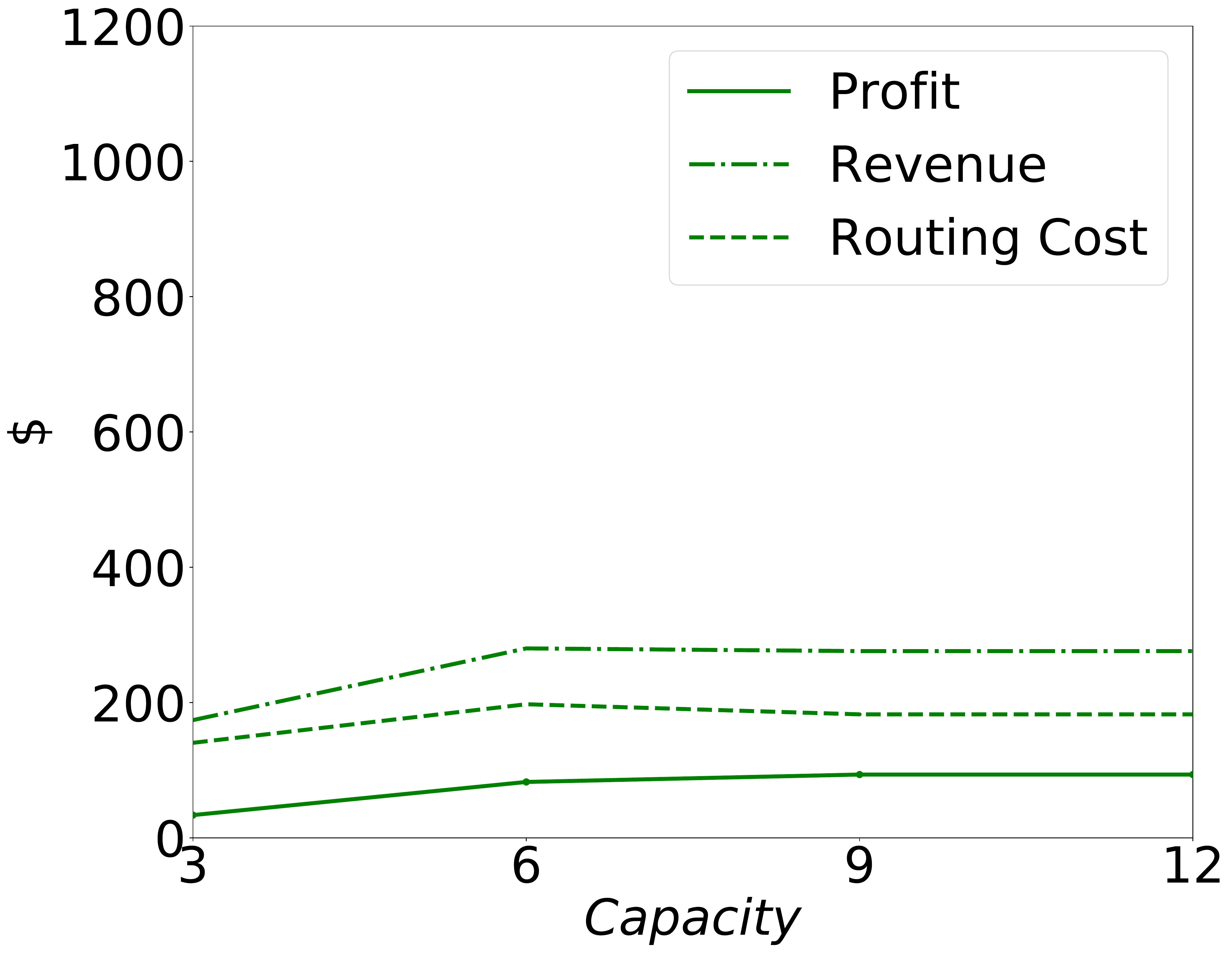}\caption{Profit, revenue and routing costs - Capacity.\label{fcp}}}
\end{subfigure}
\begin{subfigure}{0.45\linewidth}{\includegraphics[width=0.95\textwidth,valign=c]{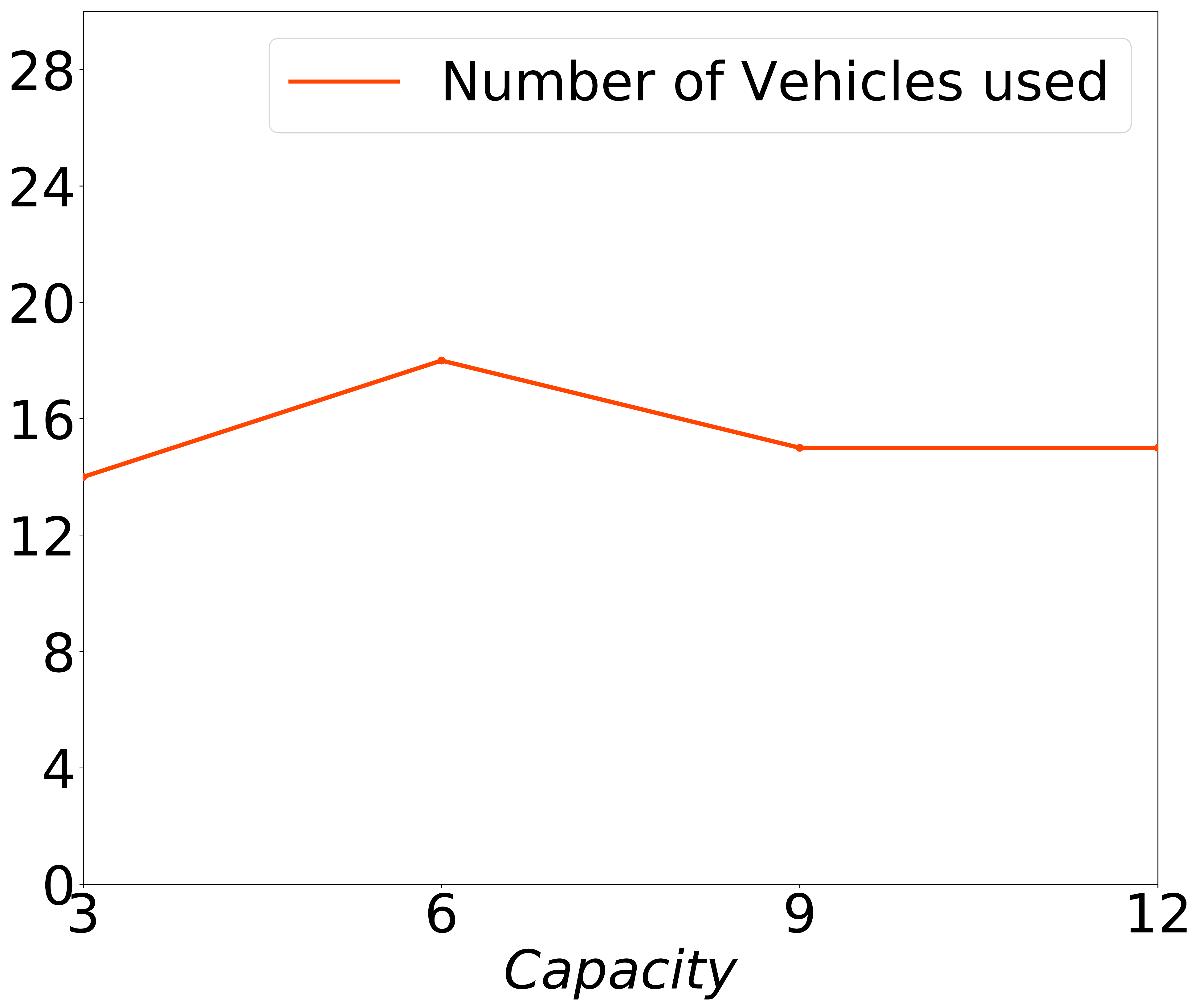}\caption{Fleet size - Capacity.\label{fcn}}}
\end{subfigure}
\caption{ Results from the NYC instance when changing $Capacity$ under a flat fare structure. }
\label{fc}
\end{figure}

Capacity of vehicles could also affect the performance of the system. In this study, we assume that all vehicles are homogeneous, and passengers from the same request are not allowed to be separated and served by different vehicles. Figure \ref{fc} reports changes in profit, ridership and number of vehicles used when we increase the capacity of the fleet while holding $p$ and $f$ constant. As shown in Figure \ref{freq}, the size of user requests varies from 1 to 6. Inevitably, using vehicles with capacities smaller than 6 will automatically filter out some requests that exceed the capacity, which leads to lower profit and ridership when capacity is low. If vehicle capacity is greater than 6, all requests are eligible for DRT services. Using vehicles with higher capacity could possibly reduce fleet size. However, when the time windows are tight, fleet size might not vary significantly even with bigger vehicles.\\

\subsubsection{Distance-based Fare Structure}

\begin{figure}[H]
\centering
\begin{subfigure}{0.48\linewidth}{\includegraphics[width=0.95\textwidth,valign=c]{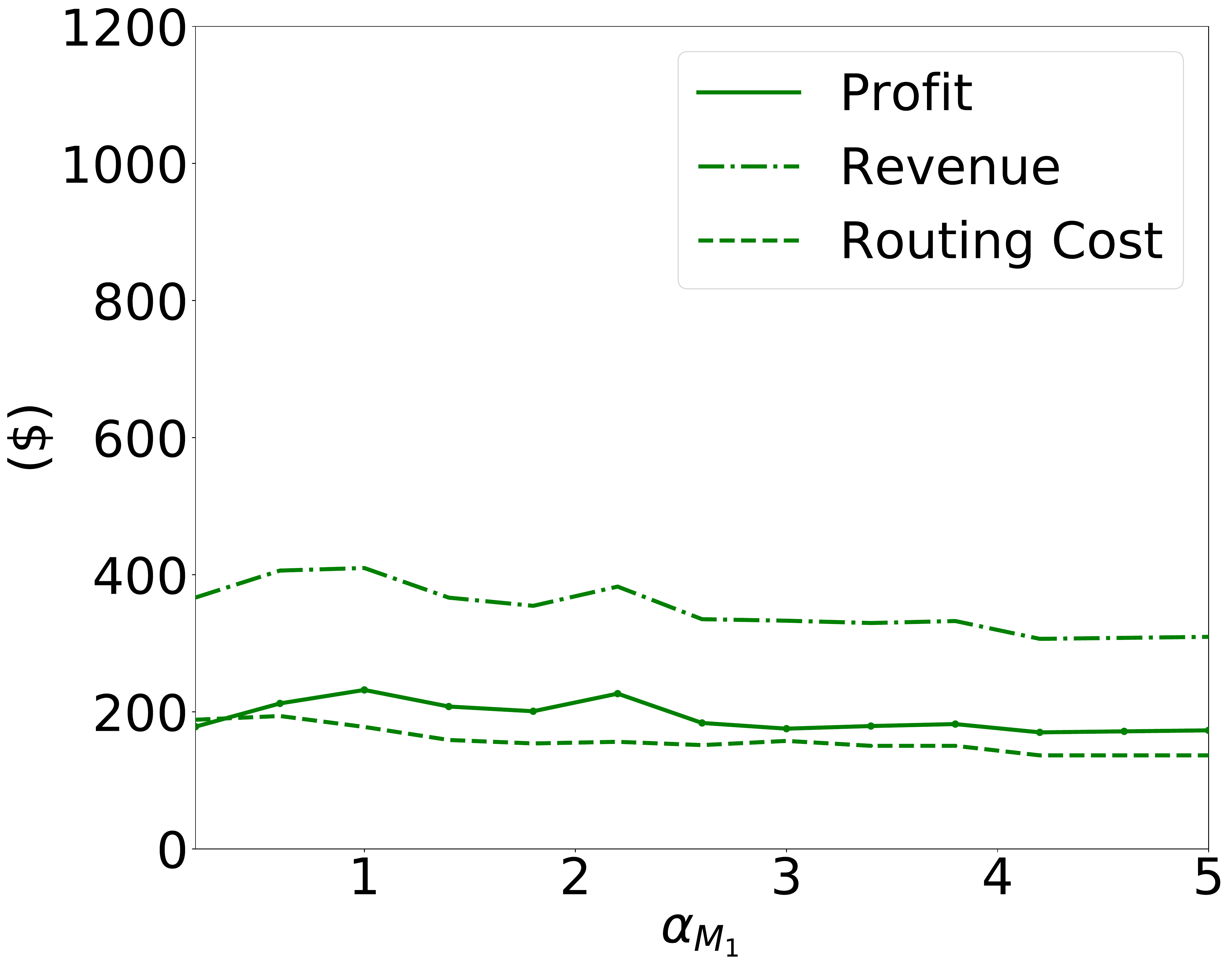}\caption{Profit, revenue and routing costs - $\alpha_{M1}$.\label{df1p}}}
\end{subfigure}
\begin{subfigure}{0.48\linewidth}{\includegraphics[width=0.95\textwidth,valign=c]{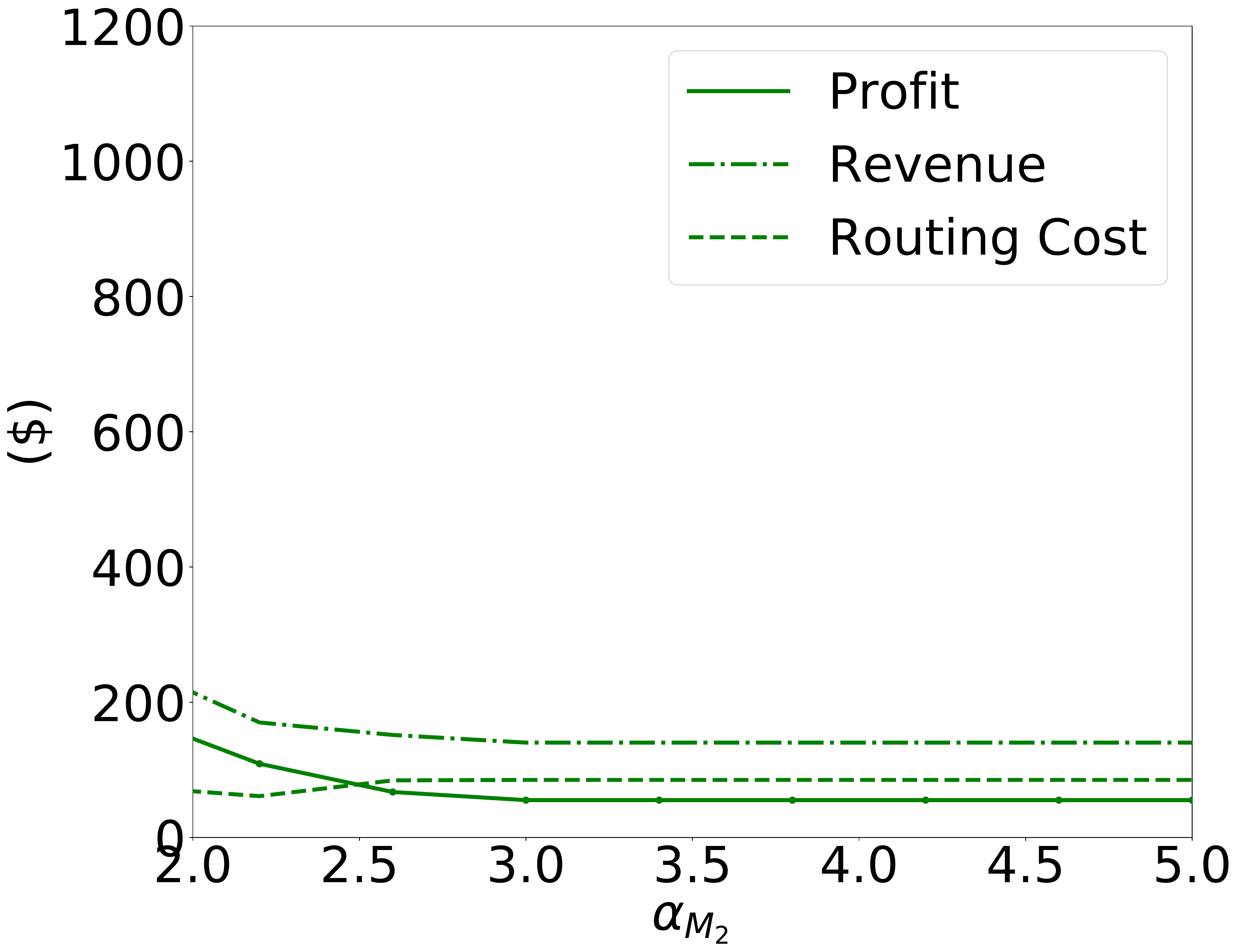}\caption{Profit, revenue and routing costs - $\alpha_{M2}$.\label{df2p}}}
\end{subfigure}\\
\begin{subfigure}{0.9\linewidth}{\includegraphics[width=0.95\textwidth,valign=c]{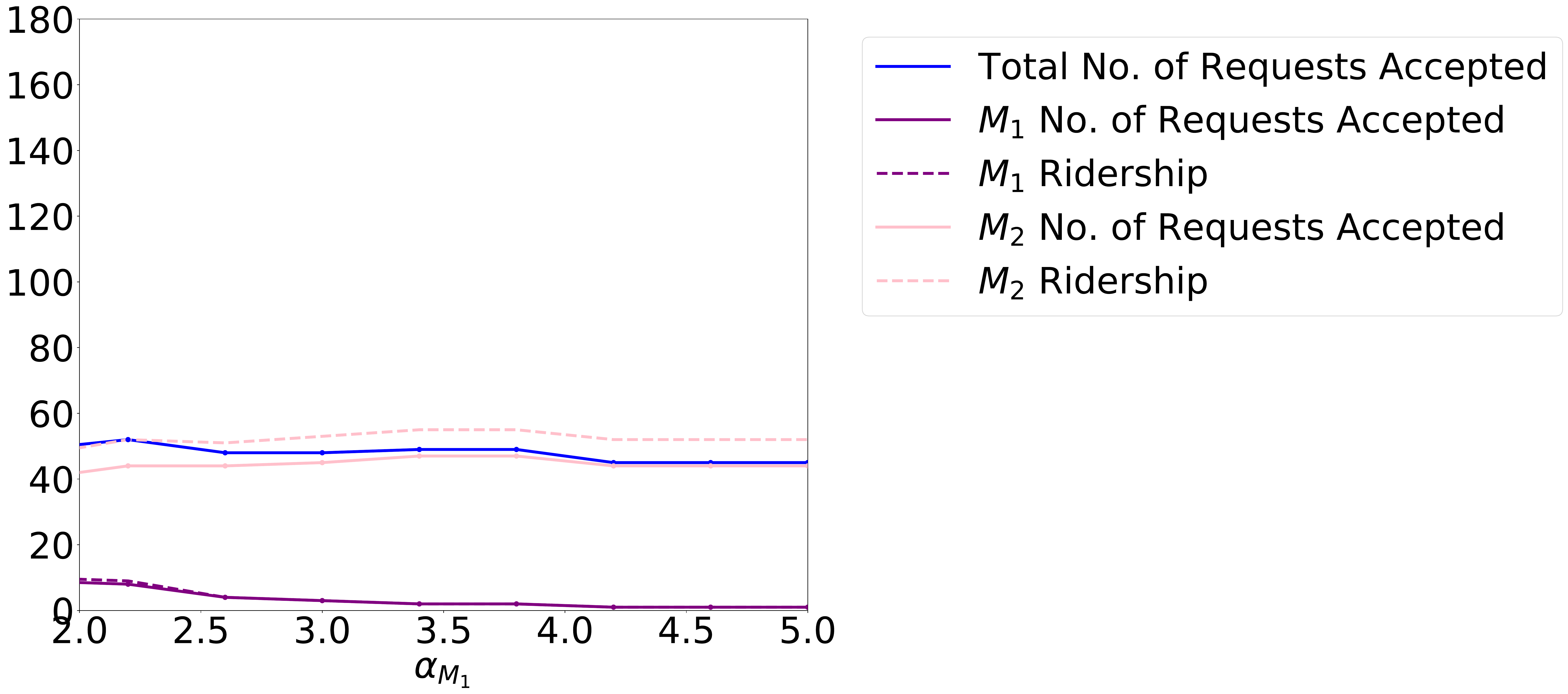}\caption{Number of accepted requests - $\alpha_{M1}$.\label{df1r}}}
\end{subfigure}\\
\begin{subfigure}{0.9\linewidth}{\includegraphics[width=0.95\textwidth,valign=c]{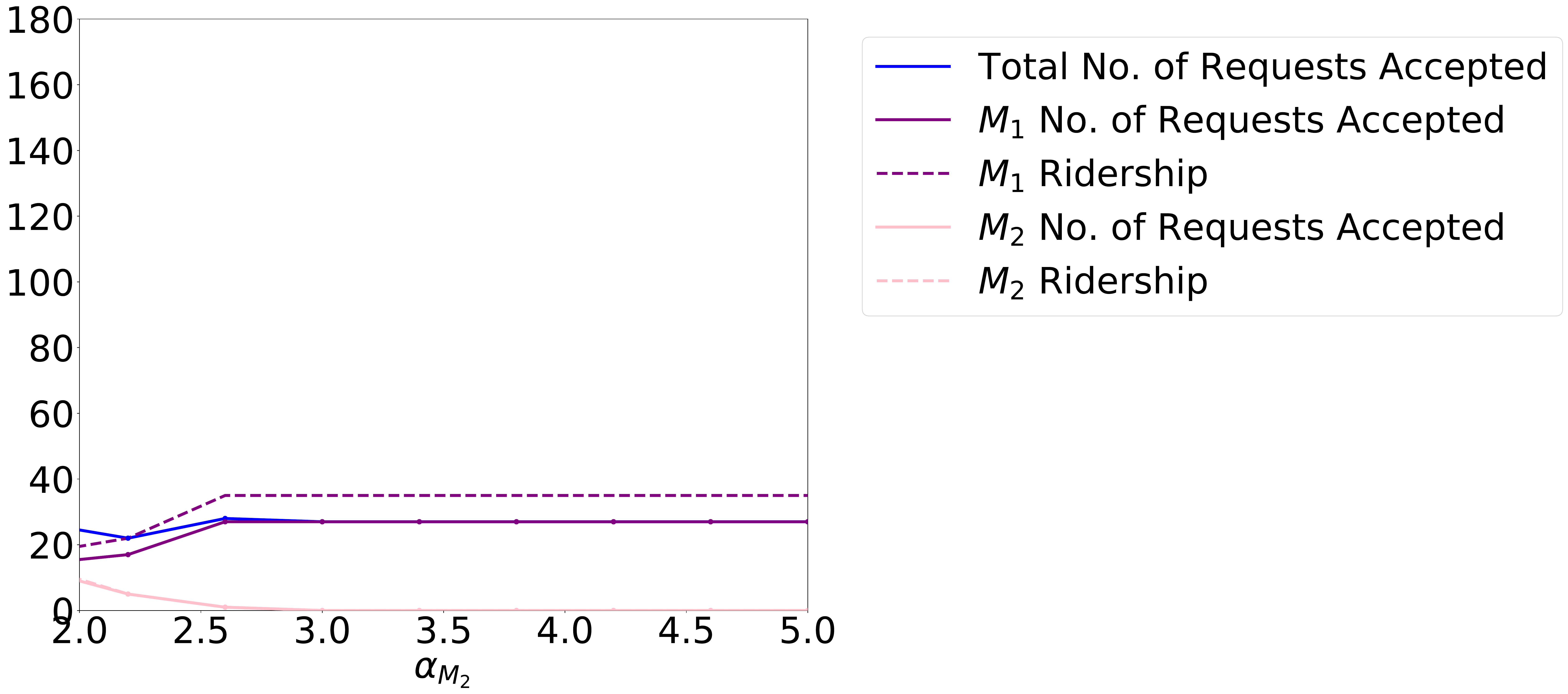}\caption{Number of accepted requests - $\alpha_{M2}$.\label{df2r}}}
\end{subfigure}\\
\caption{Results from the NYC instance when changing $\alpha_{M1}$ and $\alpha_{M2}$ under a distance-based fare structure.}
\label{df}
\end{figure}

\begin{figure}[H]
\centering
\begin{subfigure}{0.48\linewidth}{\includegraphics[width=0.95\textwidth,valign=c]{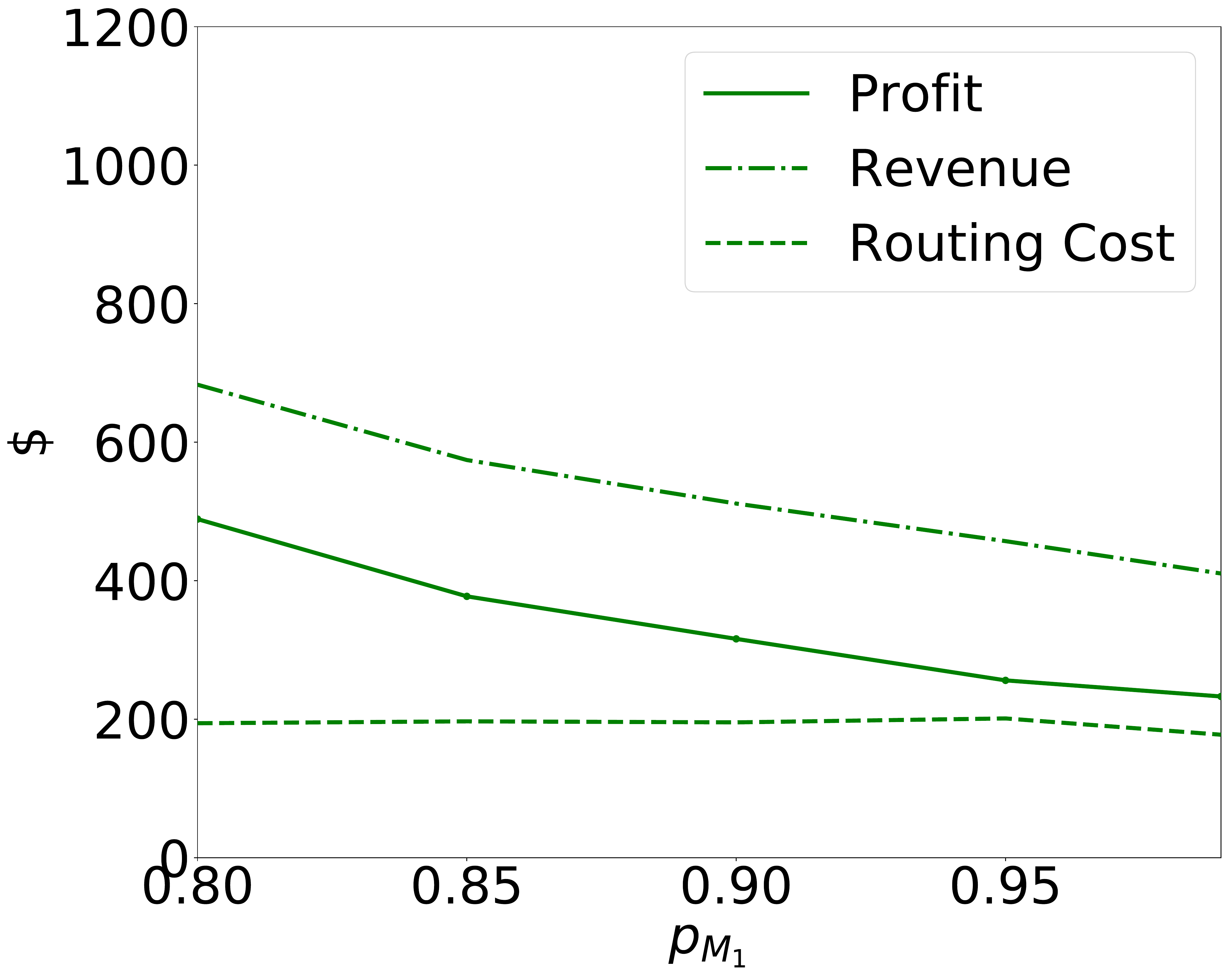}\caption{Profit, revenue and routing costs - $p_{M1}$.\label{dp1p}}}
\end{subfigure}
\begin{subfigure}{0.48\linewidth}{\includegraphics[width=0.95\textwidth,valign=c]{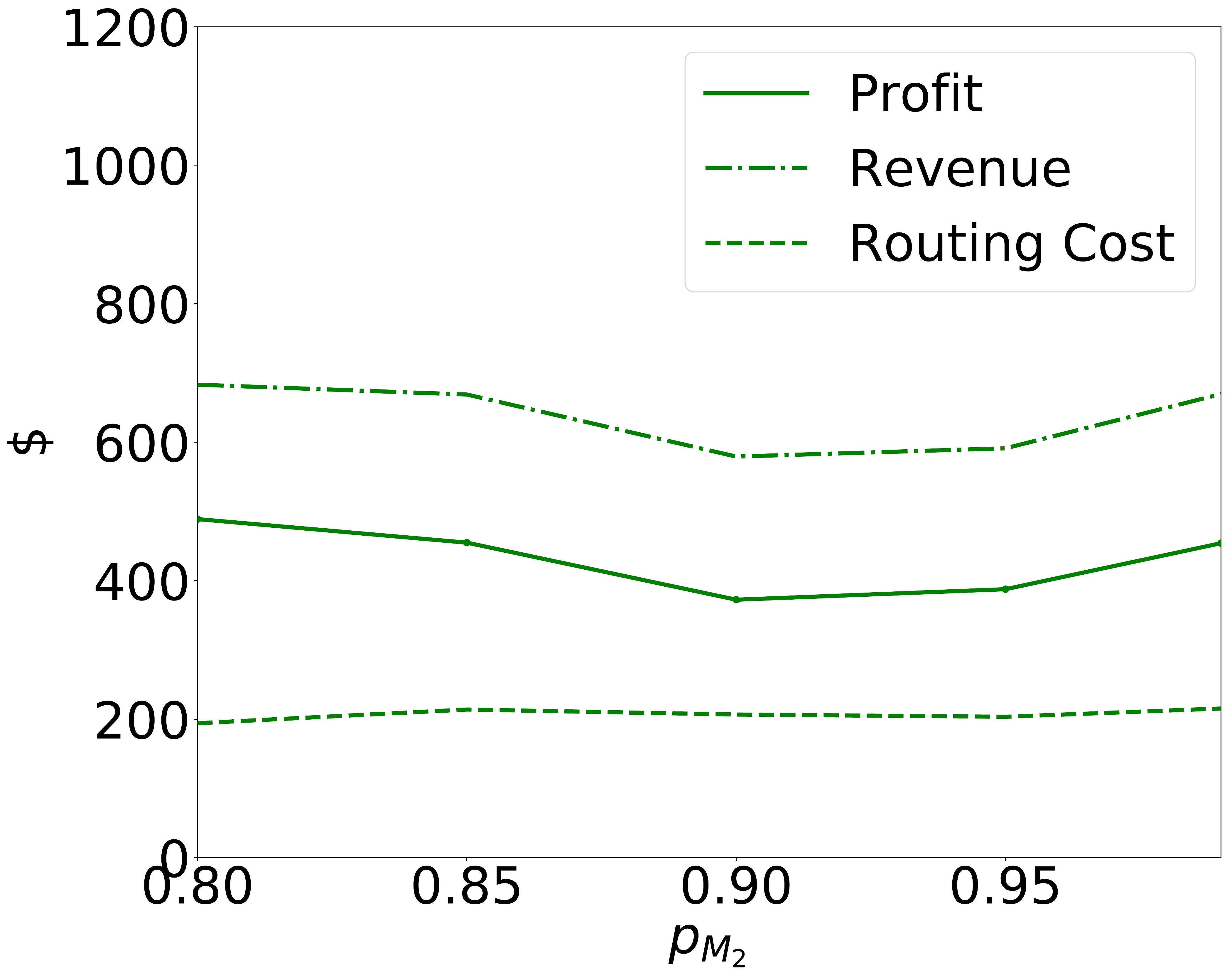}\caption{Profit, revenue and routing costs - $p_{M2}$.\label{dp2p}}}
\end{subfigure}\\
\begin{subfigure}{0.9\linewidth}{\includegraphics[width=0.95\textwidth,valign=c]{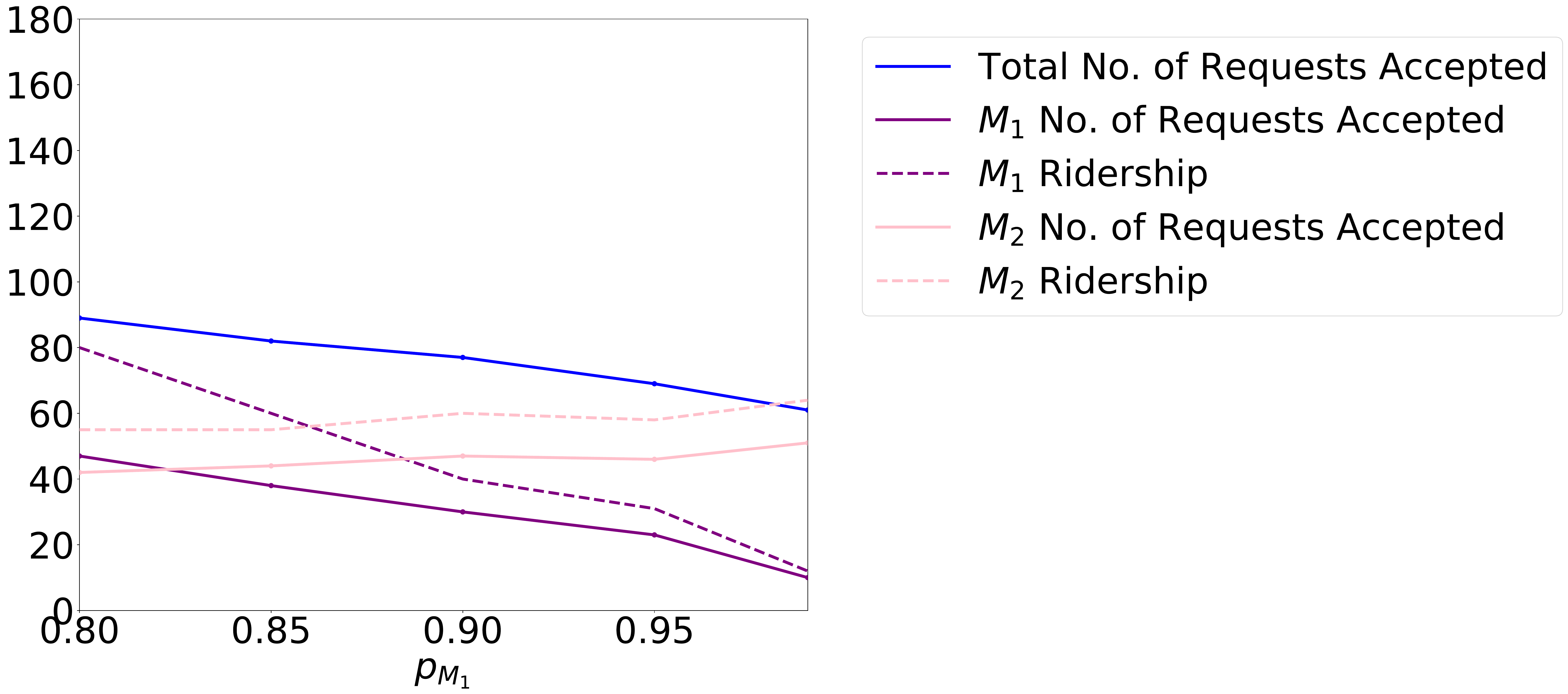}\caption{Number of accepted requests - $p_{M1}$.\label{dp1r}}}
\end{subfigure}\\
\begin{subfigure}{0.9\linewidth}{\includegraphics[width=0.95\textwidth,valign=c]{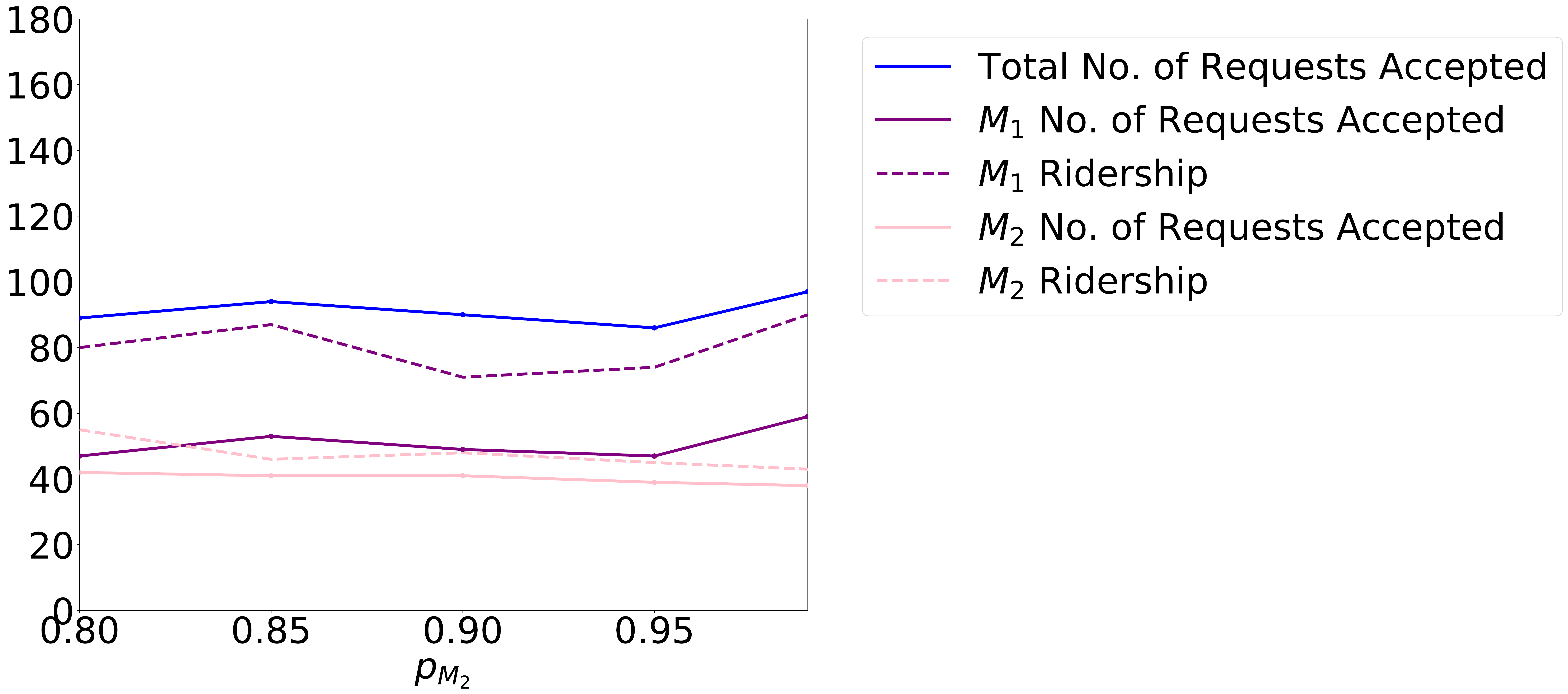}\caption{Number of accepted requests - $p_{M2}$.\label{dp2r}}}
\end{subfigure}\\
\caption{Results from the NYC instance when changing $p_{M1}$ and $p_{M2}$ under a distance-based fare structure.}
\label{dp}
\end{figure}

Distance-based fare structure is another widely adapted fare structure in transit systems, where a unique fare is established for each station or stop pair. For the NYC instance, we define the distance based fare for each class as:
\[
f_{im} = \alpha_m t_{i,n+i}, \qquad	\forall m \in \M, \forall i \in \PM.
\]

Figure \ref{df} shows the behaviours of profit and ridership of the system if we increase $\alpha_m $ in one class while holding the other one unchanged. Comparing Figures \ref{df1p} and \ref{df1r} with Figures \ref{df2p} and \ref{df2r}, it is observed that increases in $\alpha_{M_2}$ cause more drastic changes in both profit and ridership, indicating that adjusting the distance-based fare parameter of the lower VoT class has a bigger impact on the system. On the other hand, Figure \ref{dp} displays the results of the influence of increasing $p_{M_1}$ and $p_{M_2}$, respectively. Figures \ref{dp2p} and \ref{dp2r} show that tightening the chance constraints on $M_2$ may result in an increase in the overall profit and ridership. When the fleet size or the capacity of vehicles is limited, tightening the chance constraints in one class is equivalent as putting more weight on the other class(es) that are potentially more profitable.

\begin{figure}[H]
\centering
\begin{subfigure}{0.9\linewidth}{\includegraphics[width=0.95\textwidth,valign=c]{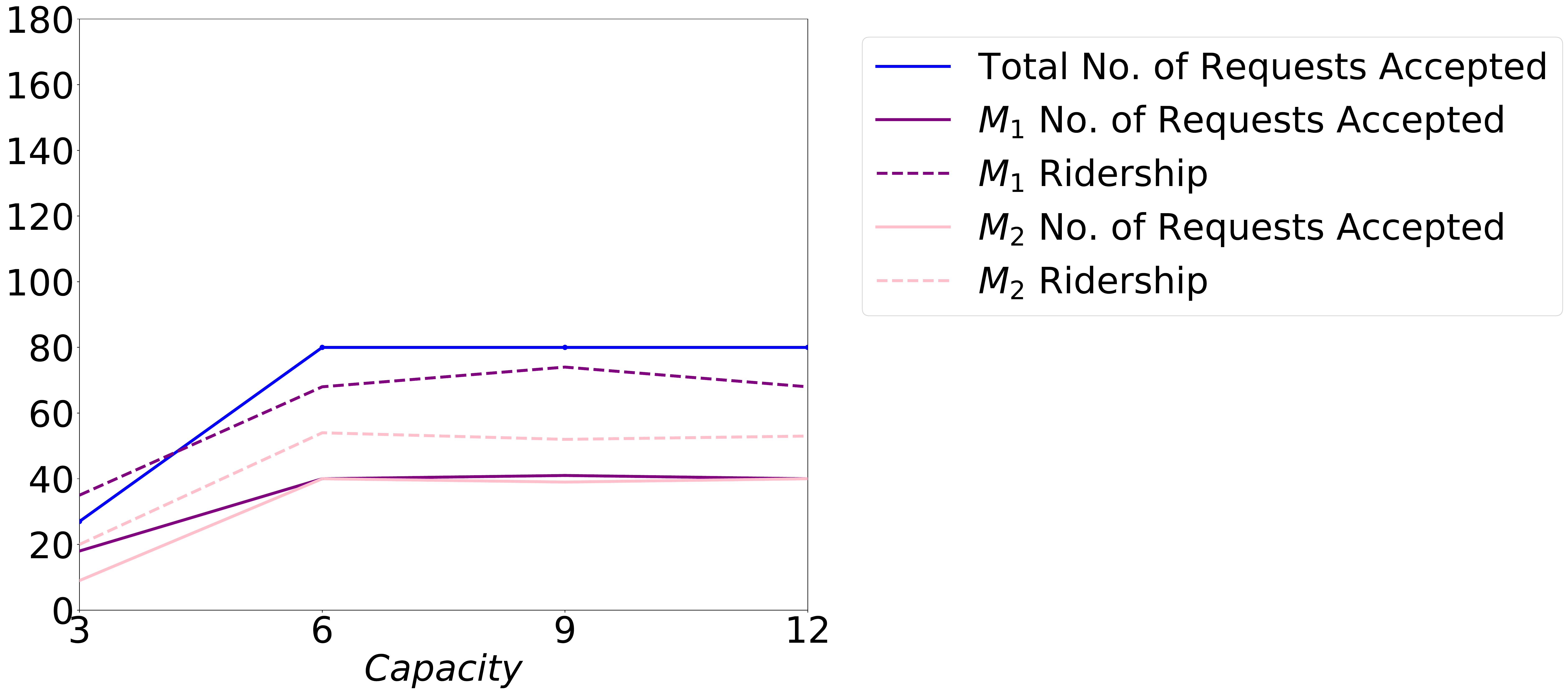}\caption{Number of accepted requests - capacity.\label{dcr}}}
\end{subfigure}
\begin{subfigure}{0.5\linewidth}{\includegraphics[width=0.95\textwidth,valign=c]{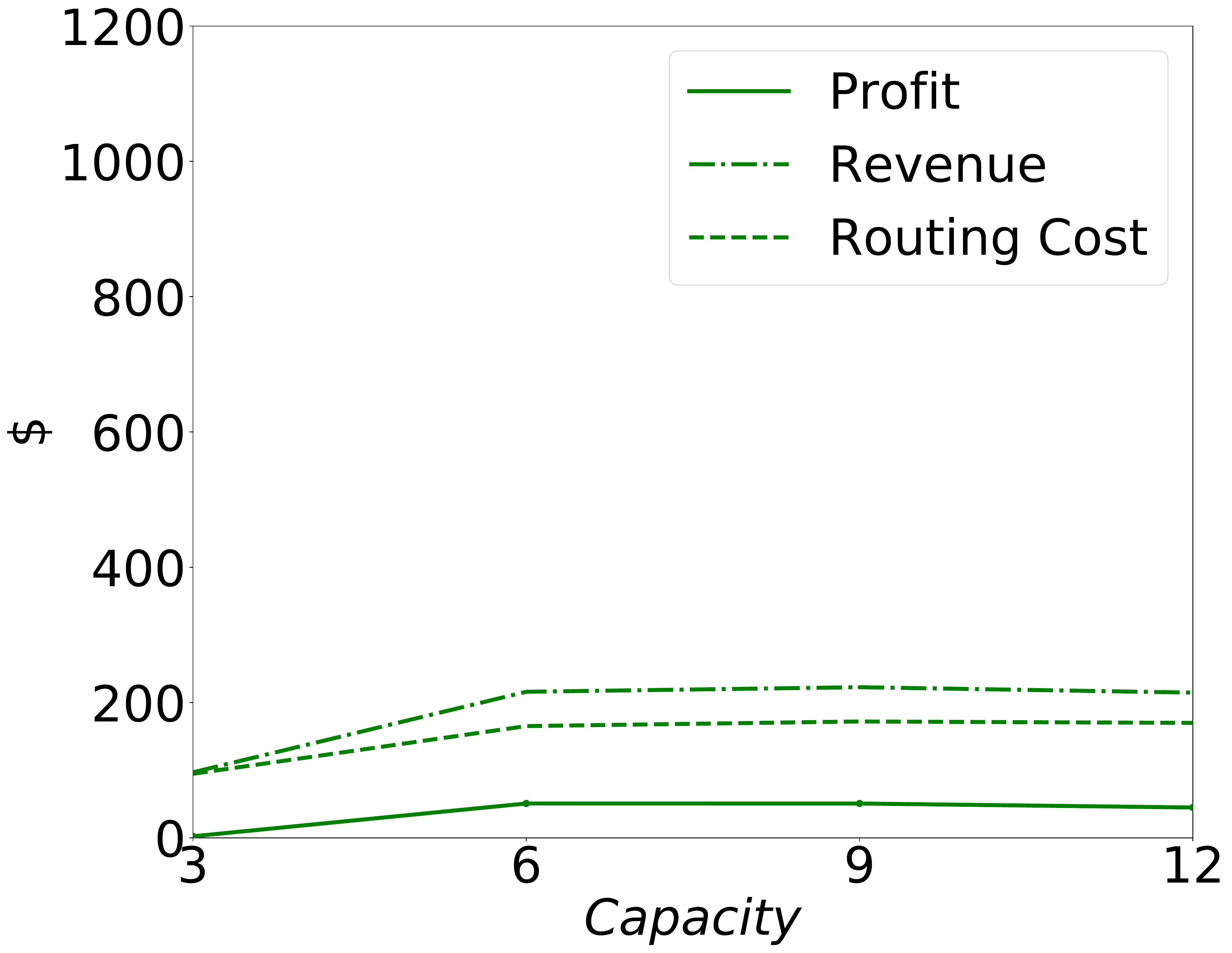}\caption{Profit, revenue and routing costs - capacity.\label{dcp}}}
\end{subfigure}
\begin{subfigure}{0.45\linewidth}{\includegraphics[width=0.95\textwidth,valign=c]{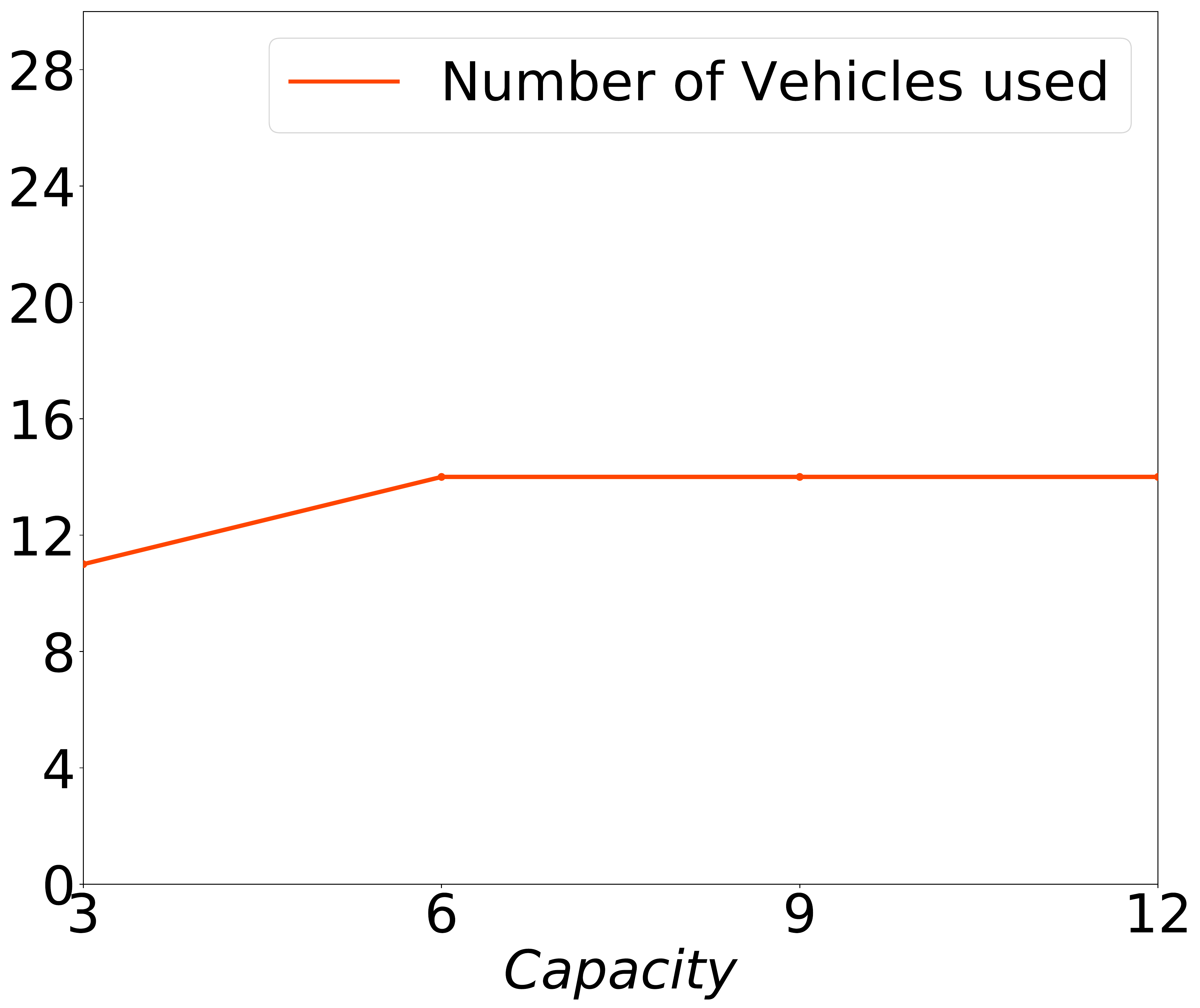}\caption{Fleet size - capacity.\label{dcn}}}
\end{subfigure}
\caption{Results from the NYC instance when changing the vehicle capacity $Q$ under a distance-based fare structure.}
\label{dc}
\end{figure}

\subsubsection{Zone-based Fare Structure}
A zone-based or zonal fare structure determine fare based on the number of zones each passenger travels through. In this section, we make some modifications on the zone-based fare structure, and simply define the fare based on the destination zone of the trip, namely $f_{M_1}$ and $f_{M_2}$. This simplified zone-based pricing rule is inspired by the implementation of congestion pricing policy to charge vehicles entering Manhattan under 60th Street (Zone 2) \citep{he2020validated}.\\

\begin{figure}[H]
\centering
\begin{subfigure}{0.48\linewidth}{\includegraphics[width=0.95\textwidth,valign=c]{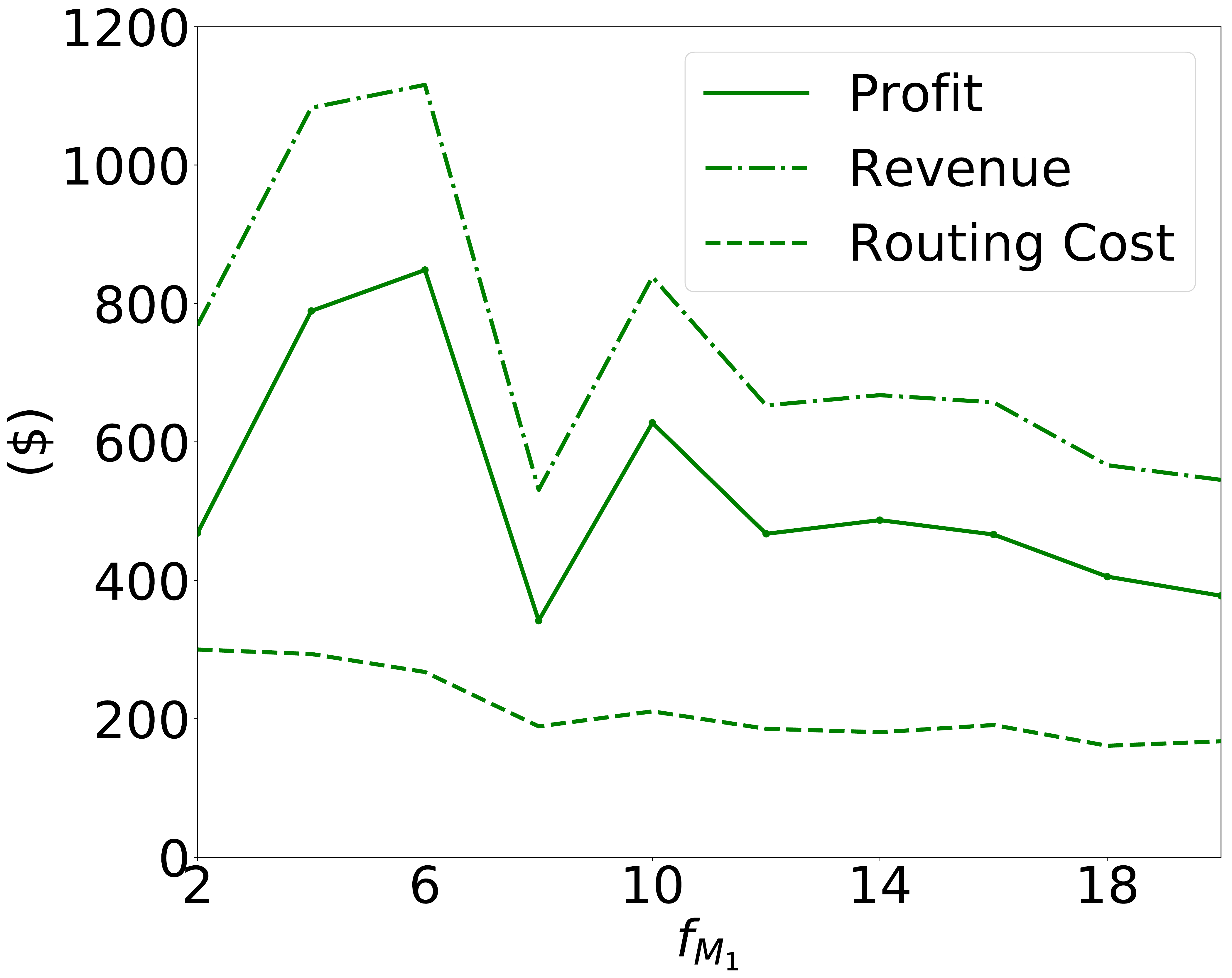}\caption{Profit, revenue and routing costs - $f_{m1}$.\label{zf1p}}}
\end{subfigure}
\begin{subfigure}{0.48\linewidth}{\includegraphics[width=0.95\textwidth,valign=c]{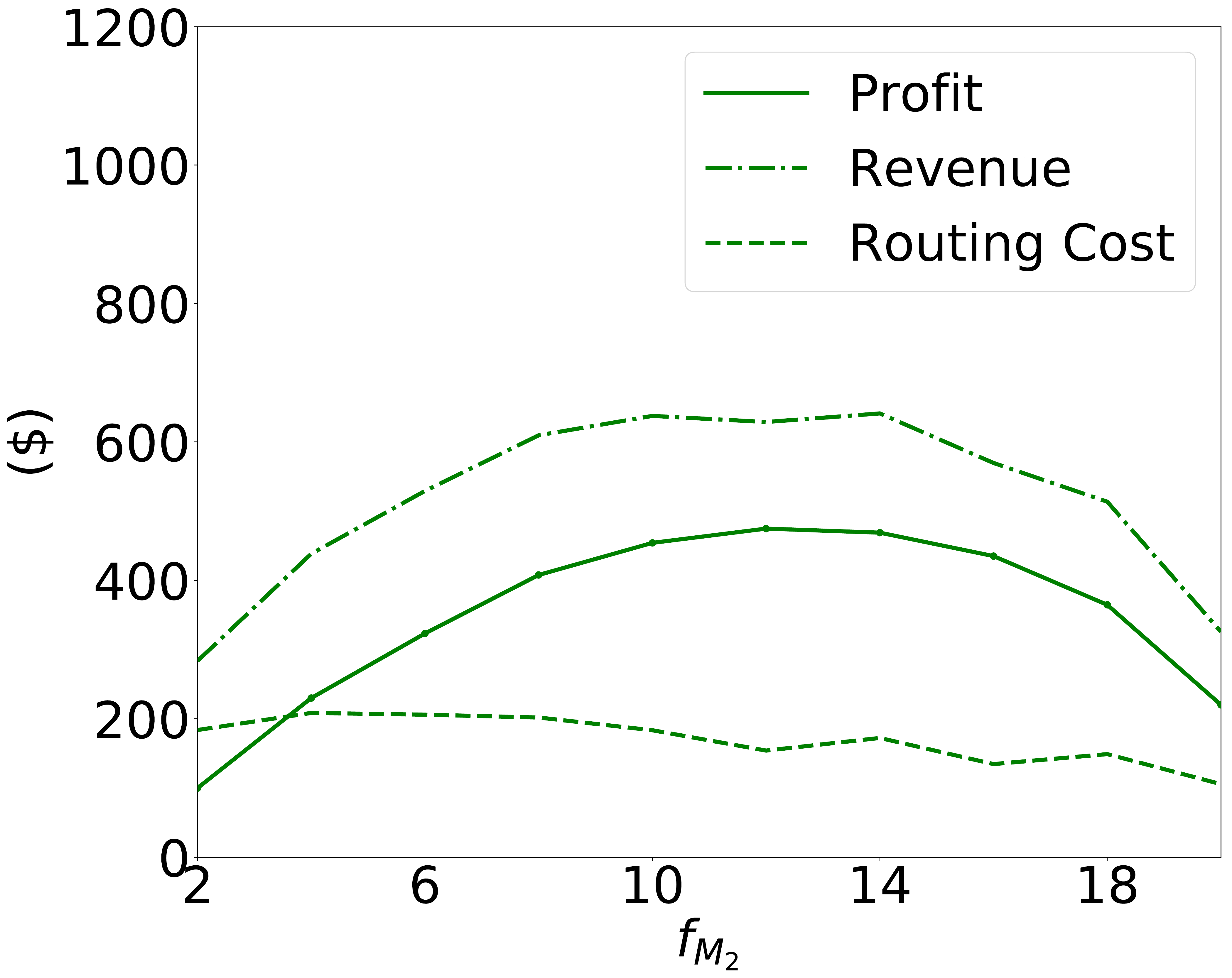}\caption{Profit, revenue and routing costs - $f_{m2}$.\label{zf2p}}}
\end{subfigure}\\
\begin{subfigure}{0.9\linewidth}{\includegraphics[width=0.95\textwidth,valign=c]{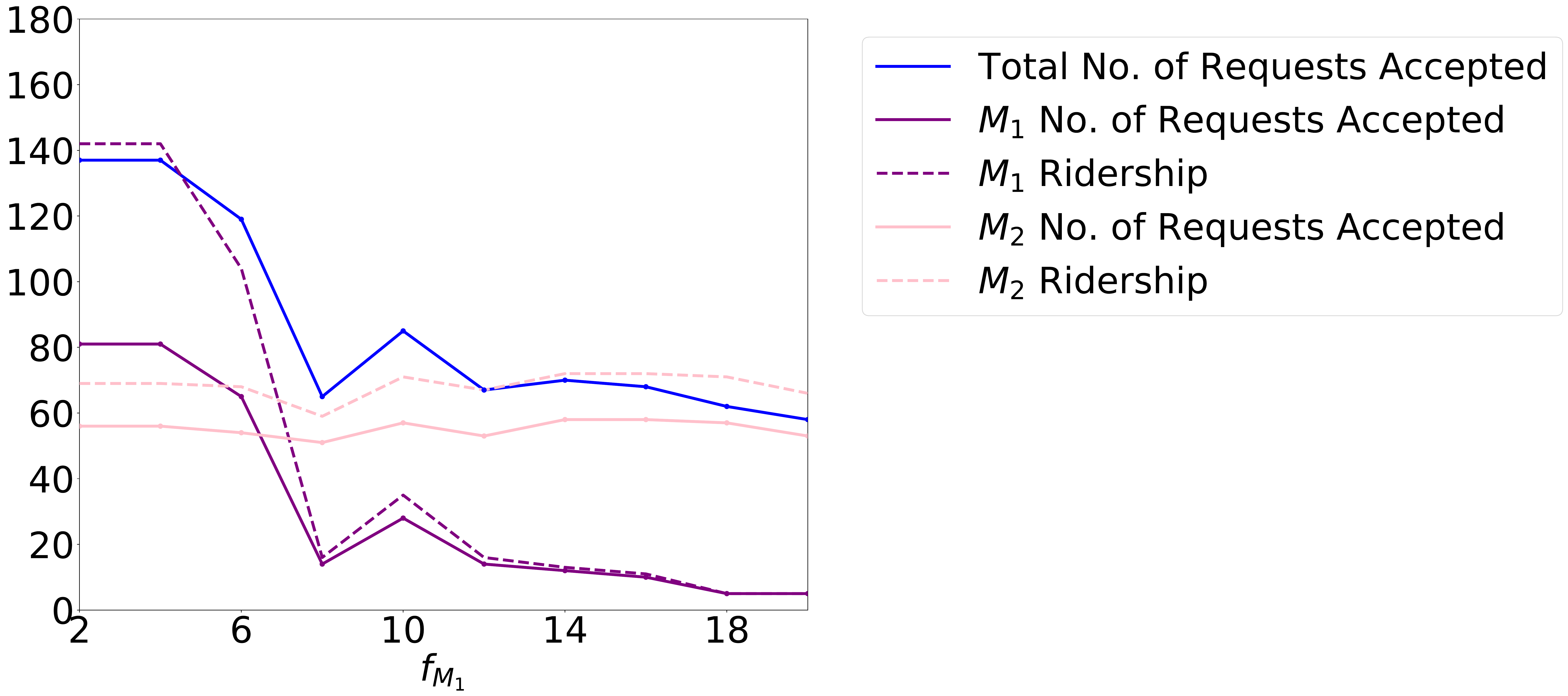}\caption{Number of accepted requests - $f_{m1}$.\label{zf1r}}}
\end{subfigure}
\begin{subfigure}{0.9\linewidth}{\includegraphics[width=0.95\textwidth,valign=c]{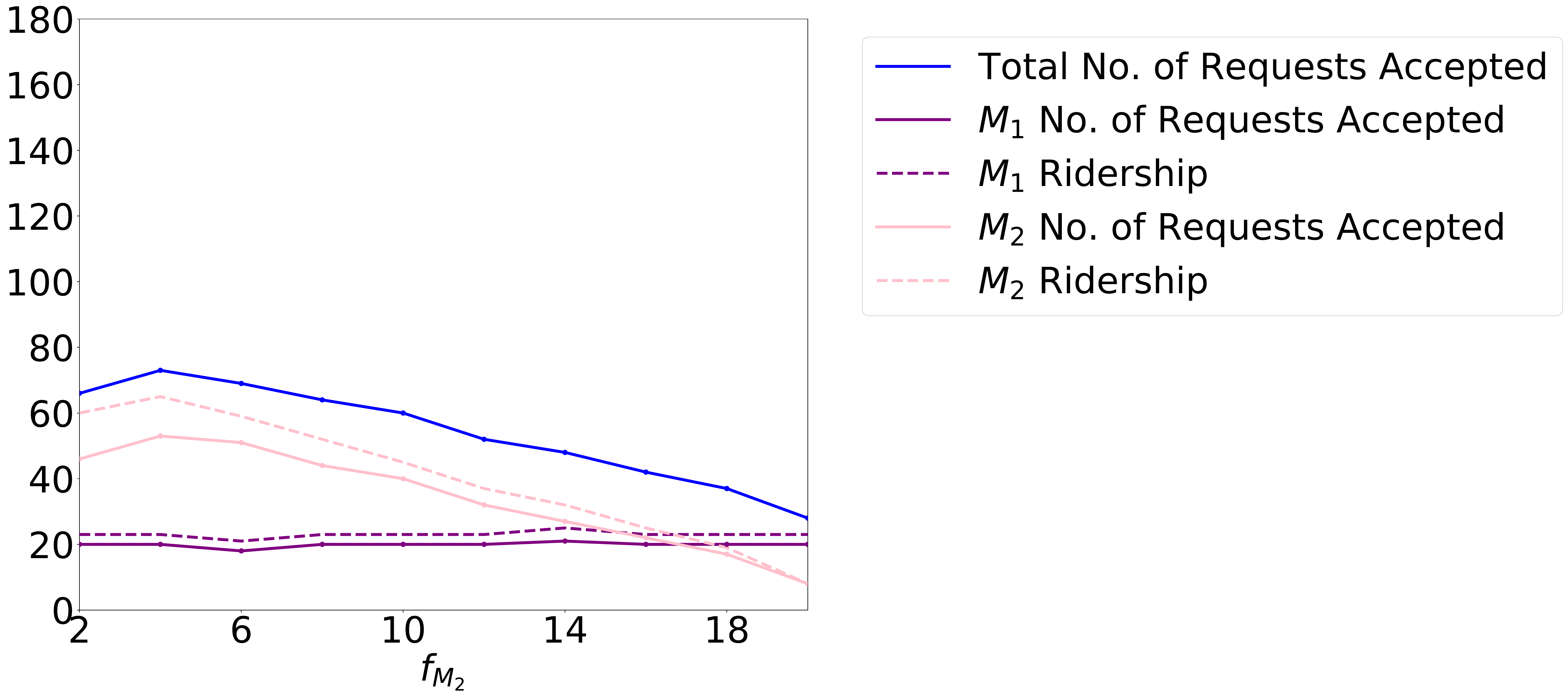}\caption{Number of accepted requests - $f_{m2}$.\label{zf2r}}}
\end{subfigure}\\
\caption{Results from the NYC instance when changing $f_{m1}$ and $f_{m2}$ under a zone-based fare structure.}
\label{zf}
\end{figure}

Similarly, we start with testing how increasing either $f_{M_1}$ or $f_{M_2}$ (while the other one is set to \$6.0) affects the overall profit and ridership. In Figures \ref{zf1p} and \ref{zf1r}, we find that the increasing $f_{M_1}$ causes more drastic fluctuations in both profit and ridership, since passengers of Class $M_1$ are more sensitive to changes in fare given their high VoT. Meanwhile, Figures \ref{zf2p} and \ref{zf2r} depict smoother transitions in curves when $f_{M_2}$ increases.

\begin{figure}[H]
\centering
\begin{subfigure}{0.48\linewidth}{\includegraphics[width=0.95\textwidth,valign=c]{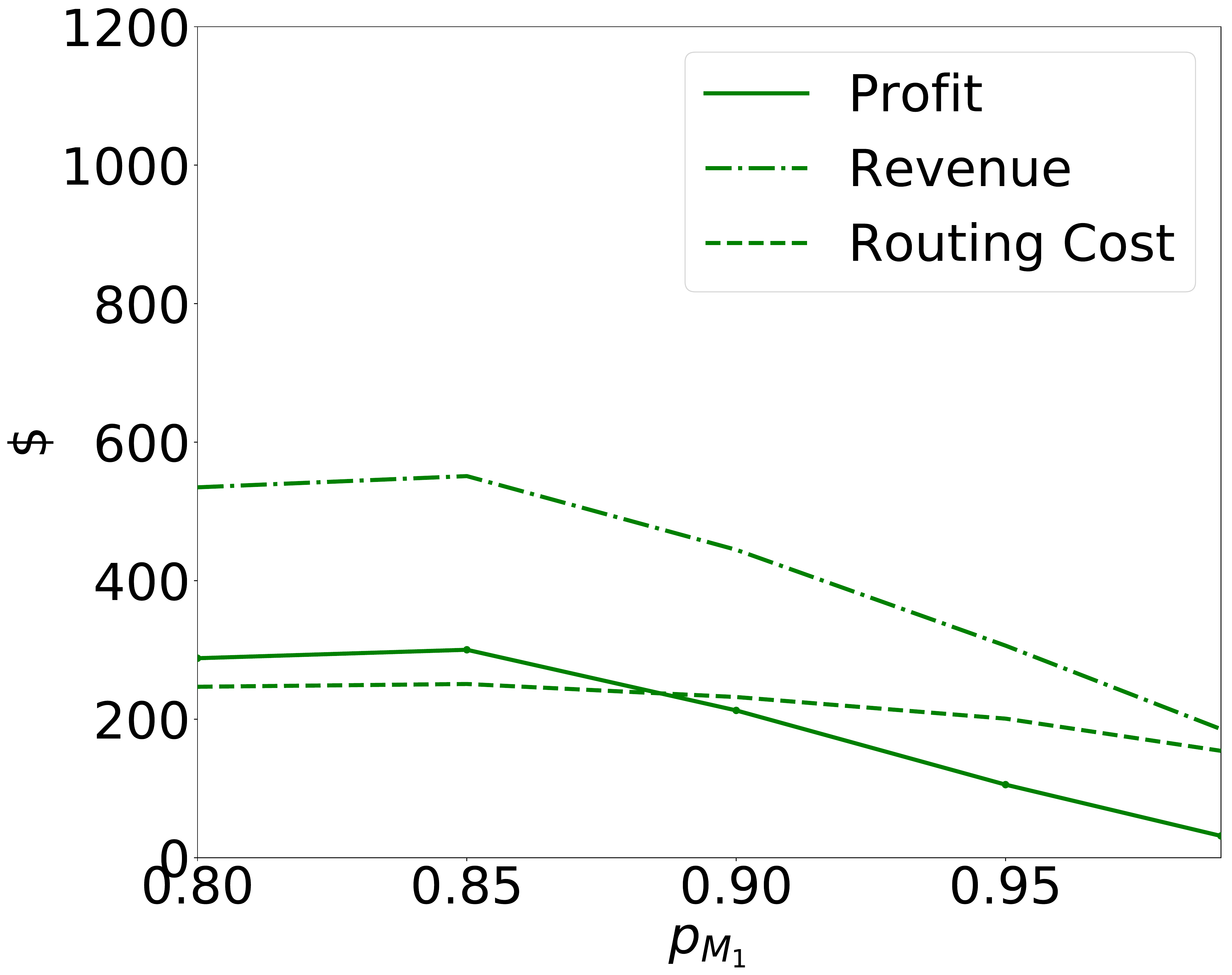}\caption{Profit, revenue and routing costs - $p_{m1}$.\label{zp1p}}}
\end{subfigure}
\begin{subfigure}{0.48\linewidth}{\includegraphics[width=0.95\textwidth,valign=c]{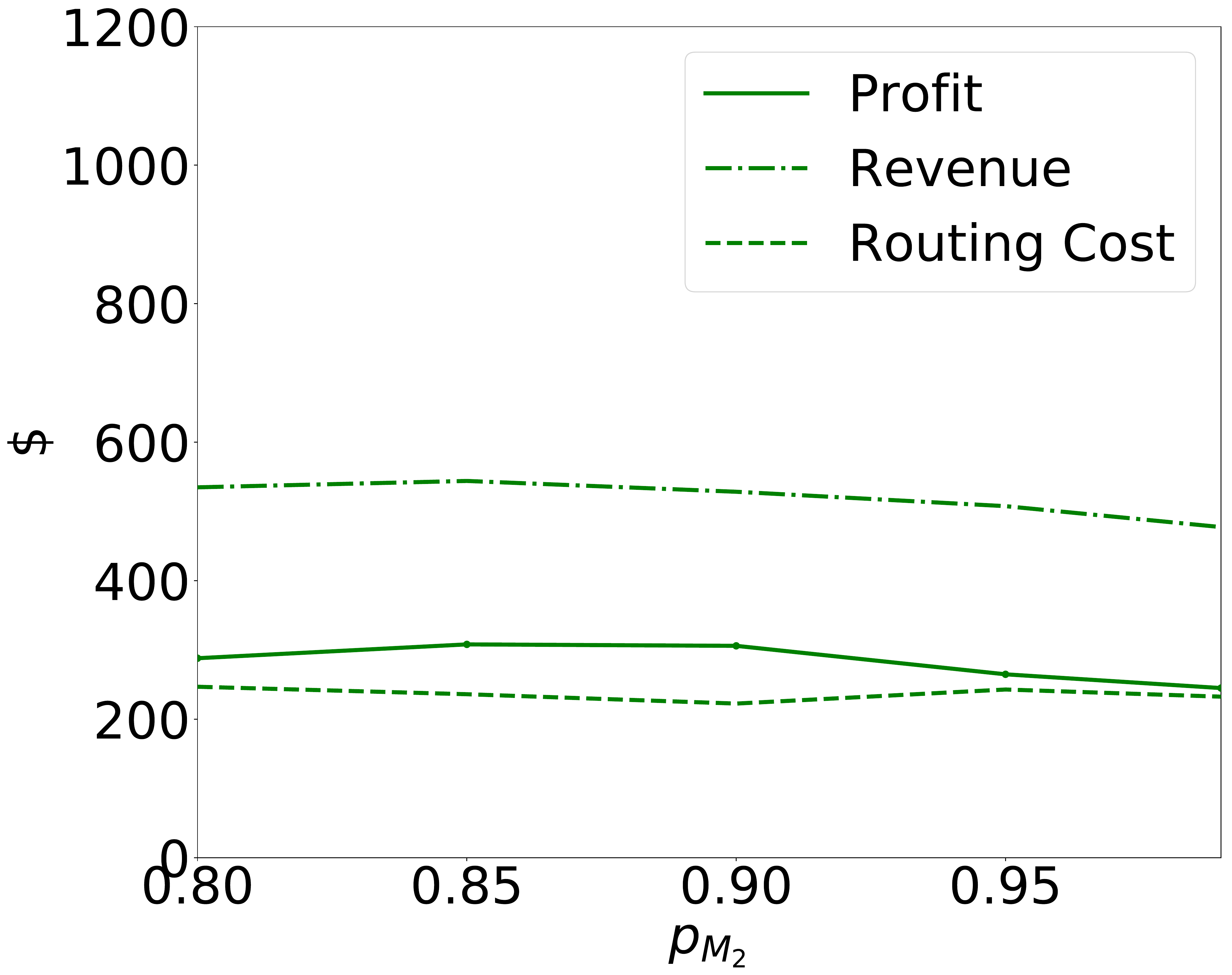}\caption{Profit, revenue and routing costs - $p_{m2}$.\label{zp2p}}}
\end{subfigure}\\
\begin{subfigure}{0.9\linewidth}{\includegraphics[width=0.95\textwidth,valign=c]{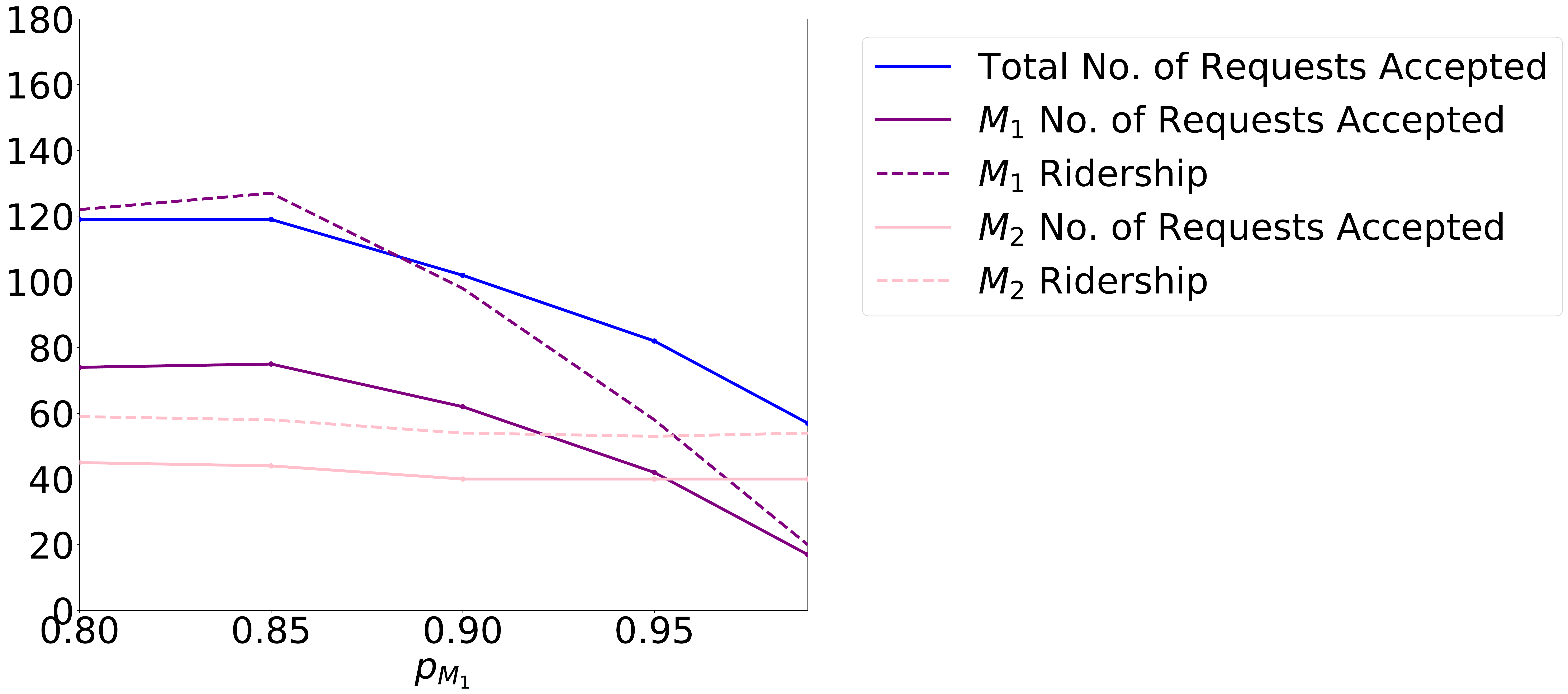}\caption{Number of accepted requests - $p_{m1}$.\label{zp1r}}}
\end{subfigure}
\begin{subfigure}{0.9\linewidth}{\includegraphics[width=0.95\textwidth,valign=c]{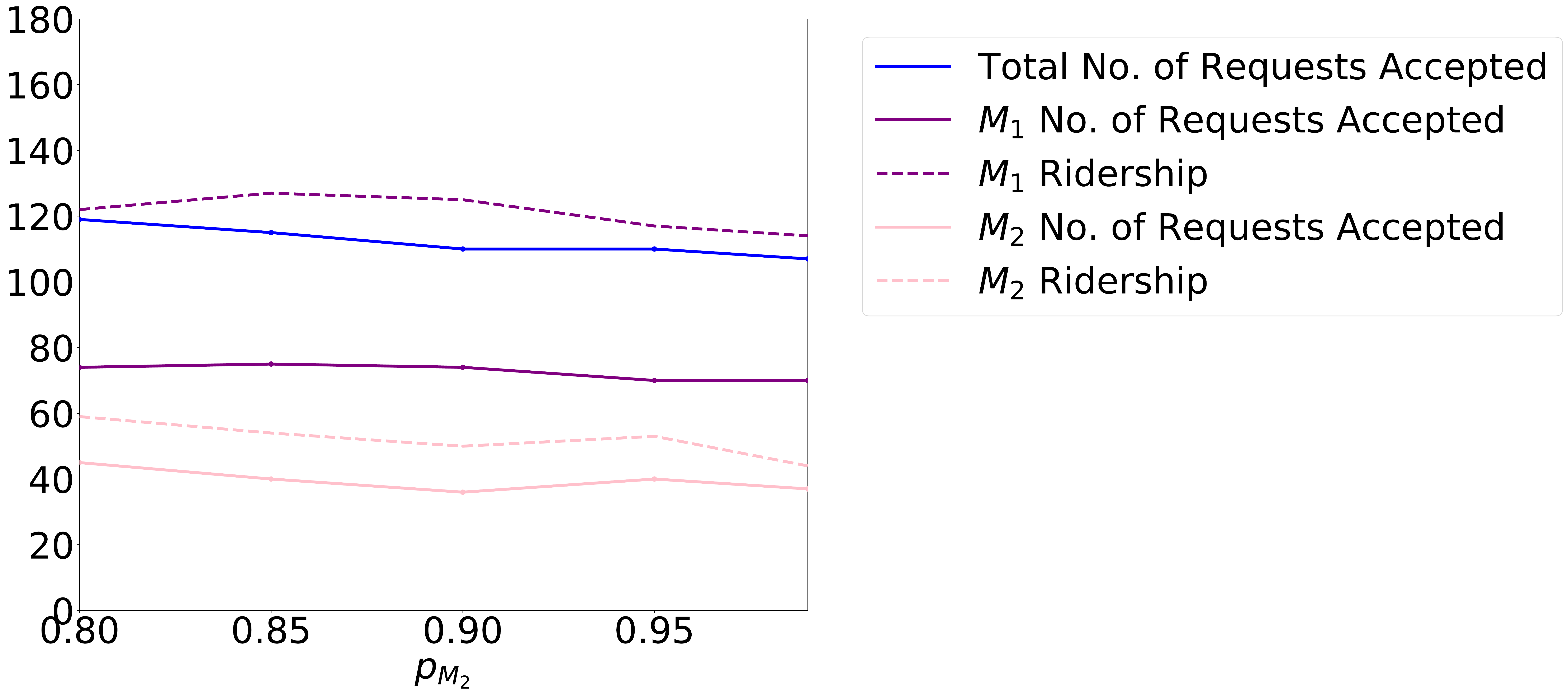}\caption{Number of accepted requests - $p_{m2}$.\label{zp2r}}}
\end{subfigure}\\
\caption{Results from the NYC instance when changing $p_{m1}$ and $p_{m2}$ under a zone-based fare structure.}
\label{zp}
\end{figure}

\begin{figure}[H]
\centering
\begin{subfigure}{0.9\linewidth}{\includegraphics[width=0.95\textwidth,valign=c]{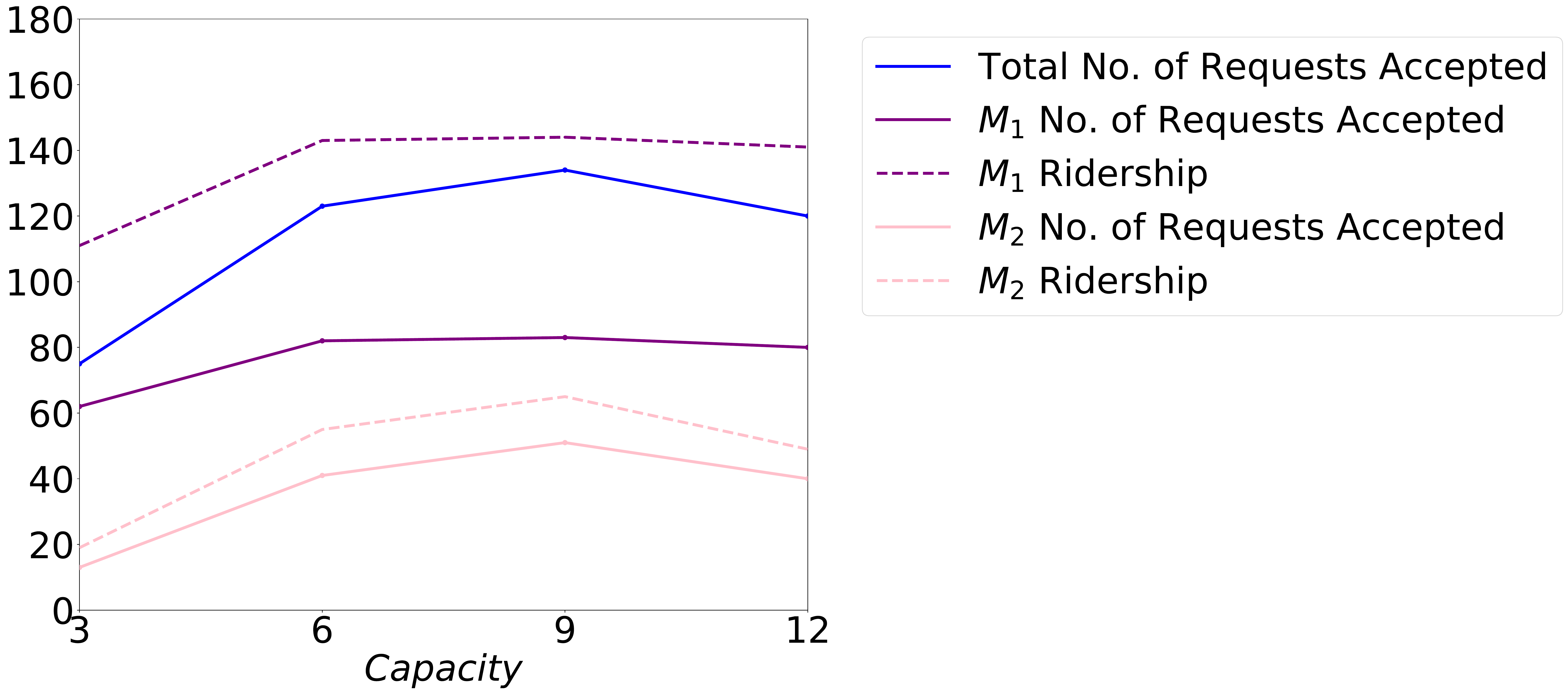}\caption{Number of accepted requests - Capacity.\label{zcr}}}
\end{subfigure}
\begin{subfigure}{0.5\linewidth}{\includegraphics[width=0.95\textwidth,valign=c]{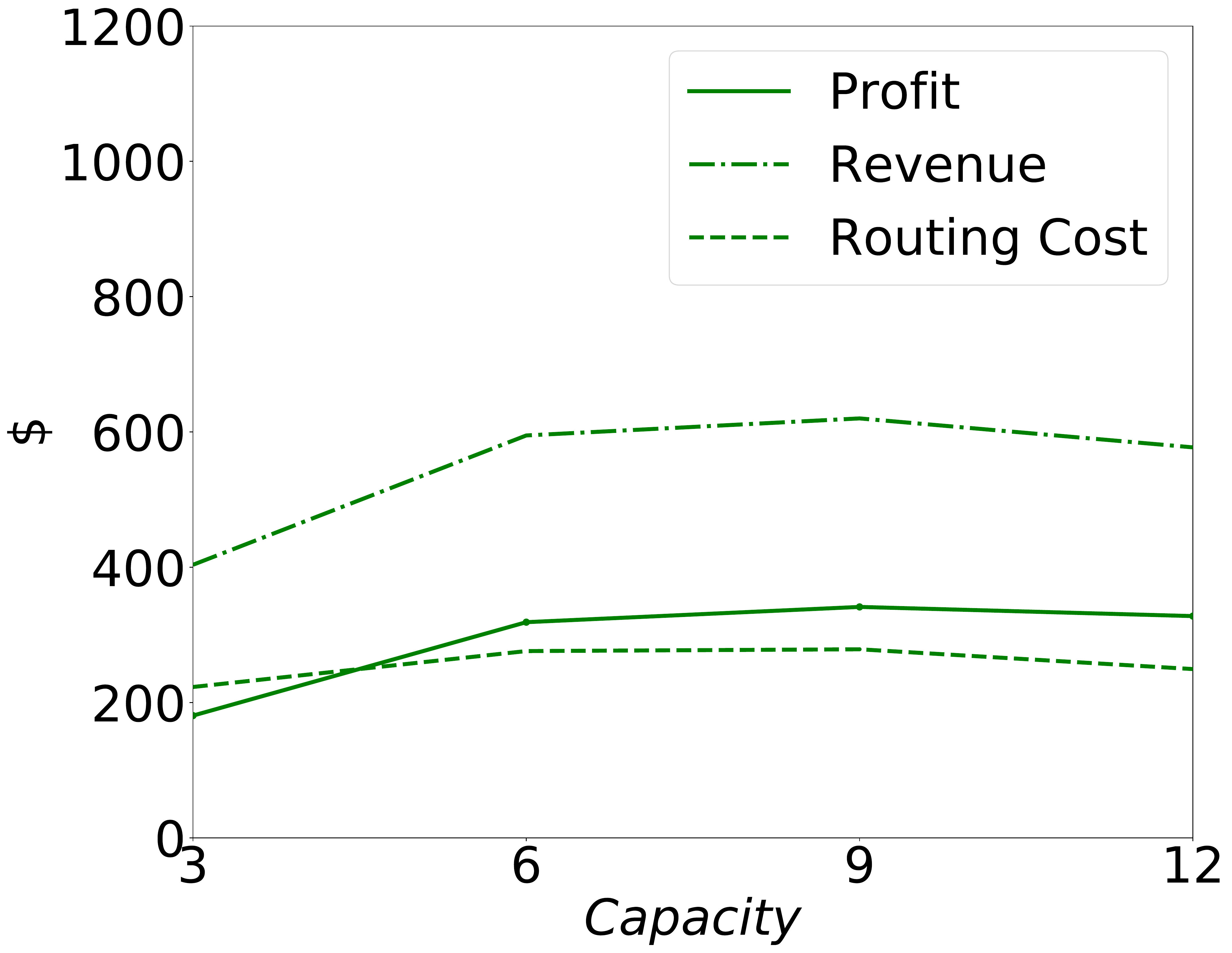}\caption{Profit, revenue and routing costs - Capacity.\label{zcp}}}
\end{subfigure}
\begin{subfigure}{0.45\linewidth}{\includegraphics[width=0.95\textwidth,valign=c]{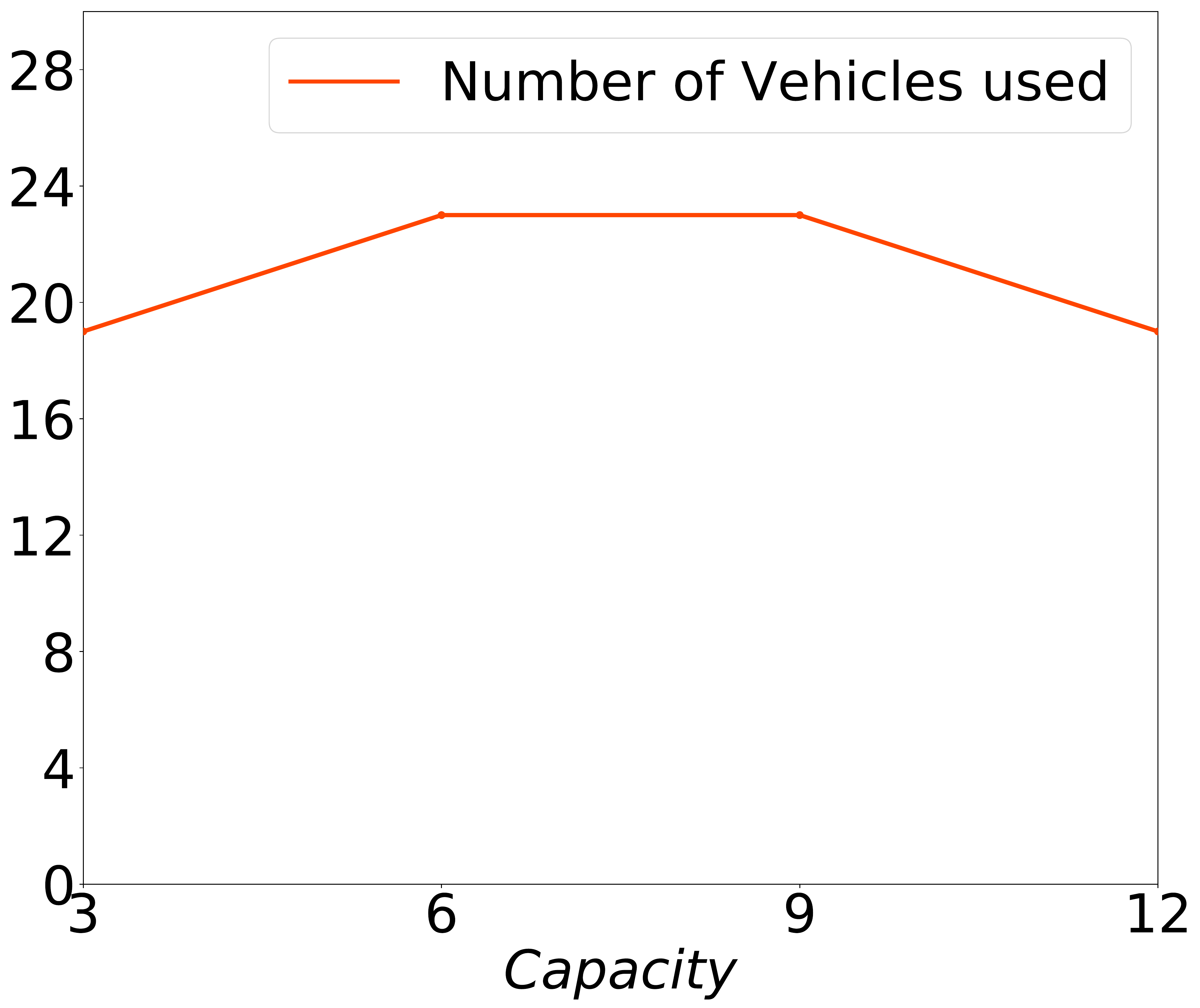}\caption{Fleet size - Capacity.\label{zcn}}}
\end{subfigure}
\caption{Results from the NYC instance when changing $Capacity$ under a zone-based fare structure.}
\label{zc}
\end{figure}

In both Figures \ref{zp} and \ref{zc}, $f_{M_1}$ is set to \$3.0 and $f_{M_2}$ to \$2.0. We then explore the impact of using a cheaper fare in Zone 2, how does tightening up chance constraints in each class impact profit and ridership (see in Figure \ref{zp}). With the combination of higher fare $f_{M_1}$ and higher VoT for Class $M_1$, the overall profit and ridership in $M_1$ drop sharply to almost zero in Figures \ref{zp1p} and \ref{zp1r} as $p_{M_1}$ increases. On the contrary, tightening the chance constraints in $M_2$ has minor impacts on the system, since $f_{M_2}$ is rather cheap and users are mostly better off when not driving.\\

Figures \ref{dc} and \ref{zc} illustrate how different vehicle capacities affect profit, ridership and fleet size under distance-based and zone-based fare structures, respectively. Note that even though the impact of capacity appears to be minor in these experiments when capacity is bigger than 6, it is not necessarily always the case. Growing capacity could have a more significant impact on routing optimisation and user selection when requests are more clustered spatially and/or temporally.

Overall, we find that the accept/reject mechanism of the proposed CC-DARP model reacts more sensitively towards the user class with a higher VoT ($M_1$), indicating that these users are more elastic when fare policy changes and the profitability of serving such users decreases more rapidly. Across three fare structures, we observe that the zone-based fare structure yields the largest profit and ridership, which gives DRT operators insight on revenue management strategies and pricing policies in the context of a class-based system. Furthermore, integrating the demographic attributes of each zone into our class-based model formulation gives zone-based fare structure the edge of designing customised pricing strategies. The flat fare structure, albeit simple, performs quite well. Additionally, its simplicity of form and low ticket issuing cost makes it more appealing to transit systems. The distance-based fare structure, on the other hand, performs poorly in maximising profit and ridership and its reactions towards variations in conditions are rather mild. Even though the distance-based fare structure is generally perceived as the fairest pricing policy, our findings suggests that it may not be a profit-maximizing strategy for DRT operators.

\section{Conclusion}
\label{con}
In this study, we introduced a novel chance-constrained DARP (CC-DARP) model that captures users' preferences in a DRT context in the long run. We model users' representative utilities for alternative transportation modes, including a Demand-Responsive Transportation (DRT) service. We consider utility-maximizing users and propose a mixed-integer programming formulation for the CC-DARP that captures users' preferences in the long run via a Logit model. Our formulation is multi-class which enables the model to explore the behaviours of various user classes under different price structures. This aims to provide insights to DRT operators on pricing policies design as well as revenue management strategies. We develop a customized heuristic solution method with layered local search structure to optimise both routing decisions and user selections. Numerical results showed that our customised local search based heuristic algorithm (LS-H) performed well on DARP benchmarking instances of the literature when compared to the best integer solution returned by an exact MILP solution. We further implemented the proposed solution method on a realistic case study derived from yellow taxi trip data from New York City (NYC). The results obtained on the realistic case study reveal that a zonal fare structure is a profit-maximizing strategy for DRT operators. 

The proposed CC-DARP formulation provides a new decision-support tool to inform on revenue and fleet management for DRT systems at a strategic planning level. Future work will aim to better integrate the proposed CC-DARP formulation with the rest of the transporation network, including fixed schedule public transport and competing DRT agents. More research is needed to develop an efficient exact solution algorithm for the CC-DARP and improve the robustness of heuristics for low fare case studies which have proven particularly challenging. Further, approaches to learn users' preferences could also provide rich areas of extensions of this research.

\clearpage

\bibliographystyle{elsarticle-harv}
\bibliography{references}
	
\end{document}